# The classification of punctured-torus groups

By Yair N. Minsky*

**Abstract**

Thurston's ending lamination conjecture proposes that a finitely generated Kleinian group is uniquely determined (up to isometry) by the topology of its quotient and a list of invariants that describe the asymptotic geometry of its ends. We present a proof of this conjecture for punctured-torus groups. These are free two-generator Kleinian groups with parabolic commutator, which should be thought of as representations of the fundamental group of a punctured torus.

As a consequence we verify the conjectural topological description of the deformation space of punctured-torus groups (including Bers' conjecture that the quasi-Fuchsian groups are dense in this space) and prove a rigidity theorem: two punctured-torus groups are quasi-conformally conjugate if and only if they are topologically conjugate.

## Contents

*This work was partially supported by an NSF postdoctoral fellowship and a fellowship from the Alfred P. Sloan Foundation.

# 1. The ending lamination conjecture and its consequences

The general classification problem for discrete groups of Möbius transformations remains tantalizingly open, although a conjectural picture of the solution has been in place since the late 70's, and is roughly as follows. In the representation space for a given group $G$ into $\mathrm{PSL}_2(\mathbf{C})$, the discrete, faithful elements are expected (barring trivial cases) to comprise the closure of an open set of *structurally stable* representations. In a component of the structurally stable set all representations are quasi-conformally conjugate and hence parametrized by a Teichmüller space, or a quotient of one. On the boundary of this set one obtains *geometrically infinite* groups and groups with new parabolics, and these are expected to be parametrized by what remains of the Teichmüller parameter, together with a combinatorial invariant known as an *ending lamination* (see Abikoff [1] for an overview).

In this paper we verify this conjectural picture for punctured-torus groups, which are the simplest of all classes of Kleinian groups with a nontrivial deformation theory. The primary component of the solution is the proof of Thurston's "ending lamination conjecture" in this case (Theorem A).

A punctured-torus group is a free, discrete, two-generator group $\Gamma$ of (orientation-preserving) Möbius transformations with the added condition that the commutator of the generators is parabolic. We should think of $\Gamma$ as the image of a representation $\rho : \pi_1(S) \to \mathrm{PSL}_2(\mathbf{C})$, where $S$ is a once-punctured torus (to keep the representation in mind we often call this a *marked* group). The commutator condition means that the loop surrounding the puncture determines a cusp of the three-manifold $\mathbf{H}^3/\Gamma$, and in general a representation of a surface group taking cusps to cusps in this way is called *type-preserving*.

To such a representation one may associate an ordered pair of *end invariants* $(\nu_-, \nu_+)$, lying in $(\overline{\mathbf{D}} \times \overline{\mathbf{D}}) \setminus \Delta$, where $\overline{\mathbf{D}}$ is the closed unit disk whose interior $\mathbf{D}$ is identified with the Teichmüller space of $S$, and whose boundary $S^1$ is identified with the space of measured laminations on $S$. We denote by $\Delta$ the diagonal of $S^1 \times S^1$. We also identify $S^1$ with $\hat{\mathbf{R}} = \mathbf{R} \cup \{\infty\}$ by stereographic projection, and let $\hat{\mathbf{Q}} = \mathbf{Q} \cup \{\infty\}$. When both $\nu_\pm$ lie in $\mathbf{D} \cup \hat{\mathbf{Q}}$ the group is *geometrically finite*; this case has been well-understood through work of Ahlfors, Bers, Kra, Marden, Maskit and others. When either invariant lies in $\mathbf{R} \setminus \mathbf{Q}$ the group is geometrically infinite (and the invariant is called a "lamination on $S$"), and the existence of the invariant in this case is due to Thurston and Bonahon. Thurston's ending lamination conjecture states that these invariants suffice to determine the group up to isometry. See Section 3 for more precise definitions, and [70], [71] for discussions of the conjecture for more general groups. In this paper we shall prove:

THEOREM A (Ending Lamination Theorem). *A marked punctured-torus group $\rho : \pi_1(S) \to \mathrm{PSL}_2(\mathbf{C})$ is determined by its end invariants $(\nu_-, \nu_+)$, uniquely up to conjugacy in $\mathrm{PSL}_2(\mathbf{C})$.*

In other words, the map $\nu : \rho \mapsto (\nu_-, \nu_+)$ is injective. It can also be shown to be surjective, as a consequence of Bers' Simultaneous Uniformization Theorem [10], and of Thurston's Double Limit Theorem [85]. We will further obtain the following theorem about the deformation space $\mathcal{D}(\pi_1(S))$ of all punctured-torus Kleinian groups, modulo conjugation in $\mathrm{PSL}_2(\mathbf{C})$:

THEOREM B (deformation space topology). *The map*
$$\nu^{-1} : (\overline{\mathbf{D}} \times \overline{\mathbf{D}}) \setminus \Delta \to \mathcal{D}(\pi_1(S))$$
*is a continuous bijection.*
*In addition, every Bers slice is a closed disk, and every Maskit slice is a closed disk with one boundary point removed.*

(See Section 12.3 for definitions of Bers and Maskit slices). Note that this does not imply $\nu^{-1}$ is a homeomorphism and in fact $\nu$ itself is discontinuous! (See Section 12.3, and Anderson-Canary [6].) However, as the interior $\mathbf{D} \times \mathbf{D}$ of the space of invariants maps precisely to the set of structurally stable (or, in particular, quasi-Fuchsian) representations, from a dynamical point of view we have proved that "structural stability is dense" for this family of representations. In particular this gives a positive answer (for punctured-torus groups) to Bers' conjecture in [11] that all degenerate groups in a Bers slice are limits of quasi-Fuchsian groups.

A final application, also with a dynamical flavor, is the following rigidity theorem:

THEOREM C (qc rigidity). *If the actions of two punctured-torus groups on the sphere are conjugate by a homeomorphism, then they are conjugate by a quasiconformal or anti-quasiconformal homeomorphism, according as the original homeomorphism preserves or reverses orientation.*

The core of the proof of Theorem A is the Pivot Theorem 4.1, which is the main step to getting quasi-isometric control of the group in terms of what amounts to the continued-fraction expansions of the end invariants. In particular, the presence of very short geodesics in the quotient manifold is predicted precisely by the presence of high coefficients in the expansion, as has long been conjectured. The statement is given (see §4) in terms of the combinatorics of the *Farey triangulation* in the disk (see §2).

1.1. *Historical comments.* Ahlfors and Bers (see [4], [10], [11]) first studied the deformation theory of quasi-Fuchsian groups (in any genus) and showed that they are parametrized by a product of Teichmüller spaces ($\mathbf{D} \times \mathbf{D}$ in our case). Maskit [60] further studied the groups that arise on the boundary of these deformation spaces when the domains of discontinuity are "pinched" and new parabolics arise – this corresponds in our discussion to the case when $\nu_+$ or $\nu_-$ is a rational point in the boundary $\hat{\mathbf{R}}$. Keen-Maskit-Series [54] gave a proof of Theorem A in this case. Jørgensen made some very careful studies of quasi-Fuchsian punctured-torus groups in [46], in particular obtaining a combinatorial description in terms of the Farey triangulation which is very closely related to the results we obtain in Theorem 4.1. He also studied some degenerate groups, in particular with Marden in [50], applying the triangulation to show that two particular degenerate groups are not quasiconformally conjugate. These ideas have been helpful to the writing of this paper. Degenerate groups were first shown to exist by Bers, and then in greater generality by Thurston, who analyzed them geometrically and introduced the ending lamination invariant (§3). Bonahon showed that Thurston's theory applied in fact to all Kleinian surface groups. Our analysis takes these developments as its starting point.

The problem has also been studied by McMullen [67] who showed that cusps (representations corresponding to a rational $\nu_+$) are dense in the boundary of a Bers slice. Wright [89] has carefully analyzed the combinatorics of limit sets for groups lying at the boundary of a punctured-torus Maskit slice, and produced some very good computer pictures of such limit sets and of the boundary itself, giving considerable evidence to support the above conjectures. Keen and Series [53] have given geometric coordinates for the interior of punctured-torus Maskit slices, in particular generalizing some of Wright's findings.

Recently, Bowditch [16] has given an analysis of trace functions on the Farey graph arising from general (not necessarily discrete) representations of $\pi_1(S)$ with parabolic commutator, using algebraic methods. In particular he has established a considerably stronger version of our Lemma 8.1, by a completely different proof. Alperin, Dicks and Porti [5] have used the Farey graph to analyze the geometry of the Gieseking manifold, which is a specific punctured-torus bundle over the circle. In particular they have given an alternate proof of the theorem of Cannon-Thurston [25] in this case.

1.2. *Summary of the proof.* To simplify this discussion let us consider here a punctured-torus group $\rho$ for which both $\nu_+$ and $\nu_-$ are irrational points in $\hat{\mathbf{R}}$. In other words, the manifold $N = \mathbf{H}^3/\rho(\pi_1(S))$ has two simply degenerate ends and no domain of discontinuity. This allows us to avoid a number of special cases in the argument having to do with boundaries of convex hulls or accidental parabolics. In fact, the reader is strongly advised to make this assumption on a first reading of the proof itself.

The problem of proving Theorem A reduces to showing that the end invariants $\nu_\pm$ describe the group, or manifold, up to *quasi-isometry* in $\mathbf{H}^3$. Then any two representations with the same invariants are conjugate by a quasi-conformal homeomorphism, and one can use Sullivan's rigidity theorem [83] to show that the conjugating map is, in fact, a Möbius transformation.

However, the end invariants only give asymptotic information: for example, simple closed curves in $S$ are represented by rational numbers in $\hat{\mathbf{R}}$ (§2.1), and $\nu_\pm$ are characterized by the property that any infinite sequence of curves whose corresponding geodesics in $N$ have uniformly bounded lengths gives rise to a sequence of rational numbers accumulating onto $\nu_\pm$. To know the quasi-isometry type of the manifold we need at the very least to determine which such sequences can arise.

*Farey graph and pivot sequence*: In Section 2 we discuss the Farey triangulation $\mathcal{C}$, a well-known triangulation of the disk with vertices in $\hat{\mathbf{Q}} \equiv \mathbf{Q} \cup \{\infty\}$, which can be interpreted in terms of slopes and intersection numbers of simple closed curves on $S$. The two irrational points $\nu_\pm$ determine, via the combinatorial structure of $\mathcal{C}$, a bi-infinite sequence $P = \{\alpha_n\} \subset \hat{\mathbf{Q}}$, closely related to continued fraction approximations, such that $\alpha_n \to \nu_\pm$ as $n \to \pm\infty$. We call these vertices *pivots*. A good starting point for reading the paper is Section 4, where we define $P$ and state the Pivot Theorem 4.1. This theorem asserts that the pivots indeed have bounded length in the manifold, and furthermore gives an explicit recipe for estimating their complex translation lengths from the combinatorial data of $P$.

*Connectivity*: The main idea that leads to the Pivot Theorem is an application of the fact that paths in $\mathcal{C}$ correspond to continuous families of simplicial hyperbolic surfaces in the 3-manifold. With this we prove Lemma 8.1, which states that the set of vertices in $\mathcal{C}$ whose geodesics in the 3-manifold satisfy a certain length bound is connected. From here it is easy to obtain an *a priori* bound on the lengths of all of the pivots, which we do in Lemma 8.2.

*Bounded homotopies and control on Margulis tubes*: The second idea is roughly this: two homotopic Lipschitz maps of the same nonelementary hyperbolic surface into a hyperbolic 3-manifold are connected by a homotopy of bounded length (the bound depending on the Lipschitz constant). We prove a version of this via the "figure-8 argument" in Section 9.5, and this allows us to constrain the geometry of Margulis tubes in the manifold. In particular, any Margulis tube in $N$ can be encased, homologically, by a pair of surfaces with controlled geometry and a homotopy between them which has bounded tracks in the complement of the Margulis tube. This discussion is carried out in Section 9 via a mechanism we call a "building block". In particular each block $\mathcal{B}$ contains a solid torus $U$ whose boundary torus $\partial U$ is mapped to the boundary of the corresponding Margulis tube. It follows, for example, that there is a uniform bound on the diameter of the boundary torus of the Margulis tube.

*Halfway surfaces*: The block construction also depends on an analysis we carry out in Lemma 7.1 of Section 7, to describe the possible geometric configurations of axes for generator pairs in punctured-torus groups. In particular when a pair of generators have bounded lengths we can find a simplicial hyperbolic map of $S$ into $N$ in which both generators are *simultaneously* bounded. We call these halfway surfaces, and they are used to begin the building block construction in Section 9. The proof of Lemma 7.1 is carried out by fairly standard use of trace identities, although one can give a more geometric argument (as was done in a previous draft of this paper).

These ingredients are put together in Section 10, where the proof of the Pivot Theorem is completed. The blocks, one per pivot, are glued end to end to produce a "model manifold" $M = \bigcup_n \mathcal{B}_n$ and a Lipschitz map $f : M \to N$, that in particular takes the solid torus $U_n$ in each block to the Margulis tube it is meant to model. A subtle point to emphasize here is that, for each solid torus individually it is not *a priori* clear that the map gives a faithful model; for example, the map restricted to each boundary torus $\partial U_n$ is not at first known to be homotopic to a homeomorphism. This issue is settled by a *global* argument showing that the map $f$ is proper and has degree 1 (Lemma 10.1).

Once this map is in place we have enough control to estimate a Teichmüller parameter for the boundary torus of each Margulis tube, with respect to a natural marking of the torus. As described in Section 6.2, this is exactly what we need to determine the quasi-isometry type of each Margulis tube, and give the statement of Theorem 4.1.

The construction of the model manifold is actually completed in Section 11, where we must face the fact that our map $f : M \to N$ is only Lipschitz, and not bilipschitz. What we actually need is that $f$ lifts to a quasi-isometry of the universal covers, and this is the goal of Theorem 11.1. To prove this theorem, we must switch to a different mode: instead of explicit constructions and plausibly computable bounds, we appeal to compactness arguments. In particular we must consider the possible geometric limits of sequences of punctured-torus manifolds with basepoints. These limits conform to the already well-known picture of "drilled holes" developed by Thurston [85] and Bonahon-Otal [15], in which an infinite sequence of new rank-2 cusps can appear in the limiting manifold. (Note that for higher-genus surface groups much more dramatic limits can occur – see Brock [18].) The geometric limits we obtain come equipped with their own model manifolds, and this control in the limit suffices to give us the uniform bounds we need.

The proofs of the main theorems A, B and C are carried out in Section 12, and are fairly straightforward given what has gone before. A subtle issue which arises in the proof of Theorem B involves the continuity of the end invariants under algebraic limit, which in fact does not hold in general. We are led back to consideration of the geometric limit apparatus of Section 11 in order to resolve this point.

Sections 2, 3, 5 and 6 are essentially expository. In Section 2 we discuss the Farey graph and some elementary cartography of the Teichmüller space of the torus, and give a lemma on quasiconformal maps. In Section 3 we state the definitions of the end invariants $\nu_\pm$. In Section 5 we discuss simplicial hyperbolic surfaces and pleated surfaces, which will play a central role in almost all of our arguments. In Section 6 we discuss Margulis tubes. In particular we state some well-known bounds on the radii of such tubes, due to Brooks-Matelski and Meyerhoff, and develop in Section 6.2 the connection between the geometry of a Margulis tube and a parameter in Teichmüller space describing its quotient torus at infinity. In Section 6.3 we give some additional constraints, due to Thurston and Bonahon, on Margulis tubes which appear in surface groups.

1.3. *Speculations on the general case.* There are a few straightforward generalizations of the ideas in this paper. The quadruply-punctured sphere can be treated almost identically: in particular the combinatorics of the set of simple closed curves are again encoded by the Farey triangulation, and the figure-8 argument in Lemma 9.3 still applies.

Furthermore, let $(M, P)$ be a "pared" manifold, that is let $M$ be a compact 3-manifold and $P$ a collection of tori and annuli on $\partial M$, and suppose that all the components of $\partial M \setminus P$ are punctured tori or 4-punctured spheres which are incompressible in $M$. Suppose that $M$ admits an embedding into a hyperbolic 3-manifold $N$ which is a homotopy-equivalence and takes each component of $P$ into a distinct parabolic cusp. Then the techniques of this paper apply directly to the ends of the resulting manifold, and a suitably restricted version of the ending lamination theorem holds (see [70] to see how this type of argument works). One may also extend the rigidity theorem (C) to this context; see e.g. Ohshika [75].

Beyond this, one must begin to consider the general problem for higher-genus surface groups. A number of very serious difficulties arise here. There is a natural simplicial complex generalizing the Farey graph, but its properties are considerably harder to understand. In Masur-Minsky [64] we study this complex from a point of view partially motivated by these ideas. The bound on the diameter of Margulis tube boundaries also fails in general, and this is related to the fact that geometric limits of sequences of general surface groups can be much more complicated than what we obtain in Section 11.2. See Brock [18] for some of the types of phenomena that can occur.

*Acknowledgements.* I am very grateful to Dick Canary, Curt McMullen and Jeff Brock, for many enjoyable and illuminating conversations on the subject of this paper, and to Caroline Series who pointed out to me the special nature of the punctured-torus case. Special thanks are due to Ada Fenick, who told me I had to finish writing it.

## 2. The Farey triangulation and the torus

Let $\mathbf{H}^2$ denote the upper half plane with boundary $\mathbf{R}$. There is a classical ideal triangulation of $\mathbf{H}^2$, defined as follows. For any two rational numbers written in lowest terms as $p/q$ and $r/s$, say they are *neighbors* if $|ps - qr| = 1$. Allow also the case $\infty = 1/0$. Joining any two neighbors by a hyperbolic geodesic, we obtain the Farey triangulation. (The proof that this is a triangulation is easy after we consider the edges incident to $\infty$, and observe that the diagram is invariant under the natural action of $\mathrm{SL}_2(\mathbf{Z})$. See also Series [82], [81] or Bowditch [16]). Figure 1 shows this triangulation, after stereographic projection to the unit disk $\mathbf{D}$. This picture is intimately related to the torus, as we shall now see.

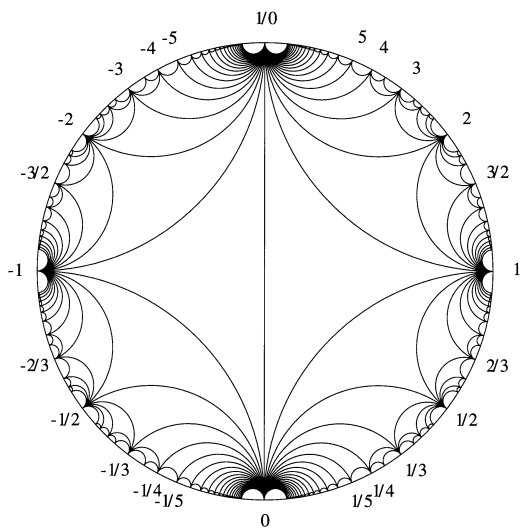

Figure 1. The Farey graph in the unit disk, with some of the vertices labeled.

2.1. *The complex of curves.* Let $\mathcal{C}$ denote the set of free homotopy classes of unoriented simple nonperipheral curves on the punctured torus $S$. These are in one-to-one correpondence with $\hat{\mathbf{Q}} \equiv \mathbf{Q} \cup \{\infty\}$, after one chooses an ordered basis for $H_1(S)$. Let us do this explicitly in order to be careful about sign conventions:

Fix an orientation for $S$ and choose a positively oriented ordered basis $(\alpha, \beta)$ for $H_1(S)$. This is equivalent to choosing two oriented simple closed curves which cut $S$ into a positively oriented rectangle. Any element of $H_1(S)$ can be written as $(p, q) = p\alpha + q\beta$ in this basis, and we associate to this the slope $-p/q \in \hat{\mathbf{Q}}$. Note that this ratio forgets the orientation of the curve, as well as integer multiples. Thus it exactly describes an element of $\mathcal{C}(S)$. (The proof that every element of $\mathcal{C}$ is obtained in this way is left to the reader.) The determinant $ps - rq$ which appeared above is easily seen to be just the oriented intersection number $i(\cdot, \cdot)$ in $S$. Let $\alpha \cdot \beta = |i(\alpha, \beta)|$ denote the unoriented intersection number, which is defined on $\mathcal{C}$.

Thus the Farey graph reflects the combinatorial structure of $\mathcal{C}$, and from now on we shall identify the two.

*Remarks.* 1. $\mathcal{C}$ is a special case of the Hatcher-Thurston complex [41], and is related closely to the complexes of curves introduced by Harvey [40] and studied by Harer [38], [39] and Ivanov [43], [42], [44]. See also Bowditch-Epstein [17] for another perspective. 2. A pair of vertices joined by an edge can also be considered as representing a pair of generators for $\pi_1(S)$, up to conjugation and inverses, as in Jørgensen [46] and Jørgensen-Marden [50]. (See also Section 7.) 3. The same construction works for the regular torus – the difference there is that we do not need to worry about peripheral curves or nonsimple curves that, say, wind around the puncture.

2.2. *Neighbors and Dehn twists.* Given $\alpha \in \mathcal{C}$, its neighbors may be indexed by the integers $\{\beta_n\}_{n \in \mathbf{Z}}$ according to their counterclockwise order around $S^1 \setminus \{\alpha\}$. Denote by $D_\alpha$ the positive Dehn twist around $\alpha$, defined for example by the convention that positive twists around a vertical curve increase slope. Note that "positive" makes sense after a choice of orientation on $S$, but without orienting $\alpha$ – see Figure 2. Then the indices are defined, after arbitrary choice of $\beta_0$, by $\beta_n = D_\alpha^n(\beta_0)$.

2.3. *Teichmüller space.* The interior $\mathbf{D}$ of the disk, or half-plane $\mathbf{H}^2$, also has a well-known interpretation in terms of the torus: it parametrizes the Teichmüller space $\mathcal{T}(S)$ of conformal, or hyperbolic, structures on $S$. (Recall that the conformal structures on the regular and once-punctured torus are the same – although the regular torus admits no hyperbolic metric). In this interpretation, the circle $\hat{\mathbf{R}} \equiv \mathbf{R} \cup \{\infty\}$ is Thurston's compactification of $\mathcal{T}(S)$ using projective measured laminations – where the rational points correspond to simple closed curves, as above, and the irrational points to laminations with infinite leaves (see [31], [79]).

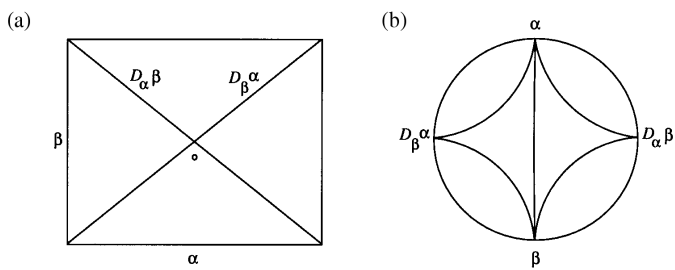

Figure 2. Positive Dehn twists (a) in the punctured torus, and (b) in the Farey triangulation.

Explicitly, fix a marking $(\alpha, \beta)$ of the torus, as above. To a point $z$ in the upper half plane we associate the lattice generated by 1 and $z$, whose quotient is a torus with induced conformal structure $\nu$. The position of the puncture is irrelevant since the torus has a transitive family of conformal automorphisms. An orientation-preserving identification of our fixed torus $S$ with this torus is determined by taking the curve $\alpha$ to the image of $[0, 1]$ and $\beta$ to the image of $[0, z]$.

For later convenience let us denote by

$$z = \tau(S, \nu, \alpha, \beta)$$

the relationship between a marked conformal torus and its Teichmüller parameter. (We also write $\tau(S, \alpha, \beta)$ if the conformal structure on $S$ is understood.)

Let us also recall the *Teichmüller metric* which is defined as $\frac{1}{2} \log K$ where $K$ is the best dilatation constant for a marking-preserving quasiconformal homeomorphism between two marked tori. We will use the fact that this is exactly equal to the hyperbolic metric on $\mathbf{H}^2$.

*Shortest curves.* Extremal lengths (see Ahlfors [3]) can be computed directly in the Euclidean metric inherited by the torus, as length squared over area. Note in particular that $\alpha$, which was identified with $\infty \in \hat{\mathbf{R}}$, has extremal length $1/\mathrm{Im}\, z$ in the structure parametrized by $z$. When this quantity is very small it gives a good estimate for the hyperbolic length of $\alpha$ (see Maskit [62]).

To get a clean picture for shortest hyperbolic lengths as a function of the Teichmüller parameter, we need a little geometry. We remark first that the shortest closed geodesic is always simple, by an easy surgery argument. Now given a simple closed geodesic $\alpha$ of length $\ell$ on a hyperbolic punctured torus $S$, cut $S$ along $\alpha$ to get a punctured cylinder with distance $h$ between its boundaries. There is a unique way to cut this cylinder along the shortest geodesic $\gamma$ between its boundaries and along geodesics orthogonal to the boundaries and heading into the cusp, to obtain two congruent pentagons, each with one ideal vertex and four right angles. A bit of hyperbolic trigonometry on this configuration yields (see Beardon [8] or Buser [20]):

(2.1) $$\cosh h/2 = \coth \ell/2.$$

If $\beta$ is a geodesic in $S$ that intersects $\alpha$ once, a symmetry argument shows that $\beta$ must intersect the segment $\gamma$ at its midpoint. Thus there are lifts $\widetilde{\alpha}$, $\widetilde{\beta}$ and $\widetilde{\gamma}$ to $\mathbf{H}^2$ which form a right triangle with legs $h/2$ and $t$, and hypotenuse $\ell'/2$, where $\ell'$ is the length of $\beta$ and $t$ is the distance along $\widetilde{\alpha}$ between its intersection with $\widetilde{\gamma}$ and with $\widetilde{\beta}$. The hyperbolic law of cosines gives us:

(2.2) $$\cosh \ell'/2 = \cosh h/2 \cosh t/2.$$

We deduce from this a number of things. First, we see that $t$ is a function of $\ell$ and $\ell'$. Thus, if $\ell = \ell'$ we would find that the configuration obtained by cutting along $\beta$ is equivalent to that obtained by cutting along $\alpha$. We conclude that if $\ell = \ell'$ there is an orientation-reversing isometry of $S$ that exchanges $\alpha$ and $\beta$.

Back in $\overline{\mathbf{H}}^2$, suppose that $\alpha = \infty$ and $\beta = n \in \mathbf{Z}$. The orientation reversing homeomorphism of $S$ interchanging them acts on $\mathbf{H}^2$ as the Möbius reflection through the axis $\{\tau : |\tau - n| = 1\}$. Thus if $\alpha$ and $\beta$ have equal hyperbolic lengths for some $\tau \in \mathbf{H}^2$ then $\tau$ lies on this semicircle. We know (by the collar lemma) that $\alpha$ is the shortest curve in $S$ when $\operatorname{Im}\tau$ is sufficiently large, and that a different curve can only become shortest at a point where it and $\alpha$ have equal lengths. It follows that the locus of $\mathbf{H}^2$ where $\alpha$ has strictly shorter hyperbolic length than any of its neighbors is exactly $H(\alpha) = \{\tau : \forall n \in \mathbf{Z}, |\tau - n| > 1\}$. Define $H(\gamma)$ for other $\gamma \in \hat{\mathbf{Q}}$ via the action of $\operatorname{SL}_2(\mathbf{Z})$. We will see momentarily that in fact $H(\alpha)$ is the locus where $\alpha$ is strictly shortest among all geodesics.

After applying this discussion to all vertices we find that, for any Farey triangle $\Delta$, if we divide up $\Delta$ into six regions by the axes of its reflection symmetries, then each vertex $u$ has minimal hyperbolic length in the pair of regions that meet $u$, and is strictly minimal in the interior of the union of the pair.

There is always, for general reasons, some constant $L_0 > 0$ bounding the length of the shortest geodesic in all finite-area hyperbolic surfaces with a fixed topology (see e.g. Buser [20]). In fact $L_0$ for the punctured torus is one of the few constants in this paper whose value we can compute precisely: If $\alpha$ and $\beta$ are as above then, after possibly Dehn-twisting $\beta$ about $\alpha$ a number of times, we may assume that $t$ in (2.2) is at most $\ell/2$ (this is easiest to see by considering the lifts of $\gamma$ that are crossed by $\widetilde{\beta}$ – successive ones are separated by at most $\ell/2$ along the translates of $\widetilde{\alpha}$). Now if $\alpha$ has minimal length among its neighbors then $\ell \le \ell'$ and so, combining (2.1) and (2.2) and using $t \le \ell/2$, we find that $\sinh \ell/2 \le \cosh \ell/4$. It follows that $\ell \le 4 \sinh^{-1}(1/2) \approx 1.9248$. Thus we let this be $L_0$ and it serves as an upper bound for the shortest $\alpha$.

We can also compute, for $\ell \le L_0$, that $h \ge h_0 \approx 1.609$. Thus any curve which crosses $\alpha$ more than once has length at least $2h_0 \approx 3.218$. We summarize our findings as follows:

LEMMA 2.1. *Let $S$ be a hyperbolic punctured torus. The number $L_0 \approx 1.9248$ bounds the length of the shortest closed geodesic on $S$. If $\alpha$ and $\beta$ have lengths bounded by $L_0$ on $S$ then they are Farey neighbors. In particular if they are both shortest on $S$ then they are Farey neighbors.*

Finally, it will be useful to know that the same value of $L_0$ bounds the shortest geodesic for *any* complete metric on $S$ with curvatures bounded above by $-1$ where the hole is conformally a puncture: By a lemma of Ahlfors (see [2]) the hyperbolic metric in the same conformal class as such a metric is pointwise bigger. The author is grateful to Curt McMullen for pointing out this lemma. In particular, for the metrics which will arise later on from simplicial hyperbolic surfaces, the curvatures are $-1$ except for isolated cone singularities with cone angle $2\pi$ or more, and these can be represented in an isothermal coordinate as zeroes of the conformal factor. Ahlfors' lemma applies in this generality.

*Teichmüller parameters for annuli.* A *marked annulus* is an oriented annulus $A$ with an arc $\beta$ whose endpoints lie on distinct boundary components. A conformal structure on such an annulus yields a Teichmüller parameter $\tau(A, \beta) \in \mathbf{H}^2$, namely the point corresponding to the marked torus obtained by placing a Euclidean metric on $A$, gluing the boundaries together with an isometry that identifies the endpoints of $\beta$, and marking with the core curve of $A$ and the image of $\beta$ under the gluing.

Let $T$ be a torus which is the union of a sequence of annuli $A_1, \ldots, A_k$ glued along their boundaries, and suppose that in a Euclidean metric on $T$ the boundaries of the $A_i$ are geodesic. Let $\alpha$ be a curve in the homotopy class of the boundary curves, and let $\mu$ be a curve crossing each annulus $A_i$ in a single arc $\mu_i$ (we do not require that $\mu$ be a geodesic). Then $\alpha, \mu$ give a marking for $T$, and we immediately have

$$(2.3) \qquad \tau(T, \alpha, \mu) = \sum_{i=1}^{k} \tau(A_i, \mu_i).$$

by considering the decomposition in the universal cover of $T$.

2.4. *Quasiconformal Lemmas.* Let us record some easy facts useful for estimating quasiconformal distortion in simple situations where we have mixed geometric and quasiconformal data.

Let $h : T_1 \to T_2$ be an $L$-Lipschitz map of degree 1 between Euclidean tori, where $T_1$ is a square torus of area 1, and $T_2$ has area at least $A_0$. Then $h$ is homotopic to a $K$-quasiconformal homeomorphism, where $K$ depends only on $L$ and $A_0$. This is easily seen by lifting to an isomorphism of lattices in $\mathbf{R}^2$.

A slightly more subtle fact is the following, which will be applied in Sections 9 and 10 to control Margulis tubes at the boundary of the convex core. A Euclidean annulus is an annulus isometric to the product of a circle with an interval. Call the length of the circle the *girth*.

LEMMA 2.2. *Let $h : A_1 \to A_2$ be a proper map between two Euclidean annuli of girth* 1, *which is $K$-quasiconformal on a Euclidean subannulus $C$ of $A_1$, whose modulus is at least $M_0$, and such that the components $B_0$ and $B_1$ of $A_1 \setminus C$ have modulus bounded by $M_0$. Suppose also that $h$ is $L$-Lipschitz on each $B_i$, and that $h$ is an embedding on the boundary of $A_1$ and $L$-bilipschitz on each boundary component.*

*Then $h$ is homotopic to a $K'$-quasiconformal homeomorphism $h' : A_1 \to A_2$, where $K'$ depends only on the constants $L, K, M_0$, and the homotopy is constant on the boundary.*

*Proof* (sketch). For $i = 0, 1$ let $B_i'$ be the Euclidean subannulus of $A_1$ containing $B_i$, such that $B_i' \setminus B_i$ has modulus $m = M_0/2$. We will replace $h|_{B_i'}$ with a quasiconformal map $h'$ which is homotopic to $h$ rel $\partial B_i'$, and whose quasiconformality constant depends only on the previous constants.

To see that this is possible, note that $h|_{B_i'}$ is an element of a family of maps $\{g : B_i' \to S^1 \times [0, \infty)\}$ which satisfy the same quasiconformality and Lipschitz conditions that $h|_{B_i'}$ does (here we are identifying $A_2$ with an appropriate initial subannulus $S^1 \times [0, T]$ of $S^1 \times [0, \infty)$).

Each such map $g$ can be deformed rel boundary to a $K(g)$-quasiconformal map $g'$, by standard methods: uniformize both $B_i'$ and $g(B_i')$ by the upper half-plane and check that the induced boundary map is quasisymmetric (this follows from the quasiconformality for one part of the boundary, and from the bilipschitz condition for the rest). Moreover, the family $\{g\}$ is compact in the compact-open topology, and after proper normalization, so is the family of lifts to the universal covers. It follows that the constants $K(g)$ have a uniform upper bound.

On $A_1 \setminus (B_0' \cup B_1')$, we simply keep the same map $h$. This concludes the proof. □

We remark also that in the argument, the process of lifting to the upper half-plane is what allows us to do the qc extension in a way that keeps proper track of the twisting of the original map. The intuition here is that the amount of twist in the original map $h$ is bounded by quasiconformality in $C$ and by the Lipschitz condition in $B_i$. This is subtle to detect directly – in fact, the twist itself can go to infinity with fixed $K$ if $\mathrm{mod}\,(A_1) \to \infty$. However, $K'$ is bounded independently of $\mathrm{mod}\,(A_1)$, and this is reflected in the quasisymmetry of the lift.

## 3. Geometric tameness and end invariants

In this section we describe how to associate to a punctured torus group an ordered pair of "end invariants" $(\nu_-, \nu_+)$, each lying in the closed disk $\overline{\mathbf{D}}$, or equivalently $\overline{\mathbf{H}}^2$. This is a special case of end invariants for general (geometrically tame) Kleinian groups, coming from the work of Ahlfors, Bers and Maskit for geometrically finite ends (where the invariant is a collection of simple closed curves on the boundary and an element of the Teichmüller space of their complement), and from Thurston, Bonahon and Canary for geometrically infinite ends (where the invariant is a geodesic lamination). We omit a general discussion of this, referring the reader to [14], [22], [71], [86] for more details.

Let us now concentrate on the case of a punctured torus group $\rho : \pi_1(S) \to \mathrm{PSL}_2(\mathbf{C})$ and associated manifold $N = \mathbf{H}^3/\rho(\pi_1(S))$.

Let $\check{N}$ denote $N$ minus the $\varepsilon_0$-Margulis tube $Q_N$ associated to the parabolic commutator (we call this the "main cusp"). This manifold has two ends; in fact, circumventing historical order we may note that Bonahon's theorem [14] implies that $N$ is homeomorphic to $S \times \mathbf{R}$, and $\check{N}$ is homeomorphic to $S_0 \times \mathbf{R}$, where $S_0$ is $S$ minus an open neighborhood of the puncture. (Remark: the ends can be defined without knowing Bonahon's theorem, for example by considering the way that $\check{N}$ is cut up by its relative Scott core [80], [65]. However, we prefer this simplified exposition). Let us name the ends $e_-$ and $e_+$, where the following orientation convention is applied:

If $M$ is an oriented manifold we orient $\partial M$ by requiring that the frame $(f, n)$ has positive orientation whenever $f$ is a positively oriented frame on $\partial M$ and $n$ is an inward-pointing vector. Now fix the orientation on $N$ induced from $\mathbf{H}^3$, and choose a fixed orientation for $S$. This determines (up to homotopy through proper maps) an identification of $N$ with $S \times (-1, 1)$ which induces the map $\rho$ on fundamental groups, and such that the orientation of $S$ agrees with that induced on $S \times \{1\}$. Let $e_+$ denote the end of $\check{N}$ whose neighborhoods are neighborhoods of $S_0 \times \{1\}$, and $e_-$ the other end.

Let $\Omega$ denote the (possibly empty) domain of discontinuity of $\Gamma$. Let $\overline{N}$ denote the quotient $(\mathbf{H}^3 \cup \Omega)/\Gamma$. Any component of the boundary $\Omega/\Gamma$ is reached by going to one of the ends $e_+$ or $e_-$, and this divides it into two disjoint pieces $\Omega_+/\Gamma$ and $\Omega_-/\Gamma$ (where $\Omega_+, \Omega_-$ are the corresponding invariant subsets of $\Omega$). There are three possibilities for each of these boundaries, corresponding to three types of end invariants (here let $s$ denote either $+$ or $-$):

1. $\Omega_s$ is a topological disk, and $\Omega_s/\Gamma$ is a punctured torus. This determines a point in the Teichmüller space of $S$, denoted by $\nu_s$.

2. $\Omega_s$ is an infinite union of disks and $\Omega_s/\Gamma$ is a thrice-punctured sphere, obtained from the corresponding boundary of $S \times (-1, 1)$ by deleting a simple closed curve $\gamma_s$. In this case $\nu_s \in \hat{\mathbf{Q}}$ denotes the slope of $\gamma_s$, as in Section 2.1. The conjugacy class of $\gamma_s$ in $\Gamma$ is parabolic.

3. $\Omega_s$ is empty. In this case $\nu_s \in \mathbf{R}\backslash\mathbf{Q}$; we describe its geometric significance below.

(This trichotomy is due to Maskit; see [59], [63]). Let $C(N)$ denote the convex core of $N$, namely the quotient by $\Gamma$ of the convex hull $CH(\Lambda)$ of the limit set $\Lambda$ of $\Gamma$. Each component of $\partial CH(\Lambda)$ corresponds to a component of $\Omega$ via orthogonal projection from $\hat{\mathbf{C}}$ to $\partial CH(\Lambda)$ (see [30]), so that $\partial C(N)$ divides naturally into $\partial_+ C(N) \cup \partial_- C(N)$ where each $\partial_s C(N)$ is a punctured torus, thrice-punctured sphere or empty according to the three cases above. Each boundary component is a convex pleated surface in $N$ with an induced hyperbolic metric.

In case 2, since all thrice-punctured sphere groups (with parabolic boundaries) are conjugate to a fixed Fuchsian group, the components of $\Omega_s$ are actually round circles, and the boundary component $\partial_s C(N)$ is totally geodesic.

To define $\nu_s$ in case 3 we need to recall the theory of ends due to Bonahon and Thurston. For a simple closed curve $\gamma$ in $S$ let $\gamma^*$ denote its geodesic representative in $N$ (more precisely $\gamma$ determines a conjugacy class taken by $\rho$ to a conjugacy class in $\Gamma$. If this class is nonparabolic it has a geodesic representative). Thurston showed [86] that if $\{\gamma_n\}$ is a sequence of simple closed curves such that $\gamma_n^*$ are eventually contained in any neighborhood of $e_s$, then the slopes of $\gamma_n$ converge in $\mathbf{R}$ to a unique irrational number. We say that such a sequence "exits the end". Bonahon [14] showed that in fact for each end $e_s$ that is not geometrically finite, that is, not in case 1 or 2, there is such a sequence of geodesics. (Thurston showed this for groups that are known to be limits of quasi-Fuchsian groups.) We define $\nu_s$ to be this limiting irrational slope. Thurston calls $\nu_s$ an *ending lamination* in this case because it describes a geodesic lamination for any hyperbolic metric on $S$, obtained as a limit of the geodesics in $S$ corresponding to $\gamma_j$.

The consequences of the existence of sequences of geodesics exiting an end will be more apparent once we introduce simplicial hyperbolic surfaces in Section 5.

## 4. The pivot theorem

In this section we associate to any end-invariant pair $(\nu_-, \nu_+)$ a *pivot sequence*, which is closely related to a continued-fraction expansion, and state our main structural theorem, Theorem 4.1, which translates the combinatorial structure of the sequence into geometric data about the associated representation $\rho$. The proof of Theorem 4.1 will be completed in Section 10, after a number of necessary tools are developed.

4.1. *The pivot sequence.* Letting $s$ denote $+$ or $-$, define a point $\alpha_s \in \hat{\mathbf{R}}$ to be *closest to* $\nu_s$ in the following sense: If $\nu_s \in \hat{\mathbf{R}}$ let $\alpha_s = \nu_s$. If $\nu_s \in \mathbf{D}$, let $\alpha_s \in \mathcal{C}$ represent a geodesic of shortest length (hence at most $L_0$ – see Lemma 2.1) in the hyperbolic structure corresponding to $\nu_s$. In particular $\nu_s$ is contained in a Farey triangle $\Delta_s$ of which $\alpha_s$ is a vertex. Note that (nongenerically) there may be two or three choices for $\alpha_s$, in which case we choose one arbitrarily. Our constructions will work for any of the choices.

Now define $E = E(\alpha_-, \alpha_+)$ to be the set of edges of the Farey graph which separate $\alpha_-$ from $\alpha_+$ in the disk. Let $P_0$ denote the set of vertices of $\mathcal{C}$ which belong to at least two edges in $E$ (except for one exceptional case, see (3) below). We call these vertices *internal pivots* (see Figure 3). The edges of $E$ admit a natural order where $e < f$ if $e$ separates the interior of $f$ from $\alpha_-$, and it is easy to see this induces an ordering on $P_0$. We therefore arrange $P_0$ as a sequence $\{\alpha_n\}_{n=\iota}^p$ where $\iota = -\infty$ if $\nu_- \in \mathbf{R} \setminus \mathbf{Q}$, and $\iota = 1$ otherwise, and $p = \infty$ if $\nu_+ \in \mathbf{R} \setminus \mathbf{Q}$, and is some finite nonnegative integer otherwise.

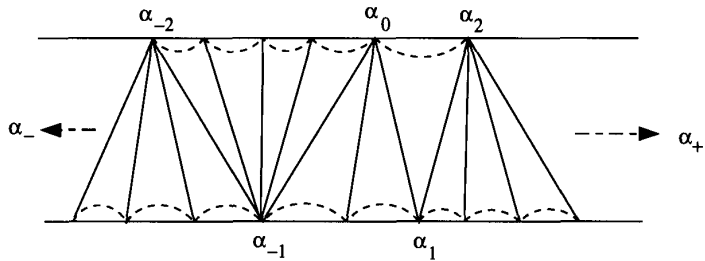

Figure 3. Sketch of a pivot sequence. For visibility the circle boundary has been stretched to two parallel lines and only the edges separating $\alpha_-$ from $\alpha_+$ are solid.

Note in particular the following special cases:

1. $E$ is bi-infinite and $P_0 = P$ is bi-infinite: both $\nu_\pm$ lie in the irrational part of the boundary.

2. $E = \emptyset$ and $P_0 = \emptyset$: $\alpha_\pm$ lie in the closure of a single Farey triangle. Note $\alpha_+$ and $\alpha_-$ may or may not be equal.

3. $E$ is a singleton $\{e\}$: $\alpha_\pm$ lie in the closures of adjacent Farey triangles. In this case we *redefine* $P_0$ to be a single point $\{\alpha_1\}$, which is chosen arbitrarily from the two endpoints of $e$.

Case (1) should be kept in mind throughout most of the paper, as it is easiest to deal with. Cases (2) and (3) are particularly simple types of geometrically finite groups, and hence in some sense we have nothing new to say about them; but they tend to complicate our exposition and notation. The geometrically finite cases are also responsible for introducing some real geometric subtleties, particularly in the behavior of geometric limits (see Section 11.2), and so we have no choice but to be careful.

We obtain the full pivot sequence $P$ by appending to the beginning of $P_0$ the vertex $\alpha_-$ if $\alpha_- \in \mathcal{C}$ (hence $\iota = 1$), and appending to the end of $P_0$ the vertex $\alpha_+$ if $\alpha_+ \in \mathcal{C}$ (hence $p < \infty$). In these cases we define $\alpha_0 \equiv \alpha_-$ and $\alpha_{p+1} \equiv \alpha_+$, respectively.

With this numbering convention, we note that for any $n \in \{\iota, \ldots, p\}$, $\alpha_{n-1}$ and $\alpha_{n+1}$ are related by

$$\alpha_{n+1} = D_{\alpha_n}^{w(n)} \alpha_{n-1},$$

and the integers $w(n)$ so defined are called the *widths* of the pivots.

It will be useful to consider this from a different perspective: For each $\beta \in \mathcal{C}$ fix an identification of $\mathbf{D}$ as $\mathbf{H}^2$ (by an orientation-preserving Möbius transformation) such that $\beta$ is identified with $\infty$, and its neighbors with the integers $\mathbf{Z}$. Such a normalization is unique only up to integer translation; the ambiguity will turn out not to matter. Let $\nu_+(\beta)$ and $\nu_-(\beta)$ denote the points of $\overline{\mathbf{H}}^2 = \mathbf{H}^2 \cup \hat{\mathbf{R}}$ to which $\nu_\pm$ are taken by this identification. Similarly let $\alpha_i(\beta)$ be the images of $\alpha_i$ by this identification, for $i \in \mathbf{Z}$ or $i = \pm$.

One easily checks that (except in case (3)) $\beta$ is in $P_0$ if and only if $\alpha_-(\beta)$ and $\alpha_+(\beta)$ are separated by at least two integers. The width $w(n)$ can be written $\alpha_{n+1}(\alpha_n) - \alpha_{n-1}(\alpha_n)$, and can also be estimated as follows:

For any $x \neq y \in \mathbf{R}$ which are separated by an integer, let $\lfloor x, y \rfloor$ be defined as $k - j$ where $j$ is the integer closest to $x$ in the closed interval spanned by $x$ and $y$, and $k$ is the integer closest to $y$ in this interval. In particular we note the sign of $\lfloor x, y \rfloor$ equals the sign of $y - x$, and $|y - x| - 2 < |\lfloor x, y \rfloor| \leq |y - x|$. With this notation we have

(4.1) $$w(n) = \lfloor \alpha_-(\alpha_n), \alpha_+(\alpha_n) \rfloor.$$

This easily gives the estimate

(4.2) $$|w(n) - (\operatorname{Re} \nu_+(\alpha_n) - \operatorname{Re} \nu_-(\alpha_n))| \leq 2.$$

We can also define widths for the first and last pivots of $P$, if these exist. Suppose $\alpha_+ \ne \alpha_-$. When $\alpha_+ = \alpha_{p+1}$ and $\nu_+ \ne \alpha_+$, we compute $w(p+1)$ using (4.1), but replacing the second term by $\operatorname{Re}\nu_+(\alpha_{p+1})$. When $\alpha_- = \alpha_0$ and $\nu_- \ne \alpha_-$, we compute $w(0)$ using (4.1), but replacing the first term by $\operatorname{Re}\nu_-(\alpha_0)$. If $\alpha_+ = \alpha_-$ we make both replacements in (4.1) to obtain $w(0)$. In all cases, the estimate 4.2 applies.

*Remark.* The connection of this to continued-fraction expansions is easiest to state if $\nu_- = \infty$ and $\nu_+ \in \mathbf{R}$. Then $|w(n)|$ are exactly the continued-fraction coefficients of the fractional part of $\nu_+$. The other cases are similar, and we omit the details.

4.2. *Statement of the Pivot Theorem.* For any element $\gamma$ of $\operatorname{PSL}_2(\mathbf{C})$, let $\lambda(\gamma) = \ell + i\theta$ denote its complex translation length; it is determined by the identity $\operatorname{tr}^2 \gamma = 4\cosh^2 \lambda/2$, and the normalizations $\ell \ge 0$ and $\theta \in (-\pi, \pi]$ (more about this choice in Section 6.2). Note that $\lambda$ is invariant under conjugation and inverse. Geometrically, $\ell$ (if positive) gives the translation length of $\gamma$ along its axis, and $\theta$ gives the rotation.

Thus, fixing a discrete faithful representation $\rho : \pi_1(S) \to \operatorname{PSL}_2(\mathbf{C})$, we obtain a function on $\mathcal{C}$ which we write $\lambda(\alpha) \equiv \lambda(\rho(\alpha))$.

The Pivot Theorem will give us quasi-isometric control of the complex translation lengths of the pivots of a punctured-torus group:

THEOREM 4.1 (Pivot Theorem). *There exist positive constants $\varepsilon, c_1$ such that, if $\rho$ is a marked punctured-torus group with associated pivot sequence,*

1. *If $\ell(\beta) \le \varepsilon$ then $\beta$ is a pivot.*

2. *If $\alpha$ is a pivot then*

$$\frac{2\pi i}{\lambda(\alpha)} \approx \nu_+(\alpha) - \overline{\nu}_-(\alpha) + i$$

*where "$\approx$" denotes a bound $c_1$ on hyperbolic distance in $\mathbf{H}^2$ between the left and right sides.*

*Remarks.* 1. The quantity $\omega(\alpha) = 2\pi i/\lambda(\alpha)$ is a convenient way to encode the geometry of $\alpha$ and its Margulis tube. In particular note that $\lambda$ lies in the right half-plane $\{\operatorname{Re} z > 0\}$, and $\omega$ lies in the upper half-plane $\{\operatorname{Im} z > 0\}$. Both are Teichmüller parameters for the torus $(\widehat{\mathbf{C}} \setminus \operatorname{Fix}(\alpha))/\alpha$, with respect to different markings; see Section 6.2 for more on this. The hyperbolic distance estimate on $\omega$ is natural because, being also a Teichmüller distance estimate, it implies a bilipschitz estimate on the action of $\alpha$ in $\mathbf{H}^3$ – see Lemma 6.2.

2. Note in particular that $|\omega(\alpha_n)|$ is bounded away from zero by (2), which implies a universal upper bound on the length of all pivots.

3. The real part of the estimate is just $\operatorname{Re}\nu_+(\alpha_n) - \operatorname{Re}\nu_-(\alpha_n)$, which is an estimate for $w(n)$, via (4.2). Furthermore, if $\alpha_n$ is an internal pivot we note that $\operatorname{Im}\nu_\pm(\alpha_n) \leq 1$, by the discussion in Section 2.3. It follows that (2) is equivalent, for internal pivots, to

$$\omega(\alpha_n) \approx w(n) + i. \tag{4.3}$$

In terms of $\lambda = \ell + i\theta$, (4.3) translates to:

$$\frac{c_2}{w(n)^2} \leq \ell(\alpha_n) \leq \frac{c_3}{w(n)^2} \tag{4.4}$$

and

$$\left| w(n) - \frac{2\pi}{\theta(\alpha_n)} \right| \leq c_4 \tag{4.5}$$

with suitable constants $c_i$ independent of $\rho$ or $\alpha_n$. It is easiest to see this by noting that the map $\omega \mapsto \lambda$ is an isometry between the Poincaré distance on the upper half-plane and that on the right half-plane.

4. For a noninternal pivot $\alpha = \alpha_\pm$, the imaginary part of $\nu_+(\alpha)$ or $\nu_-(\alpha)$ is large. One direction of the estimate (2) is then a variation of an inequality of Bers [11], reflecting the fact that a curve which is very short on the domain of discontinuity is very short in the 3-manifold. See Lemma 6.4 for more on this.

## 5. Simplicial hyperbolic surfaces

An important role in the theory of hyperbolic 3-manifolds is played by images of surfaces which are, in some sense, "hyperbolic". Thurston introduced this technique with his *pleated surfaces* (see Thurston [86], [87] and Canary-Epstein-Green [24]), and Bonahon and Canary [14], [22] have used a related construction known as simplicial hyperbolic surfaces. Sometimes the difference between these is merely technical, but in our case we find that the simplicial hyperbolic surfaces have a particular advantage: one can have better explicit control over continuous families of such surfaces. We will, however, briefly use pleated surfaces, in Sections 6 and 8, so we will discuss them here as well.

5.1. *Definitions.* Let $S$ be a (possibly punctured) surface and $N$ a hyperbolic 3-manifold. A proper map $f: S \to N$ is a *simplicial hyperbolic surface* if the following hold: A neighborhood of each puncture is mapped to a cusp of $N$. There is a triangulation $\mathcal{T}$ of $S$ (with some edges terminating at punctures) such that $f$ takes each edge to a geodesic, and each triangle of $\mathcal{T}$ to a totally geodesic immersed triangle (with punctures going to ideal vertices). The sum of corner angles around any nonideal vertex is at least $2\pi$.

In particular $f$ induces on $S$ a *singular hyperbolic metric*: that is, a complete metric $\sigma$ of finite area, smooth with curvature $-1$ away from finitely many singularities (at the vertices) where atoms of negative curvature may be concentrated.

We consider two simplicial hyperbolic surfaces to be equivalent if they differ only by precomposition with a homeomorphism of $S$ isotopic to the identity.

In the case of the punctured torus we will consider simplicial hyperbolic surfaces *adapted to a curve* $\alpha \in \mathcal{C}$ (Canary calls these surfaces with a distinguished edge). Given $\alpha$ in $\mathcal{C}$, realize it as a specific curve on $S$ (still called $\alpha$), choose a vertex $x \in \alpha$ and let $\beta$ be a simple curve meeting $\alpha$ transversely only at $x$. These curves cut $S$ into a punctured quadrilateral; adding four edges from the vertices to the puncture yields a triangulation $\mathcal{T}$ (see Figure 4).

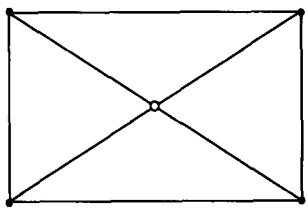

Figure 4. A triangulation of the punctured torus with one ideal vertex (in the center) and one real vertex

A simplicial hyperbolic surface $f : S \to N$ is adapted to $\alpha$ if it is simplicial with respect to this $\mathcal{T}$ and takes $\alpha$ to its geodesic representative $\alpha^*$ in $N$. It is easy to see that such maps exist in the homotopy class determined by $\rho$ [14], [22] and any two with isotopic triangulations differ only by "sliding" the vertex around $\alpha^*$ (perhaps more than once around) – that is, move (2) in Section 5.2 below.

A pleated surface is a map $f : S \to N$ which induces a *nonsingular* hyperbolic metric on $S$, with respect to which it is totally geodesic on the complement of a geodesic lamination (a closed set foliated by geodesics). The leaves of the lamination are also mapped geodesically. One can think of this heuristically as a simplicial hyperbolic surface where the triangulation has infinite-length edges and no vertices. For example one can show that, starting with a simplicial hyperbolic surface adapted to a curve $\alpha$ and performing the sliding operation infinitely many times, one obtains a pleated surface in the limit for which $\alpha$ is part of the lamination. In particular the singular hyperbolic metrics converge to the hyperbolic metric on the pleated surface.

Using simplicial hyperbolic surfaces we can obtain the following slightly stronger characterization of the ending lamination:

LEMMA 5.1. *There exists $L_0 > 0$ for which, if $\rho$ is a punctured torus group with (irrational) ending lamination $\nu_s$ (here $s$ denotes $+$ or $-$), there is a sequence $\gamma_n \in \mathcal{C}$ converging to $\nu_s$, such that the corresponding geodesics are eventually contained in any neighborhood of the end $e_s$, and $\ell(\gamma_n) \leq L_0$.*

*Proof.* Let $L_0$ be the bound on the shortest geodesic in a simplicial hyperbolic punctured torus, from Section 2.3. The Bonahon-Thurston theory gives a sequence of curves $\delta_n$ whose geodesic representatives are eventually in any neighborhood of $e_s$. Let $f_n : S \to N$ be simplicial hyperbolic surfaces adapted to $\delta_n$ (one could equally well use pleated surfaces) and let $\gamma_n$ be the shortest curve in the metric induced by $f_n$ on $S$. The statement follows for these curves. One needs only to check that $\gamma_n^*$ indeed exit the end $e_s$. Let $\hat{S}$ denote a fixed embedded cross-section of $N$. For large enough $n$, $f_n(\gamma_n)$ lies in a neighborhood of $e_s$ disjoint from $\hat{S}$. Let $A_n$ be a homotopy from $f_n(\gamma_n)$ to $\gamma_n^*$ which has geodesic tracks. The bound on the length of $f_n(\gamma_n)$ implies that either $A_n$ has bounded tracks, or it has a very short circumference for most of its length. If $\gamma_n^*$ is not contained in a neighborhood of $e_s$ disjoint from $\hat{S}$ then $A_n$ must meet $\hat{S}$, but then the translation length of $\gamma_n$ in a fixed compact set is small; this can only hold for finitely many $n$. □

5.2. *Interpolation of simplicial hyperbolic surfaces.* We recall at this point the *elementary moves* between simplicial hyperbolic surfaces discussed in Canary [23]. He uses three moves. In each case we begin with a simplicial hyperbolic surface $f_0 : S \to N$ with triangulation $\mathcal{T}_0$ and vertex $v$, adapted to a curve $\alpha_0$. A move replaces these data with $f_1, \mathcal{T}_1, v, \alpha_0$, and gives a homotopy $f_t : S \to N$, $t \in [0, 1]$, connecting the two maps so that each $f_t$ is still a simplicial hyperbolic surface.

1. Diagonal switch: Let $Q$ be a quadrilateral in $\mathcal{T}_0$ with diagonal $d$. In $\mathcal{T}_1$, $d$ is replaced by the opposite diagonal $d'$. The maps $f_t$ agree with $f_0$ everywhere but on the interior of $Q$, where for $t \in (0, 1)$ the triangulation contains both diagonals and a new vertex where their interiors intersect. (See Figure 5.)

2. Vertex slide: The vertex $v$ is "pushed" once around $\alpha_0$. The new triangulation is actually isotopic to the old one, but not rel $v$.

3. Geodesic switch: Here $\mathcal{T}_0$ and $\mathcal{T}_1$ are equal but the map $f_1$ is adapted to the other closed curve in the triangulation. (See Figure 6.)

Let us describe move (3) in more detail. Lift the image $f_0(v)$ of the vertex to $\xi \in \mathbf{H}^3$. If $\alpha_0$ and $\alpha_1$ are the two closed curves of the triangulation, the lifts of $\alpha_0$ and $\alpha_1$ based at $\xi$ determine group elements $\mathbf{A}$ and $\mathbf{B}$, respectively (compare §7). Let $P$ be the common perpendicular of their axes $T_\mathbf{A}$ and $T_\mathbf{B}$.

The homotopy $f_t$ is in three parts. Lifted to $\mathbf{H}^3$, it moves $\xi$ first along $T_\mathbf{A}$ to $P$, then along $P$ to $T_\mathbf{B}$, and then (if desired) along $T_\mathbf{B}$ to a new position. The rest of the triangulation varies accordingly, so that at each point the edges $\alpha_0$ and $\alpha_1$ map to geodesics and the map is simplicial hyperbolic.

Any two triangulations of the punctured torus (in fact any surface) can be connected by such elementary moves. In particular Figure 5 shows how a Dehn twist is effected. This is sufficient for traversing the edges of our graph $\mathcal{C}$. A consequence of this is the following lemma, proved in [23]:

LEMMA 5.2. *If $f_0, f_1$ are two homotopic simplicial hyperbolic surfaces then they may be connected by a continuous family $f_t, t \in [0,1]$ such that $f_t$ is a simplicial hyperbolic surface for each $t$.*

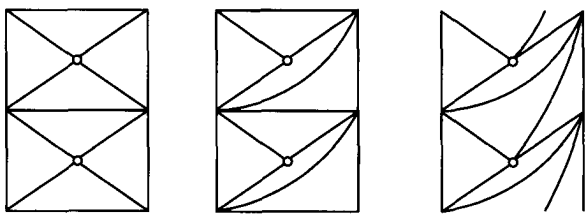

Figure 5. A sequence of two diagonal switches which effects a Dehn twist $D_{\alpha_0}$ on the triangulation. (Two fundamental domains in the abelian cover are shown, and $\alpha_0$ is the vertical curve.)

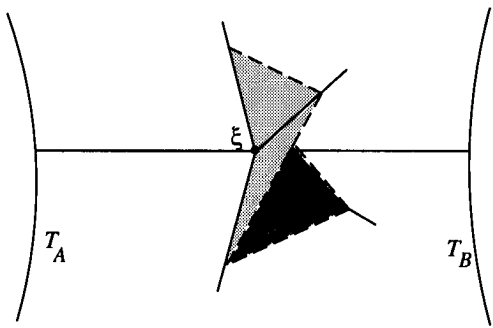

Figure 6. An intermediate stage in move (3). Pictured in the universal cover are $T_\mathbf{A}$ and $T_\mathbf{B}$, and a portion of the simplicial hyperbolic surface near $\xi$.

## 6. The structure of Margulis tubes

6.1. *Thick-thin decomposition.* The $\varepsilon$-thin part $M_{\text{thin}(\varepsilon)}$ of a manifold $M$ (see e.g. [88], [9], [7]) is the subset where the injectivity radius is at most $\varepsilon/2$. (Our epsilons always denote length of the shortest nontrivial closed path through a point, rather than injectivity radius.) We denote by $M_{\text{thick}(\varepsilon)}$ the closure of the complement of $M_{\text{thin}(\varepsilon)}$. The Margulis lemma [52], or Jørgensen's inequality [48], implies that for each dimension $n$ there is a value of $\varepsilon$, known as the Margulis constant, below which the $\varepsilon$-thin part of a hyperbolic $n$-manifold is of a standard type. In particular in dimensions 2 and 3 in the orientable case, every component, known as a Margulis tube, is either a tubular neighborhood of a closed geodesic, or the quotient of a horoball by an abelian parabolic group.

Let $\varepsilon_0$ be a Margulis constant for $\mathbf{H}^3$. Let $\mathbf{T}_\varepsilon(g)$ denote the $\varepsilon$-Margulis tube in $N = \mathbf{H}^3/\Gamma$ for a conjugacy class $g$ in $\Gamma$. Similarly with a fixed representation $\rho$ in mind and a conjugacy class $\alpha$ in $\pi_1(S)$, let $\mathbf{T}_\varepsilon(\alpha)$ denote $\mathbf{T}_\varepsilon(\rho(\alpha))$.

Let $r(\alpha)$ be the radius of $\mathbf{T}_{\varepsilon_0}(\alpha)$, that is, the distance from the geodesic core to the boundary. Recall that $\lambda(\alpha)$ denotes the complex translation length $\ell + i\theta$ of $\alpha$. Let us record some preliminary bounds relating these quantities, due to Meyerhoff [68] and Brooks-Matelski [19].

LEMMA 6.1. *There is a function $R(\varepsilon)$ with $\lim_{\varepsilon \to 0} R(\varepsilon) = \infty$, such that for an $\varepsilon_0$ Margulis tube $\mathbf{T}_{\varepsilon_0}(\alpha)$ in any hyperbolic 3-manifold, the distance $d(\varepsilon)$ between $\partial \mathbf{T}_\varepsilon(\alpha)$ and $\partial \mathbf{T}_{\varepsilon_0}(\alpha)$ for any $\varepsilon \in [\ell(\alpha), \varepsilon_0]$ is bounded below by*

$$(6.1) \qquad d(\varepsilon) \geq R(\varepsilon).$$

*In particular for $\varepsilon = \ell(\alpha)$ this gives $r(\alpha) \geq R(\ell(\alpha))$.*

*Proof* (sketch). In fact $R(\varepsilon)$ can be written as $\frac{1}{2}\log(1/\varepsilon) - c$ for some positive $c$. Consider first the lower bound on $r(\alpha)$. Brooks-Matelski apply Jørgensen's inequality to show that $r(\alpha)$ is at least $\log c'/|\lambda(\alpha)|$, for a fixed $c' > 0$.

In general, $|\lambda|$ may be much larger than $\ell$ because of the rotational part. Brooks-Matelski show that some function $R(\ell(\alpha))$ exists by a limiting argument. Meyerhoff uses a pigeonhole principle argument to find an iterate of $\alpha$ with small rotational and translational part, and applies Jørgensen's inequality to this iterate. He also gives a finer estimate using an argument of Zagier.

To get (6.1) for general $\varepsilon$, it suffices to observe that the radial projection from an $r_1$-equidistant surface of a geodesic in $\mathbf{H}^3$ to an $r_2$-equidistant surface, where $r_2 > r_1$, expands by at most $\cosh r_2 / \sinh r_1$ (and a similar bound for horospheres in the parabolic case). □

6.2. *Teichmüller parameters for Margulis tubes.* The shape of the Margulis tube $\mathbf{T}_{\varepsilon_0}(\alpha)$ is determined completely by the complex length $\lambda(\alpha)$, but this dependence is subtle in general (for example if $\ell$ is extremely small compared to $\theta$ then this shape will be sensitive to the continued-fraction expansion of $\theta/2\pi$).

To improve our picture a little, we will find it useful to define a reciprocal invariant

$$\omega(\alpha) = \frac{2\pi i}{\lambda(\alpha)} \tag{6.2}$$

whose geometric interpretation is as follows: When $\alpha$ is nonparabolic, a representative $g$ of the conjugacy class $\rho(\alpha)$ acts on $\mathbf{H}^3 \cup (\widehat{\mathbf{C}} \setminus \mathrm{Fix}(g))$ with quotient a solid torus $H$, with torus boundary $F_\infty$. Let $m$ be the unique homotopy class in $F_\infty$ of a meridian of $H$, and let $\alpha'$ be homotopic to the core $\alpha$ (more about the choice of $\alpha'$ in a moment). The pair $(\alpha', m)$ determines a marking of $F_\infty$; so together with the conformal structure inherited from $\widehat{\mathbf{C}}$, we obtain a point in the Teichmüller space of the torus, viewed as the upper half plane $\mathbf{H}^2$ (see Section 2.3). This point is exactly $\omega$, if $\alpha'$ is chosen properly:

The choice of $\alpha'$ is only determined up to a Dehn twist around $m$, which gives a parabolic action $\omega \mapsto \omega/(1+\omega)$ on the upper half-plane $\mathbf{H}^2$, corresponding to the action $\lambda \mapsto \lambda + 2\pi i$ in the right half-plane. The normalizing convention $\theta \in (-\pi, \pi]$ which we imposed in Section 4 corresponds to the condition that $\alpha'$ be of minimal length among the possible choices (strictly minimal when $|\theta| < \pi$, and the boundary choice $\theta = \pi$ taking care of the cases where there is not a unique shortest choice). This constrains $\omega$ to the region $\{\omega : |\omega - 1| \geq 1, |\omega + 1| > 1\}$, certainly the natural choice when $\ell$ is small, which will always be the case for us.

Let $F_0$ denote the torus boundary $\partial \mathbf{T}_{\varepsilon_0}(\alpha)$, with its induced Euclidean metric. The natural identification between $F_0$ and $F_\infty$ via radial projection from the axis induces a marking on $F_0$. If $r(\alpha)$ is sufficiently large, moreover, we find that $F_0$ and $F_\infty$ are uniformly quasiconformally diffeomorphic via this projection (in fact the constant is $\coth r$). Thus, letting $\omega(F_0)$ denote the Teichmüller parameter of $F_0$ with the same normalization, we find that $\omega$ and $\omega(F_0)$ are within bounded hyperbolic (i.e. Teichmüller) distance in $\mathbf{H}^2$.

This means that we can use the (marked) geometry of $F_0$ or $F_\infty$ interchangeably to describe the geometry of $\mathbf{T}_{\varepsilon_0}(\alpha)$ up to quasi-isometry. In fact we note:

LEMMA 6.2. *An approximation of $\omega(\alpha)$ up to hyperbolic distance $D$ in $\mathbf{H}^2$ yields an approximation of $\mathbf{T}_{\varepsilon_0}(\alpha)$ up to $K(D)$-bilipschitz diffeomorphism.*

*Proof.* As observed above, if two values of $\omega$ differ by a bounded distance then the corresponding loxodromic elements are quasiconformally conjugate at infinity with bounded constant. An extension theorem like Douady-Earle [29] suffices to relate the corresponding Margulis tubes. Alternatively, one can just demonstrate this by hand, working in cylindrical coordinates around the axis. □

When $\alpha$ is parabolic we have $\lambda = 0$, so we write $\omega = \infty i$. We think of this as the limit obtained when $\operatorname{Im}\omega \to \infty$. (The limit of a sequence of Margulis tubes with $\operatorname{Im}\omega$ bounded and $\operatorname{Re}\omega \to \infty$ is more subtle – see Section 11.2.) The torus $F_\infty$ becomes an infinite annulus with core homotopic to $\alpha$.

6.3. *Margulis tubes in surface groups.* As Thurston first observed, there are special constraints on the geometry of Margulis tubes in hyperbolic 3-manifolds homotopy-equivalent to a surface, which arise from the presence of pleated surfaces near every point in the convex hulls of these manifolds. In this section we fix a surface $S$ of finite type and let $\rho : \pi_1(S) \to \operatorname{PSL}_2(\mathbf{C})$ be a discrete faithful type-preserving representation, or "marked surface group". As usual $N = \mathbf{H}^3/\rho(\pi_1(S))$.

We record these facts, due to Thurston and Bonahon, in the following lemma, whose proof we sketch.

LEMMA 6.3. *There are constants $\varepsilon_1, a_0$ and $A_0$ depending only on $S$ such that, for any marked surface group $\rho : \pi_1(S) \to \operatorname{PSL}_2(\mathbf{C})$, and any primitive element $\gamma \in \pi_1(S)$ such that $\ell(\rho(\gamma)) \leq \varepsilon_1$,*

1. *$\gamma$ represents a simple curve in $S$.*

2. *The translation length of $\gamma$ on $\partial \mathbf{T}_{\varepsilon_0(\gamma)}$ is at most $a_0$.*

3. *The area of $\partial \mathbf{T}_{\varepsilon_0}(\gamma)$ is at least $A_0$.*

*Proof* (sketch). It follows from the work of Bonahon [14] and Thurston [86] that any point in the convex core $C(N)$ is within bounded distance of a pleated surface $f : S \to N$ in the homotopy class of $\rho$. (In fact Thurston's geometric tameness property, which Bonahon established for any surface group, implies that there is a continuously parametrized family of surfaces sweeping out $C(N)$, each of which is either a pleated surface or a bounded deformation of one.) Let $f$ be such a surface. If $X$ is a component of the $\varepsilon_0$-thick part of $S$ (in the induced metric) then $\pi_1(X)$ is nonabelian and thus, since $f$ is $\pi_1$-injective, the image $f(X)$ cannot be contained in $\mathbf{T}_{\varepsilon_0}(\gamma)$. Since there is a uniform bound on the diameter of $X$ depending only on the topological type of $S$, Lemma 6.1 then implies that for $\varepsilon_1$ sufficiently small (depending on the function $R(\varepsilon)$) only the $\varepsilon_0$-thin part of a pleated surface can touch the $\varepsilon_1$-thin part of $N$. Conclusion (1) follows. See also Otal [76], who shows furthermore that a sufficiently short geodesic in $N$ is actually unknotted in a natural sense.

For part (2), let $f : S \to N$ be a pleated surface which meets $\mathbf{T}_{\varepsilon_1}(\gamma)$. Since the boundary of the $\varepsilon_0$-thin part of $S$ (in the induced metric) must again map outside $\mathbf{T}_{\varepsilon_1}(\gamma)$, the translation length of $\gamma$ on $\partial \mathbf{T}_{\varepsilon_1}(\gamma)$ is at most $\varepsilon_0$. Since the distance from $\partial \mathbf{T}_{\varepsilon_1}(\gamma)$ to $\partial \mathbf{T}_{\varepsilon_0}(\gamma)$ is uniformly bounded above, we obtain a bound on the translation distance on $\partial \mathbf{T}_{\varepsilon_0}(\gamma)$ as well.

For part (3), let $\gamma'$ be the geodesic core of the Margulis tube in $S$ corresponding to $\gamma$. There exists a geodesic arc $\beta \subset S$ with both endpoints on $\gamma'$, whose intersection with the thick part of $S$ is a nonempty arc of length bounded by $D$, where $D$ depends again only on the topological type of $S$. Since $\beta$ is not deformable rel endpoints into $\gamma'$, its image cannot be deformed into $\mathbf{T}_{\varepsilon_0}(\gamma)$ in $N$.

Lifting to the universal cover, it follows that within distance $D$ of a lift $\widetilde{\mathbf{T}}$ of the tube $\mathbf{T}_{\varepsilon_0}(\gamma)$ there is a translate $\widetilde{\mathbf{T}}'$. The translates of $\widetilde{\mathbf{T}}'$ by the isometry $g$ stabilizing $\widetilde{\mathbf{T}}$ are all disjoint. Since the radii of these tubes are bounded below in terms of $\varepsilon_1$ via (6.1), one can show that their radial projections to $\partial \widetilde{\mathbf{T}}$ each have area at least $A(D, \varepsilon_1) > 0$. Since the projections are all disjoint this gives a lower bound on the area of the quotient torus. This proves (3). □

We can reinterpret these results in terms of the reciprocal invariant $\omega$.

LEMMA 6.4. *If $\rho$ is a surface group and $\alpha$ is a curve of length $\ell(\alpha) \leq \varepsilon_1$, then*

$$\operatorname{Im} \omega(\alpha) \geq c_1 \tag{6.3}$$

*where $c_1$ depends only on the topological type of the surface. Now suppose that $\rho$ is a punctured-torus group, and let $\nu_+(\alpha)$ and $\nu_-(\alpha)$ be the end invariants of $\rho$, in a normalization of $\mathbf{H}^2$ where $\alpha = \infty$. Then for any $\alpha \in \mathcal{C}$,*

$$\operatorname{Im} \omega(\alpha) \geq \operatorname{Im} \nu_+(\alpha) + \operatorname{Im} \nu_-(\alpha). \tag{6.4}$$

The last inequality (6.4) is actually just a variation of Bers' inequality $1/\ell(\alpha) \geq \frac{1}{2}(1/\ell_+(\alpha) + 1/\ell_-(\alpha))$ (see [11]). Bers' inequality holds in much greater generality, but we state it in the case of the punctured torus to indicate how it is expressed using our notation for $\nu_\pm(\alpha)$.

$\operatorname{Im} \omega$ is just the modulus of the torus $F_\infty$ cut along the Euclidean geodesic $\alpha'$ representing $\alpha$ (see the discussion in Section 6.2). Equivalently it is the reciprocal of the extremal length of $\alpha'$ in $F_\infty$.

As we will also see, $\operatorname{Im}(\omega)$ estimates the area of the torus $F_0$. In Sections 8 and 10 we will establish one of the central geometric facts about punctured torus groups – that $\operatorname{Im} \omega$ is also bounded above if the Margulis tube of $\alpha$ lies in the convex core of $N$, or equivalently that such a Margulis tube boundary has bounded area.

*Proof of Lemma* 6.4. Assume that $\ell(\alpha) \leq \varepsilon_1$. The Euclidean translation length of $\alpha$ on $F_0$ is bounded below by $\varepsilon_0$, by definition, and above by $a_0$, as a consequence of Part (2) of Lemma 6.3. Since $\alpha'$ is chosen to be the shortest representative of $\alpha$ in $F_0$, the same bounds hold for $\alpha'$.

The modulus $\operatorname{Im}\omega(F_0)$ is exactly $\operatorname{Area}(F_0)/\ell^2_{F_0}(\alpha')$. Thus, the length bounds on $\alpha'$ imply that $\operatorname{Im}\omega(F_0)$ and $\operatorname{Area}(F_0)$ are within bounded ratio.

Furthermore the bounded hyperbolic distance between $\omega(F_0)$ and $\omega(\alpha)$ implies that $\operatorname{Im}\omega(\alpha)$ and $\operatorname{Im}\omega(F_0)$ are within bounded ratio. The lower bound on $\operatorname{Im}\omega(\alpha)$ now follows from Part (3) of Lemma 6.3.

To conclude, let $\rho$ be a punctured torus group, and let us summarize the proof of our variation (6.4) of Bers' inequality. Recall that $g$ is a representative of $\alpha$. If $\operatorname{Im}\nu_+(\alpha) > 0$ then $\Omega_+$ is nonempty, and its quotient by $g$ is an annulus of modulus at least $\operatorname{Im}\nu_+$ which embeds in the torus $F_\infty$. The same is true for $\nu_-$, and monotonicity of moduli (by the method of extremal length) implies that the modulus of the torus with distinguished curve $\alpha$, namely $\operatorname{Im}\omega$, is at least the sum of the moduli of the two annuli. (Note that actually Bers' inequality is slightly stronger since it involves the hyperbolic lengths $\ell_\pm(\alpha)$, which satisfy $\pi/\ell_\pm(\alpha) \geq \operatorname{Im}\nu_\pm(\alpha)$. However this inequality is nearly an equality for large $\operatorname{Im}\nu_\pm$, and so this form will suffice for us.)

We also remark that the same argument applies in the case that $\alpha$ is parabolic in $\Omega_+$ or $\Omega_-$, so that $\operatorname{Im}\nu_+ + \operatorname{Im}\nu_- = \infty$. In this case the conclusion is that the quotient $\widehat{\mathbf{C}} \setminus \operatorname{Fix}(g)$ must have been an annulus rather than a torus, so that $\rho(\alpha)$ is parabolic. □

A final corollary of this is the inequality

$$\tag{6.5} |\operatorname{tr}(\alpha)| \geq c_2$$

for a uniform $c_2 > 0$ over all curves $\alpha$ in a surface group, since $\operatorname{tr} \to 0$ implies $\lambda \to \pm\pi i$, which means $\omega \to 2$, contradicting the lower bound on $\operatorname{Im}\omega$. We will use this briefly in Section 10.

## 7. Geometry of Farey neighbors

The purpose of this section is to investigate some constraints on the geometric relationship between a pair of generators for a punctured torus group. In particular it will follow from Lemma 7.1 that if both of a pair of generators have bounded translation lengths, there is a point in $\mathbf{H}^3$ where their actions are simultaneously bounded.

A pair $\alpha, \beta \in \mathcal{C}$ with $\alpha \cdot \beta = 1$ determines a pair of generators $\mathbf{A}, \mathbf{B}$ for $\Gamma$, up to conjugation and inversion, as follows: Realize $\alpha$ and $\beta$ as curves with a single intersection point $x$, choose an orientation for each and let $\mathbf{A}, \mathbf{B}$ be the images under $\rho$ of the resulting elements in $\pi_1(S, x)$. Fixing a positive number $D$, let $\mathrm{core}_D(\mathbf{A})$ (similarly $\mathrm{core}_D(\mathbf{B})$) denote the subset $\{p \in \mathbf{H}^3 : d(p, \mathbf{A}p) \leq \ell(\mathbf{A}) + D\}$. Note that if $\ell(\mathbf{A}) > \varepsilon_0$ then $\mathrm{core}_D(\mathbf{A})$ is a bounded tubular neighborhood of the axis of $\mathbf{A}$, whereas if $\ell(\mathbf{A}) \leq \varepsilon_0$ then $\mathrm{core}_D(\mathbf{A})$ is a horoball or tubular neighborhood of the axis of $\mathbf{A}$, contained in a bounded neighborhood of the Margulis tube of $\mathbf{A}$.

The parabolic commutator condition serves to force these cores together:

LEMMA 7.1. *There exists $D_2$ independent of the punctured-torus group $\rho$ such that, for any generator pair $\mathbf{A}, \mathbf{B}$, $\mathrm{core}_{D_2}(\mathbf{A}) \cap \mathrm{core}_{D_2}(\mathbf{B}) \neq \emptyset$.*

*Proof.* Assume first that neither $\mathbf{A}$ nor $\mathbf{B}$ is parabolic. The following familiar construction appears for example in Parker-Series [78]. Let $\mathcal{P}$ and $\mathcal{Q}$ be the axes of $\mathbf{A}$ and $\mathbf{B}$ respectively, and $\mathcal{M}$ their common perpendicular. Let $M$ denote the half-turn around $\mathcal{M}$. Then there are unique lines $\mathcal{N}$ orthogonal to $\mathcal{P}$ and $\mathcal{L}$ orthogonal to $\mathcal{Q}$ such that $\mathbf{A} = NM$ and $\mathbf{B} = ML$, where $N$ and $L$ are the corresponding half-turns. Hence $\mathbf{AB} = NL$, and its axis is the common perpendicular $\mathcal{R}$ of $\mathcal{L}$ and $\mathcal{N}$ (see Figure 7). Note that $\mathcal{L}$ and $\mathcal{N}$ cannot coincide since the group is not elementary.

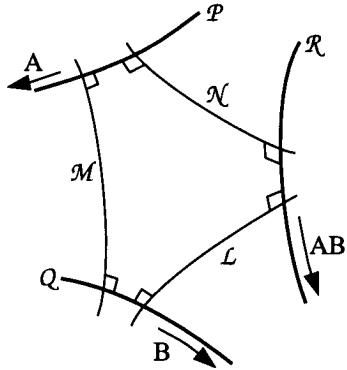

Figure 7. The hexagonal configuration of axes determined by $\mathbf{A}$ and $\mathbf{B}$.

The complex distance between $\mathcal{N}$ and $\mathcal{M}$ is given by $\mu = \lambda(\mathbf{A})/2$ and between $\mathcal{L}$ and $\mathcal{M}$ by $\tau = \lambda(\mathbf{B})/2$. Similarly $\sigma = \lambda(\mathbf{AB})/2$ gives the complex distance between $\mathcal{L}$ and $\mathcal{N}$. Let $\delta$ denote the complex distance between $\mathcal{P}$ and $\mathcal{Q}$. Note that there is an ambiguity of $\pi i$ in these. However (see Fenchel [32] or Kourouniotis [56]), one can give orientations for all the lines, which resolve the ambiguities in such a way that the hexagon cosine law holds:

$$(7.1) \qquad \cosh \sigma = \cosh \mu \cosh \tau + \sinh \mu \sinh \tau \cosh \delta.$$

Furthermore, these choices determine matrix representatives for $L, M$ and $N$, and hence for $\mathbf{A}, \mathbf{B}$ and $\mathbf{AB}$, for which the identities $\mathrm{tr} = 2\cosh\lambda/2$ hold without sign ambiguity. The standard trace identity

$$(7.2) \qquad \mathrm{tr}^2 \mathbf{A} + \mathrm{tr}^2 \mathbf{B} + \mathrm{tr}^2 \mathbf{AB} - \mathrm{tr}\mathbf{A}\,\mathrm{tr}\mathbf{B}\,\mathrm{tr}\mathbf{AB} = 2 + \mathrm{tr}[\mathbf{A},\mathbf{B}]$$

can therefore be rewritten in terms of $\mu, \tau$ and $\sigma$. Solving this simultaneously with (7.1), we arrive at the relation

$$(7.3) \qquad \sinh^2\delta\,\sinh^2\mu\,\sinh^2\tau = \frac{\mathrm{tr}[\mathbf{A},\mathbf{B}] - 2}{4}.$$

Furthermore we have $\mathrm{tr}[\mathbf{A},\mathbf{B}] = -2$, as a consequence of the parabolic commutator condition. More precisely, the sign of $\mathrm{tr}[\mathbf{A},\mathbf{B}]$ is well-defined since it is a commutator, and it must be $-2$ rather than $+2$ so that the group will be nonelementary (see e.g. [46], [50], [16]). It follows that

$$(7.4) \qquad \sinh^2\delta\,\sinh^2\mu\,\sinh^2\tau = -1$$

(cf. (2.2) of Parker-Series). This will give us an estimate for the length $\mathrm{Re}\,\delta$ of the segment $G$ of $\mathcal{M}$ between $\mathcal{P}$ and $\mathcal{Q}$.

We will next determine the subsegments of $G$ belonging to $\mathrm{core}_D(\mathbf{A})$ and $\mathrm{core}_D(\mathbf{B})$ respectively, for appropriate $D$, and use the estimate on $\delta$ to show that these segments cover $G$, and in particular have the desired intersection point.

Let $\mathcal{N}^+$ denote the ray of $\mathcal{N}$ with endpoint on $\mathcal{P}$ whose projection to $\mathcal{M}$ is a subinterval $G_1$ of $G$, and let $Z_1$ be the geodesic from the endpoint of $\mathcal{N}^+$ at infinity to its projection on an endpoint of $G_1$ (see Figure 8). Let $\delta_1$ denote the complex distance between $\mathcal{P}$ and $Z_1$. Another application of hyperbolic trigonometry to the quadrilateral $X_1$ bounded by $G_1, Z_1, \mathcal{N}^+$ and $\mathcal{P}$ (see any of [32], [56] or [78] for how to derive this from the hexagon cosine law by considering degenerate hexagons) yields

$$(7.5) \qquad \begin{aligned} \sinh^2\delta_1\sinh^2\mu &= 1, \\ \sinh^2\delta_2\sinh^2\tau &= 1, \end{aligned}$$

where $G_2, \delta_2$ are defined similarly with respect to $\mathcal{Q}$ and $\mathcal{L}$.

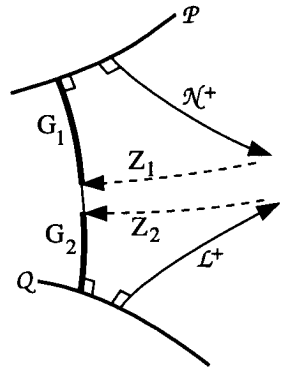

Figure 8. Projections of $\mathcal{N}$ and $\mathcal{L}$ to $G$ yield segments $G_1, G_2$ which are in the cores of **A** and **B**.

We claim that

(7.6) $$G_1 \subset \mathrm{core}_D(\mathbf{A})$$

where $D = 4\sinh^{-1}(1)$, and similarly for $G_2$ and $\mathrm{core}_D(\mathbf{B})$. It is easy to see (an ideal triangle is the extremal case) that for any triangle $\Delta(xyz)$ and point $p$ on a side $xy$, the distance from $p$ to at least one of the other two sides is no more than $\sinh^{-1}(1)$. (In fact this is the proof that hyperbolic space is "delta-hyperbolic" in the sense of Gromov). Similarly, if the triangle is a right triangle and $p$ is on a leg, the same bound holds for the distance from $p$ to the hypotenuse. Applying this to the quadrilateral $X_1$ suitably subdivided into two triangles, we find that for any $p \in G_1$, $\mathrm{dist}(p, \mathcal{P} \cup \mathcal{N}) \le 2\sinh^{-1}(1)$. If this holds for the distance to $\mathcal{P}$ then clearly $\mathrm{dist}(p, \mathbf{A}p) \le 4\sinh^{-1}(1) + \ell(\mathbf{A})$. If it holds for the distance to $\mathcal{N}$ then, since $\mathbf{A} = NM$ and $M$ fixes $p$, we have $\mathrm{dist}(p, \mathbf{A}p) \le 4\sinh^{-1}(1)$. Either way $p \in \mathrm{core}_D(\mathbf{A})$, so that (7.6) is established. The same argument applies to $G_2$.

Finally, we show that $G \setminus (G_1 \cup G_2)$ is (if nonempty) a segment of bounded length. Indeed, (7.4) and (7.5) together imply

(7.7) $$\sinh^2 \delta = -\sinh^2 \delta_1 \sinh^2 \delta_2.$$

Now taking absolute values and applying the elementary inequality $\frac{1}{2}(e^{|\mathrm{Re}\,z|} - 1) \le |\sinh z| \le e^{|\mathrm{Re}\,z|}$ and a bit of algebra, we deduce

(7.8) $$|G| \le |G_1| + |G_2| + \log 3$$

where $|\cdot|$ denotes length in $\mathbf{H}^3$.

Thus the $\frac{1}{2}\log 3$-neighborhoods of $G_1$ and $G_2$ meet, and are contained in $\mathrm{core}_{D_2}(\mathbf{A})$ and $\mathrm{core}_{D_2}(\mathbf{B})$, respectively, where $D_2 = 4\sinh^{-1}(1) + \log 3$. This concludes the proof when **A** and **B** are loxodromic.

If one or both of them are parabolic, we can make a limiting argument to obtain a $G$ which may be an infinite ray or all of $\mathcal{M}$, and $G_1$ and $G_2$ which may be subrays. The $G_i$ are still contained in $\mathrm{core}_D(\mathbf{A})$ and $\mathrm{core}_D(\mathbf{B})$, respectively, and the length bound on $G\setminus(G_1\cup G_2)$ holds as before. Thus the lemma follows in general. □

## 8. Connectivity and *a priori* bounds

In this section we exploit the connection between topological properties of hyperbolic surfaces, and combinatorial properties of the graph $\mathcal{C}$, to obtain our initial bound on the length of pivots. We first prove Lemma 8.1, which shows that the set of curves in $\mathcal{C}$ satisfying a length bound is connected. As a consequence we obtain Lemma 8.2 which gives a preliminary bound on the lengths of all vertices in the pivot sequence.

8.1. *Connectivity in the curve complex.* Let $L_0$ be the bound on the length of the shortest geodesic in a punctured torus with a singular hyperbolic metric (see Lemma 2.1). Fixing a punctured-torus group $\rho$, for any $L \geq L_0$, let $\mathcal{C}(L)$ denote the subgraph of $\mathcal{C}$ spanned by the vertices $\{\alpha \in \mathcal{C} : \ell(\alpha) \leq L\}$.

LEMMA 8.1. $\mathcal{C}(L_0)$ *is connected.*

*Remark.* Compare Theorem 1 of Bowditch [16], which has a similar conclusion but applies to the more general setting of punctured torus groups that are not necessarily discrete or faithful. Bowditch's proof is completely different, employing a beautiful algebraic/combinatorial method.

*Proof.* Consider $\alpha_0, \alpha_1 \in \mathcal{C}(L_0)$. Let $g_i : S \to N$ for $i = 0, 1$ be simplicial hyperbolic surfaces in the homotopy class of $\rho$, adapted to $\alpha_i$. In order to streamline the proof (and constants), we make a brief appeal to the technique of pleated surfaces. As mentioned earlier, if an infinite sequence of vertex slide moves around $\alpha_i$ is performed starting with $g_i$, the resulting family of induced metrics converges to a smooth hyperbolic metric $\sigma_i$ induced by a limiting map $g'_i$ (a pleated surface), which still takes $\alpha_i$ to its geodesic representative.

Let $\beta_i$ be a shortest curve in $(S, \sigma_i)$. Then $\ell(\beta_i) \leq \ell_{\sigma_i}(\beta_i) \leq L_0$, so that $\beta_i \in \mathcal{C}(L_0)$, and by Lemma 2.1 we have $\alpha_i \cdot \beta_i = 1$, so that it suffices to show that $\beta_0$ and $\beta_1$ can be connected by a path in $\mathcal{C}(L_0)$.

Applying Lemma 5.2, $g_0$ and $g_1$ can be connected by a continuous family $g_t$ of maps, for $t \in [0,1]$, such that each $g_t$ is a simplicial hyperbolic surface. Because of the way we constructed $g'_i$, we may extend this to a family (still called $g_t$) of simplicial hyperbolic surfaces interpolating between $g'_0$ and $g'_1$. Let $\sigma_t$ denote the (possibly singular) hyperbolic metric induced on $S$ by $g_t$.

Let $\tau_t$ denote the hyperbolic metric in the conformal class of $\sigma_t$. Then $\tau_t = \sigma_t$ for $t = 0, 1$ and $\tau_t \geq \sigma_t$ pointwise, by Ahlfors' lemma (see §2.3). Note also that the lengths $\ell_{\tau_t}$ of homotopy classes vary continuously with $t$. For each $t$ consider the curves of minimal $\tau_t$-length (necessarily no longer than $L_0$). If there are two such curves $\beta^1$ and $\beta^2$, then by Lemma 2.1 they are Farey neighbors.

Since $\{\tau_t\}$ is a compact family of metrics, the set $B$ of curves $\beta$ appearing as shortest curve for some $\tau_t$ is finite. If $X(\beta)$ is the closed set of $t$ values for which $\beta \in B$ is shortest in $\tau_t$ we have shown that the adjacency graph of $\{X(\beta)\}_{\beta \in B}$ is isomorphic to the subgraph of $\mathcal{C}$ spanned by $B$. Connectivity of the interval $[0, 1]$ now implies that $B$ is connected, and in particular there is a path between $\beta_0$ and $\beta_1$ of curves whose lengths in some $\tau_t$, and hence in $N$, are at most $L_0$. □

8.2. *Bounds on pivots.* As a consequence of the Connectivity Lemma 8.1, we may now obtain the simplest part of the Pivot Theorem:

LEMMA 8.2.　*There is a fixed $L_1$ such that $\ell(\alpha_n) \leq L_1$ for all pivots $\alpha_n$.*

*Proof.* Recalling the notation from Section 4, let $e$ be an edge of $E$, the set of edges separating $\alpha_-$ from $\alpha_+$. If $\nu_-$ is an irrational boundary point ($\nu_- = \alpha_-$), Lemma 5.1 implies there exists some $\gamma_- \in \mathcal{C}(L_0)$ on the same side of $e$ as $\nu_-$. If $\nu_-$ is not irrational then $\alpha_- \in \mathcal{C}$ is already in $\mathcal{C}(L_0)$; now let $\gamma_- = \alpha_-$. The same applies for $\nu_+$ and $\gamma_+$. Now $e$ separates $\gamma_-$ from $\gamma_+$, and by Lemma 8.1 there is a path between them in $\mathcal{C}(L_0)$. This implies that $e$ contains a vertex of $\mathcal{C}(L_0)$.

An internal pivot $\alpha_n$ by definition is the common point of at least two edges $e_0$ and $e_1$ of $E$ (barring the exceptional case (3) in Section 4), which we may choose to be edges of a common Farey triangle. Thus, if $\ell(\alpha_n) > L_0$, the other vertices $\beta_0$ and $\beta_1$ of $e_0, e_1$ have lengths bounded by $L_0$. Representing $\alpha_n$ and $\beta_0$ by a generator pair $A, B$, we know that $\beta_1$ can be represented by $A^{\pm 1}B$. Now the trace identity (7.2) gives a bound for $|\text{tr}A|$ in terms of $|\text{tr}B|$ and $|\text{tr}AB|$, and in particular bounds $\ell(A) = \ell(\alpha_n)$ by some $L_1$.

In the exceptional case where $E$ is a singleton $\{e\}$, we note that one vertex of $e$ is in $\mathcal{C}(L_0)$, and the other has as neighbors both $\alpha_-$ and $\alpha_+$ and so it is also in $\mathcal{C}(L_1)$ by the same argument. □

## 9. Controlled surfaces and building blocks

In this section we begin the construction of a geometric model for the quotient manifold of a punctured torus group. In Sections 9.2 and 9.3 we construct particular geometrically controlled maps of surfaces based only on some length bounds and end-invariant data. In Sections 9.4 and 9.5 we use these surfaces to construct "building blocks" which will later serve, in particular, to control the shape of Margulis tubes. Lemma 9.3 will summarize the basic properties of these blocks.

We remark that the cases of most interest here are those of "internal" blocks, and these are also the simpler cases. On a first reading one can suppress the discussion of boundary blocks (e.g. all of § 9.3), which are intended to describe the boundaries of the convex core. There are at most two boundary blocks, and typically an infinite number of internal blocks.

To streamline the discussion in this section and those that follow, we frequently use terms like "bounded" or "Lipschitz" to denote bounds which are uniform in the sense that they are independent of the particular punctured torus group.

9.1. *Standard metrics.* Let $\nu$ be a conformal structure on the punctured torus $S$, with marking $\alpha, \beta$ such that $\alpha$ is shortest in $\nu$. Then $\text{Im}\,(\tau(S, \nu, \alpha, \beta)) \geq \sqrt{3}/2$ by the discussion in Section 2.3. Let $\sigma^e$ be the Euclidean metric in the conformal class $\nu$ (incomplete at the puncture) such that $\alpha$ has length 1. We can cut $S$ into an open annulus $A$ and a closed punctured annulus $B$, with $\sigma^e$-geodesic boundaries homotopic to $\alpha$, so that the puncture is on the center circle of $B$ which has width $1/2$, while $A$ has width at least $(\sqrt{3}-1)/2$.

Let $\sigma^h$ be the complete hyperbolic metric in the conformal class $\nu$. By standard techniques of conformal mappings (see e.g. [26, Thm. 4.3]), we know that $\sigma^h$ and $\sigma^e$ are within a bounded ratio of each other in $B$ minus a disk of $\sigma^e$-radius $1/8$ (say) around the puncture. It will be useful to fix a hybrid metric $\sigma^m$ which equals $\sigma^e$ outside a disk of $\sigma^e$-radius $1/4$ around the puncture, equals $\sigma^h$ inside a disk of $\sigma^e$-radius $1/8$ around the puncture, and is within bounded ratio of both in the remaining annulus.

Let us also assume $\varepsilon_0$ is sufficiently small that the radius $1/8$ disk contains the $\varepsilon_0$-Margulis tube for the cusp of $\sigma^h$.

9.2. *Halfway surfaces.* We can apply Lemma 7.1 on cores to conclude that, for any two neighbors in $\mathcal{C}$ with bounded lengths, there is a simplicial hyperbolic surface that maps them both to curves of bounded length.

Given $\alpha, \beta$ which are neighbors in $\mathcal{C}$, consider the conformal structure $\nu_{\alpha,\beta}$ whose Teichmüller parameter is $\tau(S, \nu_{\alpha,\beta}, \alpha, \beta) = i$. Let $\sigma^e_{\alpha,\beta}$, $\sigma^h_{\alpha,\beta}$ and $\sigma^m_{\alpha,\beta}$ be the associated metrics as above. Note that $\sigma^e_{\alpha,\beta}$ makes $S$ a square torus in which both $\alpha$ and $\beta$ have length 1 and are orthogonal, and that $\sigma^h_{\alpha,\beta}$ and $\sigma^m_{\alpha,\beta}$ are within bounded ratio of each other. (The choice of a square torus is for convenience – any fixed modulus would do as well.)

The following lemma gives conditions for mapping this surface in a controlled way into $N$:

LEMMA 9.1. *Fix $L \geq L_0$ and let $\alpha, \beta$ be a pair of Farey neighbors contained in $\mathcal{C}(L)$. Then there exists a simplicial hyperbolic surface $f_{\alpha,\beta} : S \to N$ in the homotopy class of $\rho$, which is $L'$-Lipschitz with respect to the metrics $\sigma^m_{\alpha,\beta}$ and $\sigma^h_{\alpha,\beta}$.*

*In particular, $f_{\alpha,\beta}(S)$ does not meet any $\varepsilon_2$-Margulis tube other than that of the main cusp.*

*The constants $L', \varepsilon_2$ depend only on $L$.*

*Proof.* Apply the notation given in Section 7 for the lift to $\mathbf{H}^3$ and the group elements $\mathbf{A}$ and $\mathbf{B}$ corresponding to $\alpha$ and $\beta$. Let $\mathcal{T}$ denote the triangulation whose closed curves are $\alpha$ and $\beta$. Let $f_1$ and $f_2$ denote simplicial hyperbolic surfaces with triangulation $\mathcal{T}$, adapted to $\alpha$ and $\beta$ respectively. The type (3) elementary move from $f_1$ to $f_2$ involves a step in which the vertex of the triangulation (or rather its lift to $\mathbf{H}^3$) travels from $T_\mathbf{A}$ to $T_\mathbf{B}$ along their common perpendicular. By Lemma 7.1, there is an intermediate map $f_t$ for which the vertex is in the intersection of $\mathrm{core}_{D_2}(\mathbf{A})$ and $\mathrm{core}_{D_2}(\mathbf{B})$. At this point both $\mathbf{A}$ and $\mathbf{B}$ have translation distance bounded by $L + D_2$, and hence the images of the corresponding triangulation edges have this length bound.

If we realize $\mathcal{T}$ so that its edges are geodesic with respect to $\sigma^h_{\alpha,\beta}$, and then parametrize $f_t$ with constant speed on the edges, we find that it is Lipschitz on the 1-skeleton (this can be made to work on the edges going to the cusp, as well). The extension to the whole surface can easily be made (uniformly) Lipschitz with respect to $\sigma^h_{\alpha,\beta}$. Thus we let this be the definition of $f_{\alpha,\beta}$.

Since in particular there is a bound on the diameter of the noncuspidal part of $f_{\alpha,\beta}(S)$, it follows from Lemma 6.1 that there is a corresponding $\varepsilon_2$ so that the map avoids all $\varepsilon_2$-Margulis tubes other than the main cusp. □

9.3. *Boundary surfaces.* In order to control the shape of Margulis tubes that are not contained in the convex core, we will need to construct a version of the controlled maps of the previous section, associated to boundary components of the convex core. In what follows we consider the $+$ end, but of course the same considerations apply to both ends.

Suppose that $e_+$ is geometrically finite and let $S_+$ denote the surface at infinity $\Omega_+/\rho(\pi_1(S))$. Unless $\alpha_+$ is parabolic, $S_+$ is naturally identified with $S$ and has conformal structure $\nu_+$. Noting that $\alpha_+$ is shortest in $\nu_+$, let $\sigma^h_+$, $\sigma^e_+$ and $\sigma^m_+$ denote the associated metrics as in Section 9.1. Let $A_+$ and $B_+$ be the annulus and punctured annulus described in that section, with $\sigma^e$-geodesic boundaries isotopic to $\alpha_+$.

If $\alpha_+$ is parabolic then $S_+$ is a thrice-punctured sphere, $\sigma^e_+$ makes $S_+$ a biinfinite Euclidean annulus with puncture and $A_+$ is a union of two semi-infinite annuli. The other parts of the definition are unchanged.

Let $C_r(N)$ denote the $r$-neighborhood of the convex core of $N$. Fixing a positive value for $r$, let $\widehat{C}(N)$ denote the union of $C_r(N)$ with all $\varepsilon_0$ Margulis tubes other than that of the main cusp. (Note that there are at most two such tubes that are not already contained in $C(N)$, corresponding to $\alpha_+$ and $\alpha_-$ if these are sufficiently short.)

Let $\partial_+\widehat{C}(N)$ and $\partial_+C_r(N)$ denote the components of $\partial\widehat{C}(N)$ and $\partial C_r(N)$, respectively, that face the end $e_+$. Note that either $\partial_+\widehat{C}(N) = \partial_+C_r(N)$, or $\partial_+\widehat{C}(N) = \mathcal{U}_+ \cup \mathcal{R}_+$, where $\mathcal{U}_+$ is an annulus in the Margulis tube boundary $\partial \mathbf{T}_{\varepsilon_0}(\alpha_+)$, and $\mathcal{R}_+$ is a punctured annulus isotopic to the complement of $\alpha_+$ in $S$.

The following lemma gives the properties of the boundary maps we shall build:

LEMMA 9.2. *If $e_+$ is geometrically finite, there exists a homeomorphism*

$$\Pi_+ : S_+ \to \partial_+\widehat{C}(N),$$

*in the homotopy class determined by $\rho$, such that: $\Pi_+$ is $K$-bilipschitz with respect to $\sigma_+^m$ on $S_+$ and the induced metric from $\mathbf{H}^3$ on $\partial_+\widehat{C}(N)$. It takes $A_+$ to $\mathcal{U}_+$ by a $K$-quasiconformal homeomorphism, where the constant $K$ is independent of $\rho$. The corresponding statement holds for $e_-$.*

*Proof.* We construct $\Pi_+$ as follows. Recalling that we have fixed some small $r$, let $\Pi_r : S_+ \to \partial_+C_r(N)$ be the orthogonal projection from infinity (see Epstein-Marden [30]). For any $x \in S_+$ there is a geodesic ray $m_x$ based at $y = \Pi_r(x)$ and terminating in $x$ at infinity (it may be helpful to consider this to be happening in the universal cover). This ray meets $\mathbf{T}_{\varepsilon_0}(\alpha_+)$ if and only if $y \notin \text{int}(\mathcal{R}_+)$, and in such a case let $z$ be the (unique) point of intersection of $m_x$ with $\mathcal{U}_+ \subset \partial \mathbf{T}_{\varepsilon_0}(\alpha_+)$. We define $\widehat{\Pi}_r(x) = z$ in this case, and $\widehat{\Pi}_r(x) = \Pi_r(x)$ if $\Pi_r(x)$ lies in $\mathcal{R}_+$.

By Theorem 2.3.1 of [30], $\Pi_r$ is a quasiconformal and bilipschitz homeomorphism between $(S_+, \sigma_+^h)$ and $\partial_+C_r(N)$, where the constants depend only on the choice of $r$. It is also easy to see that on any point $z = \widehat{\Pi}_r(x)$ of $\mathcal{U}_+$ which is a distance $r$ or more from $\partial_+C_r(N)$, the tangent plane to $\mathcal{U}_+$ at $z$ makes an angle with the ray $m_x$ which is bounded away from 0 (this follows from the fact that $\partial \mathbf{T}_{\varepsilon_0}(\alpha_+)$ is an equidistant surface from the core geodesic of $\alpha_+$, which lies inside $C(N)$). Thus, by an elementary calculation, the map $\widehat{\Pi}_r$ at $x$ is quasiconformal with a bound again depending on $r$.

On the other hand the part of $\mathcal{U}_+$ which is within distance $r$ consists of two annuli of bounded radius adjacent to $\partial \mathcal{U}_+$, and they project to annuli $Y_1, Y_2$ in $\partial_+C_r(N)$ by a distance-decreasing map.

Since $\Pi_r$ is bilipschitz with respect to $\sigma_+^h$, and the surfaces $\mathcal{R}_+$ and $Y_1, Y_2$ are of standard shape ($\mathcal{R}_+$ is bilipschitz equivalent to the complement of a Margulis tube in a punctured torus, and $Y_i$ are bilipschitz equivalent to Eu-clidean annuli of fixed height and girth), there is a bilipschitz homeomorphism $k : S_+ \to S_+$ isotopic to the identity so that $\Pi_r^{-1}(\mathcal{R}_+) = k(B_+)$ and $\Pi_r^{-1}(Y_j) = k(X_j)$ where $X_j$ are Euclidean annuli in $A_+$ of bounded modulus ($j = 1, 2$). Then, letting $\Pi' = \widehat{\Pi}_r \circ k$, we have a map which is bilipschitz as desired from $B_+$ to $\mathcal{R}_+$, and from $A_+$ to $\mathcal{U}_+$ is Lipschitz on bounded Euclidean neighborhoods of the boundary, and quasiconformal on the rest. Lemma 2.2 implies that $\Pi'|_{A_+}$ is homotopic to a quasiconformal map which has constants bounded in terms of the previous choices, and which is bilipschitz with respect to the Euclidean metrics on the annuli. We let $\Pi_+$ be the resulting map. □

9.4. *Building blocks.* A building block $\mathcal{B}$ is a copy of $S \times [0,1]$, with a certain metric defined on a subset of $\mathcal{B}$. Blocks are associated to local configurations of pivots in $\mathcal{C}(L)$, and come in several varieties (there are also blocks with parabolics, which are topologically different).

*Internal blocks.* Let $\alpha, \beta, \gamma$ be distinct vertices of $\mathcal{C}(L)$ such that $\beta$ and $\gamma$ are neighbors of $\alpha$. We define the associated *internal block* $\mathcal{B} = \mathcal{B}_{\alpha,\beta,\gamma}$ as follows:

The definitions of two metrics $\sigma_0 = \sigma_{\alpha,\beta}^m$ and $\sigma_1 = \sigma_{\alpha,\gamma}^m$ are as in Sections 9.1 and 9.2. After isotopy we may assume that the open annulus $A$ from Section 9.1 is the same for both, and that the two are equal on the complementary annulus $B$. Let $A'$ be an open $1/8$-neighborhood of $A$ in both metrics.

Within $\mathcal{B}$ let $U = A \times (1/4, 3/4)$ and $U' = A' \times (0, 1)$ be nested, open solid tori, and let $\mathcal{B}^0 = \mathcal{B} \setminus U$. Now let $\partial_0 \mathcal{B}$ and $\partial_1 \mathcal{B}$ denote $S \times \{0\}$ and $S \times \{1\}$, respectively.

Place the metrics $\sigma_i$ on the boundaries $S \times \{i\}$, respectively, and extend to a metric on $\mathcal{B}^0$ as follows: let $S_0 = S \setminus Q$ where $Q$ is the $\varepsilon_0$-cusp of either metric, and $B_0 = B \setminus Q$. Let $Q_\mathcal{B} = Q \times [0,1]$.

On $S_0 \times [0, 1/4]$ we place the product metric $\sigma_0 \times dt^2$ where $t \in [0, 1/4]$, and on $S_0 \times [3/4, 1]$ we similarly place $\sigma_1 \times dt^2$. On $B_0 \times [1/4, 3/4]$, where both metrics are equal, we can place a product metric in the same way, using either metric. On $Q_\mathcal{B}$ we may place the geometry of a three-dimensional parabolic cusp: namely, identify $Q_\mathcal{B} = Q \times [0,1]$ with a segment of the $\varepsilon_0$-Margulis tube of a rank-1 parabolic group, so that $\partial Q \times [0,1]$ is an annular section of length 1 of the boundary, and each slice $Q \times \{t\}$ is a totally geodesic 2-dimensional Margulis tube embedded orthogonally to the boundary.

*Boundary blocks.* Suppose $e_+$ is geometrically finite, let $\alpha_+ \in \mathcal{C}$ be shortest in $\nu_+$, as in Section 4, and let $\beta$ be any neighbor of $\alpha_+$ which is also in $\mathcal{C}(L)$. To this pair we associate a *boundary block* $\mathcal{B} = \mathcal{B}_{\alpha_+,\beta}$ as follows.

Suppose first that $\alpha_+$ is not parabolic. Let the metric $\sigma_0$ be $\sigma_{\alpha_+,\beta}^m$, just as before. The metric $\sigma_1$ is $\sigma_+^m$ from Section 9.3. Again after isotopy we may assume that the annuli $A$ and $B$ of Section 9.1 are the same for both metrics, as well as the cusps $Q$, and the metrics are equal on $B$ except possibly on a small collar around the cusp, where they are in bounded ratio. Define $A'$ so that it excludes the cusp, and contains a $1/8$ neighborhood of $A$ in each metric (it can't be exactly the $1/8$ neighborhood on both because they are not quite equal), and define the solid torus $U'$ as before. However, we let $U$ be $A \times (1/4, 1)$. The product metric on $\mathcal{B} \setminus U$ is defined as before, except on the neighborhood of the puncture where we interpolate between the two slightly different metrics.

If $N$ has a parabolic cusp in its $e_+$ end (so $\nu_+ = \alpha_+ \in \mathcal{C}$) then we *redefine* the topological structure of $\mathcal{B}_{\alpha_+,\beta}$ to be $(S \times [0,1]) \setminus (a \times \{1\})$, where $a$ is a curve representing $\alpha$, and note that in this case $\sigma_+^m$ on the punctured cylinder $S_+ \cong S \setminus a$ has infinite area. The rest of the construction is the same. In this case we say that $\mathcal{B}_{\alpha_+,\beta}$ is a *boundary block with parabolic*.

For the opposite end $e_-$ we construct $\mathcal{B}_{\alpha_-,\beta}$ in just the same way, except that the roles of 0 and 1 are reversed, and $U$ for example is $A \times (0, 3/4)$.

Finally, suppose $\alpha_+ = \alpha_-$, that is, the same curve is short on both boundaries. We define a *double boundary block* $\mathcal{B} \equiv \mathcal{B}_{\alpha_\pm}$ as follows: Let $\sigma_0$ and $\sigma_1$ be $\sigma_-^m$ and $\sigma_+^m$ respectively, again with a common decomposition into annuli $A$ and $B$. Now we define $U = A \times (0, 1)$ and continue as before.

We can again have a block with parabolics in this case, but note that $\alpha_+ = \alpha_-$ can be parabolic at most on one of the two ends.

9.5. *Block maps and the figure-8 argument*. The following lemma gives us controlled maps of our blocks into the manifold $N$.

LEMMA 9.3. *Let $\mathcal{B}$ denote a block constructed as in the previous section. That is, $\mathcal{B} = \mathcal{B}_{\alpha,\beta,\gamma}$ for $\alpha, \beta, \gamma$ in $\mathcal{C}(L)$; or $\mathcal{B} = \mathcal{B}_{\alpha_+,\beta}, \mathcal{B}_{\alpha_-,\beta}$ or $\mathcal{B}_{\alpha_\pm}$. In the latter cases denote $\alpha = \alpha_+$ or $\alpha = \alpha_-$.*

*There is a map $H_\mathcal{B} : \mathcal{B} \to \widehat{C}(N)$, in the homotopy class determined by $\rho$, with the following properties:*

1. $H_\mathcal{B}$ *is Lipschitz on $\mathcal{B}^0$, with uniform constant $K$.*

2. $H_\mathcal{B}$ *respects the Margulis tube structure*:

   *If $\ell(\alpha) \le \varepsilon_3$ then $H_\mathcal{B}$ maps $(U, \partial U)$ to $(\mathbf{T}_{\varepsilon_0}(\alpha), \partial(\mathbf{T}_{\varepsilon_0}(\alpha)))$.*

   *Furthermore, $H_\mathcal{B}(\mathcal{B}^0)$ is outside every $\varepsilon_0$-Margulis tube whose core length is no more than $\varepsilon_3$, except for the main cusp.*

3. *For boundary blocks, $H_\mathcal{B}$ is boundary-preserving*:

   *If $\alpha = \alpha_+$ (respectively $\alpha = \alpha_-$) then $H_\mathcal{B}$ takes $\partial_1 \mathcal{B}$ (resp. $\partial_0 \mathcal{B}$) to $\partial_+ \widehat{C}(N)$ (resp. $\partial_- \widehat{C}(N)$) by a $K$-bilipschitz homeomorphism of degree 1.*

4. $H_\mathcal{B}$ *respects the main cusp structure*:

> *It takes $(Q_\mathcal{B}, \partial Q_\mathcal{B})$ to $(Q_N, \partial Q_N)$ by a proper map, and takes $\mathcal{B} \setminus Q_\mathcal{B}$ to $N \setminus Q_N$, where $Q_N$ is the $\varepsilon_0$ Margulis tube of the main cusp of $N$.*

*The constants $K$, $\delta$ and $\varepsilon_3$ depend only on $L$.*

*Proof.* We begin by defining $H = H_\mathcal{B}$ on $\partial \mathcal{B}$. Take the case of an internal block. Then by construction the boundary metrics $\sigma_0$ and $\sigma_1$ are $\sigma^m_{\alpha,\beta}$ and $\sigma^m_{\alpha,\gamma}$ respectively, and Lemma 9.1 gives us homotopic maps $h_0 = f_{\alpha,\beta}$ and $h_1 = f_{\alpha,\gamma}$ that are Lipschitz with respect to $\sigma_0$ and $\sigma_1$ respectively, and which avoid all $\varepsilon_2$-Margulis tubes except for the main cusp. It is easy to show (see e.g. [70, Lemma 4.1]) that, since $\pi_1(B)$ is nonabelian and the maps are $\pi_1$-injective, there is a *unique* normalized geodesic homotopy $H_1$ between $h_0|_B$ and $h_1|_B$, that is, a homotopy whose tracks $H_1|_{x \times [0,1]}$, $x \in B$, are geodesics parametrized at constant speed.

We can define a preliminary extension $H_2 : \mathcal{B} \setminus U' \to N$ to be the unique map which restricts to $h_i$ on $S \times \{i\}$ and which restricts on $B' \times [0,1]$ to the geodesic homotopy $H_1$. (Here $B' = S \setminus A'$).

We can do the same thing for a noninternal block. Suppose for example we have $\mathcal{B} = \mathcal{B}_{\alpha_+,\beta}$. The map $h_0$ is defined to be $f_{\alpha,\beta}$ as before (with $\alpha = \alpha_+$), but $h_1$ is taken to be the map $\Pi_+ : S \to \partial_+ \widehat{C}(N)$, constructed in Section 9.3. The case with parabolic is handled similarly.

If $\alpha = \alpha_+ = \alpha_-$ and $\mathcal{B} = \mathcal{B}_{\alpha_\pm}$, we define $h_0 = \Pi_-$ and $h_1 = \Pi_+$ and continue as above. In each case we extend $h_0$ and $h_1$ to the rest of $B' \times [0,1]$ with the geodesic homotopy.

We now prove that $H_2$ is uniformly Lipschitz on $\mathcal{B} \setminus U'$:

Let $\eta$ be a "figure-8" curve in $B'$, as in Figure 9. The unoriented free homotopy class of $\eta$ only depends on the curve $\alpha$, and there is a uniform bound on the length of the geodesic representative of $\eta$ in either $\sigma_0$ or $\sigma_1$. Since $\sigma_0$ and $\sigma_1$ are within a uniformly bounded ratio in $B'$, we may choose $\eta$ to have length bounded by some $\ell_1$ in both.

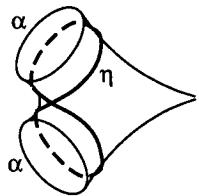

Figure 9. The figure-8 curve in $S \setminus \alpha$

Think of $\eta$ as the union of two simple loops $\eta_1$ and $\eta_2$ with a common point $y$. These determine distinct elements of $\pi_1(S, y)$ in the conjugacy class of $\alpha$. Choosing a lift of $H_2(\{y\} \times \{0\})$ to $\mathbf{H}^3$, we represent $\eta_1$ and $\eta_2$ as distinct isometries $\mathbf{A}_1, \mathbf{A}_2$. The endpoints of the lift of the geodesic $H_2(\{y\} \times [0,1])$ give us two points $p, q \in \mathbf{H}^3$ where the translation distances of both $\mathbf{A}_1$ and $\mathbf{A}_2$ are bounded above by $\ell_1$. Bounding the length of the homotopy now reduces to bounding $\mathrm{dist}(p, q)$.

Suppose $\mathrm{dist}(p, q) \geq D$ for some $D$. Let $\varepsilon(D)$ denote a function of the form $ce^{-D}$ where $c$ is some uniform constant. The geodesic $\overline{pq}$ has image $\mathbf{A}_i(\overline{pq})$ ($i = 1, 2$) whose endpoints are at most $\ell_1$ from those of $\overline{pq}$. It follows by a bit of hyperbolic trigonometry that if, say, $D > 3\ell_1$, then $\mathbf{A}_i(\overline{pq})$ comes within $\varepsilon(D)$ of the midpoint $m$ of $\overline{pq}$, and the two geodesics have directions differing by at most $\varepsilon(D)$ at their closest point, as measured in the tangent bundle to $\mathbf{H}^3$. In particular $\mathbf{A}_i(m)$ is $\varepsilon(D)$-close to $\overline{pq}$. It follows that $[\mathbf{A}_1, \mathbf{A}_2]$, which is nontrivial, translates $m$ by no more than $\varepsilon(D)$. However since $\mathbf{A}_1$ and $\mathbf{A}_2$ translate $m$ no more than $\ell_1$, when $D$ is sufficiently large this violates Jørgensen's inequality, and this gives us our upper bound. (See [70, Lemma 4.2] for a similar argument in a somewhat different setting.)

This bounds the length of the track $H_2(\{y\} \times [0,1])$ and because of the bound on the length of $\eta$, we obtain a bound for the other tracks on $\eta$ as well. Now since $\eta$ is a deformation-retract of $B$ and the thick part $B' \setminus Q$ of $B'$ in either metric has uniformly bounded diameter, we obtain a bound on the track length of $H_2$ in $B' \setminus Q$. Since $H_2$ is parametrized at constant speed, this gives the desired Lipschitz bound on $H_2$ in $(B' \setminus Q) \times [0,1]$. In the cusp region $Q$, a bound is immediate from the fact that $h_0$ and $h_1$ are standard embeddings.

We next extend $H_2$ to the rest of $\mathcal{B}$.

Consider first the internal case. Let $V'$ be the boundary torus of $U'$. Consider the solid-torus covering space $\widetilde{N}$ of $N$ whose geodesic core is a homeomorphic lift of $\alpha^*$. The map $H_2|_{V'}$ lifts to $\widetilde{H}_2 : V' \to \widetilde{N}$, and its image lies in a bounded neighborhood of the boundary torus $F_0$ of the $\varepsilon_0$-Margulis tube of $\widetilde{N}$ (a homeomorphic lift of $\mathbf{T}_{\varepsilon_0}(\alpha)$).

If $\ell(\alpha)$ is at most some $\varepsilon$ given by Lemma 6.1, the radius of the tube is at least $1 + \mathrm{diam}(\widetilde{H}_2(V'))$, and hence the image is disjoint from a 1-neighborhood of the core of $\widetilde{N}$. Thus the radial projection taking $\widetilde{N}$ minus its core to the torus $F_0$ is uniformly Lipschitz on $\widetilde{H}_2(V')$.

We can therefore use this projection, parametrized as a deformation retraction, to extend $\widetilde{H}_2$ to a collar $V' \times [0,1]$ so that $\widetilde{H}_2$ takes $V' \times \{0\}$ to $F_0$, and so that the map is Lipschitz with uniform constant. Projecting downstairs to $N$ we obtain an extension of $H_2$ to the collar $U' \setminus U$, which is still Lipschitz with some uniform constant since the collar is bilipschitz equivalent to $V' \times [0, 1)$. The map is now defined on $\partial U$ and extends to map the rest of $U$ into $\mathbf{T}_{\varepsilon_0}(\alpha)$ because it is in the right homotopy class. (Note that there is no Lipschitz condition to be satisfied on $U$, which in fact does not yet have a metric). Let $H_3$ denote the resulting map.

If $\ell(\alpha) \geq \varepsilon$, then the lifted map lies in a bounded neighborhood of the axis. Radial projection to the axis is then uniformly Lipschitz and we use this to extend the map to the the collar so that $\partial U$ maps to the axis. Again, the map extends to some $H_3$ on the rest of $U$ for topological reasons.

The extension for boundary blocks works similarly – we lift to a cover in which $\mathbf{T}_{\varepsilon_0}(\alpha)$, which is still a solid torus, lifts homeomorphically, and apply radial projection as before to show that the map on the subannuli of $V'$ which do not map to $\partial \widehat{C}(N)$ can be pushed rel boundary to the boundary of the Margulis tube. We use this homotopy to define the map on $U' - U$ so that $\partial U$ maps as required. It is easy to see that this can be made Lipschitz with uniform constants.

By construction, for some uniform $\varepsilon_3$, $H(\mathcal{B} \setminus U)$ stays out of the $2\varepsilon_3$-thin part (minus the main cusp). This is because $\mathcal{B} \setminus (U \cup Q_\mathcal{B})$ (the thick part of the block) has bounded diameter and hence can penetrate only a bounded amount into any thin part. Thus for any curve $\gamma$ of length less than $\varepsilon_3$, the region $\mathbf{T}_{\varepsilon_0}(\gamma) \setminus \mathbf{T}_{2\varepsilon_3}(\gamma)$ retracts onto $\partial \mathbf{T}_{\varepsilon_0}(\gamma)$ by radial projection, which is uniformly Lipschitz. We may compose $H_3|_{\mathcal{B} \setminus U}$ with this retraction to obtain a new map $H_4$ which stays out of this Margulis tube, and still satisfies uniform Lipschitz bounds. Thus we have Part (2).

A similar argument using the radial projection in the main cusp gives us an additional Lipschitz deformation of $H_4$ to get $H_5$ which satisfies Part (4). This will be our final map $H_\mathcal{B}$.

Part (3) of the lemma follows from Lemma 9.2, with the statement about degree following from our orientation conventions applied to $\mathcal{B}$ and to $N$. □

It is worth comment at this point that the restriction $H|_{\partial U} : \partial U \to \partial \mathbf{T}_{\varepsilon_0}(\alpha)$ has *not* been shown to be homotopic to a homeomorphism. In fact one can construct blocks for which this map has degree different from 1. This will be remedied in the next section by a global argument.

## 10. Proof of the pivot theorem

We begin by piecing together a selected sequence of the blocks from the previous section.

Recall from Section 4 that the pivot sequence $P$ is written $\{\alpha_n\}_{n \in J}$ where $J$ is a finite, half-infinite or bi-infinite interval in $\mathbf{Z}$. The subset $P_0$ of internal pivots leaves out $\alpha_- = \alpha_0$ and/or $\alpha_+ = \alpha_{p+1}$, if these occur in $P$.

Hence for each internal pivot $\alpha_n$ the pivots $\alpha_{n-1}$ and $\alpha_{n+1}$ are defined. Since, by Lemma 8.2, all of $P$ is contained in $\mathcal{C}(L_1)$, we may construct a block $\mathcal{B}_n = \mathcal{B}_{\alpha_n, \alpha_{n-1}, \alpha_{n+1}}$ and associated map $H_n : \mathcal{B}_n \to \widehat{C}(N)$, as in Sections 9.4 and 9.5.

Note that $L_1$ determines values for the constants $\varepsilon_2, \varepsilon_3$, as described in the previous section. Furthermore, a small positive constant $r$ was chosen in Section 9.5, which determines the definition of $\widehat{C}(N)$.

If $\alpha_-$ is in $P$ and different from $\alpha_+$ we let $\mathcal{B}_0 = \mathcal{B}_{\alpha_-, \alpha_1}$, and if $\alpha_+$ is in $P$ and different from $\alpha_-$ we let $\mathcal{B}_{p+1} = \mathcal{B}_{\alpha_+, \alpha_p}$. If $\alpha_- = \alpha_+$ (hence $P_0$ is empty) we have just one (double boundary) block which we label $\mathcal{B}_0 = \mathcal{B}_{\alpha_\pm}$. Again we define the maps $H_n$ appropriately, as in Section 9.5.

Next we observe that for consecutive $n, n+1 \in J$, there is a natural identification between $\partial_1 \mathcal{B}_n$ and $\partial_0 \mathcal{B}_{n+1}$ – namely, the identity on $S$ induces an isometry of the metrics on these surfaces. Thus we can glue together all the blocks in sequence to obtain

$$M = \bigcup_{n \in J} \mathcal{B}_n$$

with a map $H : M \to \widehat{C}(N)$, restricting to $H_n$ on each $\mathcal{B}_n$. On the subset $M^0 = M \setminus \cup U_n$, the map $H$ is Lipschitz. (In fact on the complement of $M^0$ we have not yet defined a metric.)

Note that in the case of two degenerate ends the sequence is bi-infinite and $M$ is homeomorphic to $S \times \mathbf{R}$. In other cases it will have one or two boundary components, and in all cases $M$ is homeomorphic to $\widehat{C}(N)$.

Let us also denote by $Q_M$ the "main cusp" of $M$, that is, the union of cusps $Q_{\mathcal{B}_n}$ over all blocks $\mathcal{B}_n$.

LEMMA 10.1. *The map $H : M \to \widehat{C}(N)$ is proper and has degree 1.*

*Proof.* $H$ takes boundary to boundary, the image of the thick part of each block has bounded diameter, and the pivot geodesics $\alpha_n^*$ leave every compact set as $|n| \to \infty$. It follows that $H$ restricted to $M \setminus Q_M$ is proper to $\widehat{C}(N) \setminus Q_N$. On each cusp block $Q_\mathcal{B}$, $H$ is already proper and in fact by construction the function $x \mapsto \mathrm{dist}(H(x), \partial Q_N)$ is proper on $Q_\mathcal{B}$. It follows that $H$ is proper as a map from $M$ to $\widehat{C}(N)$. Therefore the degree of $H$ is well-defined.

If $\widehat{C}(N)$ has nonempty boundary then, by the last part of Lemma 9.3, $H$ has degree 1 on the boundary and we are done. If not, then $M = S \times R$, and we know that $H$ acts as a homotopy equivalence in the homotopy class of an orientation preserving homeomorphism on $S$. Thus it has degree 1 if the order of the ends is preserved, degree $-1$ if it is reversed, and degree 0 if both ends of $M$ are mapped to the same end of $N$. Because the pivot geodesics $\alpha_n^*$ exit the end $e_+$ when $n \to \infty$ and $e_-$ when $n \to -\infty$, our orientation convention applied to both $M$ and $N$ implies that $H$ has degree 1. □

With this map in place we can complete the proof of the Pivot Theorem 4.1.

Let us show that only the pivots can have length less than $\varepsilon_3$: Since the map $H$ has nonzero degree, in particular it covers all of $\widehat{C}(N)$ and hence meets the Margulis tube of every short geodesic. However, by the construction it can only meet the $\varepsilon_3$-Margulis tubes of the pivots. This proves part (1) of Theorem 4.1

We now consider the proof of part (2) of Theorem 4.1, which we recall asserts

(10.1) $$\omega(\alpha_n) \approx \nu_+(\alpha_n) - \overline{\nu}_-(\alpha_n) + i$$

for any pivot $\alpha_n$, where $\omega(\alpha_n) = 2\pi i/\lambda(\alpha_n)$.

Let $\alpha_n$ be a pivot of length greater than $\varepsilon_3$ (and no greater than $L_1$). Then $\lambda(\alpha_n)$ is bounded hyperbolic distance from 1 in the right half-plane, so $\omega(\alpha_n)$ is bounded hyperbolic distance from $2\pi i$ in the upper half plane $\mathbf{H}^2$. Thus, to establish (10.1) it suffices to find an upper bound for $\mathrm{Im}\,(\nu_+(\alpha_n)) + \mathrm{Im}\,(\nu_-(\alpha_n))$, and for $|\mathrm{Re}\,(\nu_+(\alpha_n)) - \mathrm{Re}\,(\nu_-(\alpha_n))|$. The former comes directly from the variant of Bers' inequality, given in part (6.4) of Lemma 6.4.

The difference of real parts is just $w(n)$ up to an error of at most 2 by (4.2), where if $\alpha_n$ is internal then $w(n)$ is the number such that $\alpha_{n+1} = D_{\alpha_n}^{w(n)} \alpha_{n-1}$. Representing $\alpha_n, \alpha_{n-1}$ by a generating pair $A, B$, we may represent $\alpha_{n+1}$ by $A^{\pm w(n)} B$. An easy computation shows that $\mathrm{tr} A^k B$ is of the form $\mathrm{tr}_k = ae^{k\lambda/2} + de^{-k\lambda/2}$, where $\lambda = \lambda(A)$. Since we have an upper bound on length (hence trace) of all pivots, we have an upper bound on $\mathrm{tr}_0$ and $\mathrm{tr}_{w(n)}$. On the other hand we have a positive lower bound $|\mathrm{tr}_k| \geq c_2 > 0$ by (6.5), and a lower bound of $\varepsilon_3$ on $\mathrm{Re}\,\lambda$, which gives $|\mathrm{tr}_k|$ definite exponential growth as a function of $|k|$. This implies an upper bound on $|w(n)|$.

If $\alpha_n$ is first or last in $P$, suppose without loss of generality it is first. If it has length at least $\varepsilon_3$, then in fact there is a second curve of bounded length on $\partial_- C(N)$ besides $\alpha_n$, and this plays the role of $\alpha_{n-1}$. The rest of the computation goes through in the same way.

Now we come to the heart of the matter: consider a pivot $\alpha_n$ with length at most $\varepsilon_3$. In this case our strategy will be to estimate the Teichmüller parameter $\omega(F_0)$, where $F_0$ is the boundary torus of $\mathbf{T} \equiv \mathbf{T}_{\varepsilon_0}(\alpha_n)$. Recall from Section 6 that $\omega(F_0)$ is within bounded distance of $\omega(\alpha_n)$.

We estimate $\omega(F_0)$ by comparing it to the Teichmüller parameter for the torus $\partial U_n$ in the model manifold $M$, which we can compute explicitly from the end invariants. Thus the argument is completed by showing that a quasiconformal map exists from one torus to the other, in the homotopy class preserving the markings.

Let us then fix a block $\mathcal{B} = \mathcal{B}_n$, and compute a Teichmüller parameter for the torus $\partial U = \partial U_n$.

We subdivide $\partial U$ into a union of four annuli $A_0 \cup A_R \cup A_1 \cup A_L$, as follows: Let $t_0 = 1/4$ if $\mathcal{B}$ is internal or of type $\mathcal{B}_{\alpha_+,\beta}$; in other cases let $t_0 = 0$. Let $t_1 = 3/4$ if $\mathcal{B}$ is internal or of type $\mathcal{B}_{\alpha_-,\beta}$; in other cases let $t_1 = 1$. With this convention, $A_0 = A \times \{t_0\}$, $A_1 = A \times \{t_1\}$, and $A_L$ and $A_R$ are the components of $\partial A \times [t_0, t_1]$.

We call $A_0, A_1$ the bottom and top respectively, and $A_L, A_R$ the sides of $\partial U$. Note that each annulus is Euclidean with geodesic boundary, so the union $\partial U$ is a Euclidean torus (which is PL-embedded in $\mathcal{B}$).

Next, we construct an explicit curve $\mu$ which represents the meridian of $\partial U$. Fix a neighbor $\beta \in \mathcal{C}$ of $\alpha = \alpha_n$. Realize $\beta$ as a curve $\beta_0$ which intersects $B = S \setminus A$ in a geodesic arc orthogonal to the boundary with respect to the Euclidean metric $\sigma_0^e$, and intersects $A$ in a $\sigma_0^e$-geodesic (recall from §9.1 that $\sigma_0^e$ and $\sigma_0$ are conformal, and agree except on a small disk around the puncture). Define $\beta_1$ similarly with respect to $\sigma_1^e$. Since $\sigma_0^e$ and $\sigma_1^e$ are equal on $B$, we may choose the arcs $\beta_0 \cap B$ and $\beta_1 \cap B$ to be equal. Join the endpoints of $(\beta_0 \cap B) \times \{t_0\}$ to those of $(\beta_1 \cap B) \times \{t_1\}$ by vertical arcs $\mu_L, \mu_R$ in $A_L, A_R$ respectively, of length $t_1 - t_0$. Let $\mu_0 = (\beta_0 \cap A) \times \{t_0\}$, and $\mu_1 = (\beta_1 \cap A) \times \{t_1\}$. The closed curve $\mu = \mu_0 \cup \mu_R \cup \mu_1 \cup \mu_L$ is then null-homotopic in $U$.

By (2.3), the Teichmüller parameter of $\partial U$ can be written

(10.2) $$\tau(\partial U, \alpha, \mu) = \tau(A_0, \mu_0) + \tau(A_L, \mu_L) + \tau(A_1, \mu_1) + \tau(A_R, \mu_R).$$

By construction we have

$$\tau(A_L, \mu_L) = \tau(A_R, \mu_R) = (t_1 - t_0)i$$

which we note is between $i$ and $i/2$.

For $A_0$ and $A_1$, we claim

(10.3) $$\begin{aligned} \tau(A_1, \mu_1) &= \tau(S, \sigma_1, \alpha, \beta) - i/2, \\ \tau(A_0, \mu_0) &= -\overline{\tau}(S, \sigma_0, \alpha, \beta) - i/2. \end{aligned}$$

For the first of these, (2.3) implies that $\tau(S, \sigma_1, \alpha, \beta) = \tau(A, \mu_1) + \tau(B, \beta_1 \cap B)$, and by construction $\tau(B, \beta_1 \cap B) = i/2$.

For the second, the proof is the same, except that we note the orientations of $A_0$ as a subannulus of $S$ and as a subannulus of $\partial U$ are opposite, and hence $\tau$ is replaced with $-\overline{\tau}$.

We next claim that

(10.4) $$\tau(S, \sigma_1, \alpha, \beta) - \overline{\tau}(S, \sigma_0, \alpha, \beta) \approx \nu_+(\alpha) - \overline{\nu}_-(\alpha) + i,$$

where "$\approx$" again denotes uniformly bounded hyperbolic distance in $\mathbf{H}^2$.

Recall that $\nu_\pm(\alpha)$ is the value associated to $\nu_\pm$ under some fixed normalization of $\mathbf{H}^2$ taking $\alpha$ to $\infty$. The difference $\nu_+ - \overline{\nu}_-$ is invariant under real translations, so we may as well assume the normalization is one that takes $\beta$ to $0$ (that is, $\nu_\pm(\alpha) = \tau(S, \nu_\pm, \alpha, \beta)$).

If $\alpha = \alpha_+$, then by construction $\sigma_1$ is in the conformal class of $\nu_+$, and $\tau(S, \sigma_1, \alpha, \beta)$ is exactly $\nu_+(\alpha)$.

If $\alpha = \alpha_n \ne \alpha_+$, then $\sigma_1 = \sigma_{\alpha_n, \alpha_{n+1}}$, which by construction has Teichmüller parameter $\alpha_{n+1} + i$, where $\alpha_{n+1}$ is an integer in this normalization. Now, $\nu_+(\alpha_n)$ by definition of the pivot sequence must be within Euclidean distance 1 of $\alpha_{n+1}$ – either $\nu_+$ is separated from $\alpha_n = \infty$ by a Farey edge with integer endpoints adjacent to $\alpha_{n+1}$, or it is in the horoball of diameter 1 tangent to $\mathbf{R}$ at $\alpha_{n+1}$ if $\alpha_{n+1} = \alpha_+$. It follows that $\alpha_{n+1} + i$ and $\nu_+(\alpha_n) + i$ are within bounded hyperbolic distance.

A similar argument applies to $\nu_-$, and the approximation (10.4) follows. Note that we need not be precise here about multiples of $i$.

Putting these estimates together with (10.2), we have the desired estimate on $\tau(\partial U, \alpha, \mu)$.

*Remark on the internal case.* This computation may seem intricate because it takes care of all special cases at once, but for an internal pivot $\alpha_n$ we simply note that $\mu_1$ and $\mu_0$ differ by exactly $w(n)$ twists around $A$, and that each of the four annuli making up $\partial U$ have modulus exactly $1/2$. Thus we obtain $\tau(\partial U, \alpha, \mu) = w(n) + 2i$.

It remains to relate the conformal structure on $\partial U$ to that on $F_0$.

We claim first that the map $H|_{\partial U_n} : \partial U_n \to \partial \mathbf{T}_{\varepsilon_0}(\alpha_n)$ has degree 1. This follows from the fact that $H$ itself has degree 1, takes the pair $(U_n, \partial U_n)$ to $(\mathbf{T}_{\varepsilon_0}(\alpha_n), \partial \mathbf{T}_{\varepsilon_0}(\alpha_n))$, and takes nothing else into $\mathbf{T}_{\varepsilon_0}(\alpha_n)$.

The map $H$ is also uniformly Lipschitz on the torus $\partial U_n$, and the target torus $F_0$ is itself Euclidean with area no less than $A_0$ by part (3) of Lemma 6.3. When $\alpha_n$ is an internal pivot, $\partial U_n$ is uniformly bilipschitz equivalent to a square torus. By the discussion in Section 2.3, this suffices to imply that $H|_{\partial U_n}$ is homotopic to a uniformly quasiconformal homeomorphism.

When $\alpha_n$ is $\alpha_+$ or $\alpha_-$ (or both), the torus $\partial U_n$ contains high modulus annuli $A_+$ and/or $A_-$. On these, the map $H$ is already uniformly quasiconformal. Since it is already Lipschitz on the remaining annuli, which have bounded modulus, we can apply Lemma 2.2 to show $H$ is homotopic to a quasiconformal map, with a uniform bound on distortion.

Thus, the Teichmüller parameter $\omega(F_0)$ is approximated by $\nu_+(\alpha_n) - \overline{\nu}_-(\alpha_n) + i$ up to bounded distance in $\mathbf{H}^2$, and therefore so is $\omega(\alpha_n)$. This concludes the proof of the pivot theorem.

## 11. Model manifolds and geometric limits

The pivot theorem and the constructions that led to its proof give us what looks like a fairly complete picture of a punctured torus manifold, as a sequence of standard blocks encasing Margulis tubes whose geometry is controlled by the pivot sequence. However, the map from the model to the hyperbolic manifold is so far only Lipschitz, with no bounds in the opposite direction. Such bounds are needed in order to determine the quasi-isometry type of our punctured torus group, and this is the purpose of the following theorem:

THEOREM 11.1 (model manifold). *Let $\rho : \pi_1(S) \to \mathrm{PSL}_2(\mathbf{C})$ be a marked punctured-torus group with quotient manifold $N$. The end invariants $(\nu_-, \nu_+)$ determine a manifold $M$ homeomorphic to $\widehat{C}(N)$, and a homotopy equivalence $f : M \to \widehat{C}(N)$ whose lift $\widetilde{f}$ to the universal covers is a quasi-isometry. The constants for the quasi-isometry are independent of $\rho$.*

For us a map $F : X \to Y$ between metric spaces is a quasi-isometry if

$$\frac{1}{K}d(x,y) - \delta \le d(F(x), F(y)) \le Kd(x,y) + \delta$$

for all $x, y \in X$, where $K$ and $\delta$ are positive constants.

11.1. *Building the model.* Much of the work of building the model has already been done by the block construction in Sections 9 and 10. In particular we have a manifold $M$ and a proper map $H : M \to \widehat{C}(N)$ which is a degree-1 homotopy equivalence, and a metric on the submanifold $M^0$ with respect to which $H|_{M^0}$ is Lipschitz (with uniform constant).

What remains is to define the metric on the rest of $M$ – that is, on the solid tori $U_n$, and adjust $H$ on those parts to extend the Lipschitz condition. Each torus $\partial U_n$ is Euclidean, and we have fixed a marking by $\alpha_n$ and the meridian curve $\mu_n$. It has a Teichmüller parameter $\tau_n \equiv \tau(\partial U_n, \alpha_n, \mu_n)$ which we computed in the previous section. In particular for internal blocks $\tau_n$ is just $w(n) + 2i$.

The length of $\alpha_n$ on $\partial U_n$ is exactly 1, and it is an easy exercise to show that there is a unique Margulis tube whose boundary is isometric to $\partial U_n$ with the same marking. That is, the Euclidean structure on $\partial U_n$ determines a complex translation length $\lambda'$ and a radius $r$. Furthermore the Teichmüller parameter at infinity $\omega' = 2\pi i / \lambda'$ for this tube will be within bounded Teichmüller distance of $\tau_n$, and the distance becomes smaller for large $w(n)$ (the radius $r$ goes to $\infty$ as $w$ does). Note that this is not quite an $\varepsilon_0$ tube because the translation distance at the boundary is not exactly $\varepsilon_0$, but it is an $\varepsilon'$-tube where $\varepsilon'$ is within bounded ratio of $\varepsilon_0$.

We therefore extend the metric on $U_n$ to make it isometric to this Margulis tube. The pivot theorem assures us, for a pivot $\alpha_n$ with length at most $\varepsilon_3$, that this Margulis tube admits a uniformly bilipschitz homeomorphism to the $\varepsilon_0$-Margulis tube of the pivot $\alpha_n$ in $N$ itself. Our map $H$ restricted to $\partial U_n$ is only a Lipschitz map of degree 1, but admits a bounded homotopy to an affine isomorphism of Euclidean tori, and after realizing this homotopy in a collar we can extend to the rest of $U_n$, obtaining a map homotopic to $H$ rel boundary on the solid torus. (The discussion is similar in the case of boundary blocks with parabolics, where $\tau_n$ is infinite in that case.)

Note, by the way, that because of the corners in $\partial U_n$ the resulting metric is not smooth, but this does not affect any of our arguments.

For pivots whose length is bounded *below* by $\varepsilon_3$, the widths $w(n)$ and hence the meridian lengths are bounded, and we can obtain a metric on $U_n$ in which the radius is bounded away from 0 and $\infty$. Therefore $H$ on $U_n$ can be homotoped rel boundary to a uniformly Lipschitz map taking $U_n$ into a neighborhood of the geodesic $\alpha_n^*$.

Let $f$ be the resulting map on $M$ with its finished metric. We have so far proven that $f$ is a degree 1 Lipschitz homotopy equivalence. To show that $f$ lifts to a quasi-isometry, we will need some uniform estimates, which we will obtain by considering geometric limits. Our first goal then is Lemma 11.2, which gives us a uniform description of geometric limits of punctured-torus groups.

11.2. *Geometric limits.* We assume familiarity with the geometric topology on metric spaces with basepoints $(N, x)$ (see [9], [24], [51], [86] and also [69] for applications in the spirit of this paper). For our purposes we use these basic definitions: a sequence $(N_i, x_i)$ converges geometrically to $(N, x)$ if and only if there is a sequence of numbers $R_i \to \infty$ and $K_i \to 1$, and maps $h_i : (B(x, R_i), x) \to (N_i, x_i)$ which are $K_i$-bilipschitz homeomorphisms to their images where $B(x, R)$ denotes an $R$-neighborhood of $x$. Similarly a sequence of mappings $f_i : (M_i, x_i) \to (N_i, y_i)$ converges geometrically to $f : (M, x) \to (N, y)$ if we have convergence of $(M_i, x_i)$ to $(M, x)$ and $(N_i, y_i)$ to $(N, y)$ as above (with bilipschitz maps $h_i$ and $k_i$, respectively), and if for any fixed $R$ the maps $h_i^{-1} \circ f_i \circ k_i : B(x, R) \to N$ (defined for large enough $i$) converge in $C^0$ to $f$.

An important phenomenon in hyperbolic geometry is that the set of pointed hyperbolic $n$-manifolds with a fixed positive lower bound on injectivity radius at the basepoint is compact in the geometric topology.

In this subsection, we try to understand geometric limits of punctured-torus manifolds, by considering geometric limits of their models and the maps $f$.

*Convergence of combinatorial data.* The model manifolds are determined by combinatorial data which it is convenient to describe as follows. Recalling that the pivot sequence is $P$ is indexed by an interval $J \subset \mathbf{Z}$, we have a sequence of integer *widths* $W = \{w(n)\}_{n \in J}$, determined by the pair $(\nu_-, \nu_+)$ as in Section 4. Furthermore, if $P$ has a first (or last) element we define $v_-$ (or $v_+$) as a complex number with $|\operatorname{Re} v_\pm| \leq 1$ which equals $\nu_\pm(\alpha_\pm) \bmod 1$. More specifically, with notation as in Section 4, if $\alpha_- \neq \alpha_+$ we let $v_+$ be such that $w(p+1) + v_+ = \nu_+(\alpha_+) - \alpha_p(\alpha_+)$, and $v_-$ be such that $w(0) + v_- =$ $\alpha_1(\alpha_-) - \overline{\nu}_-(\alpha_-)$. If $\alpha_- = \alpha_+$ then we want $w(0) + v_- + v_+ = \nu_+(\alpha_+) - \overline{\nu}_-(\alpha_+)$.
In case $e_+$ or $e_-$ has an accidental parabolic we write $v_+ = \infty i$ or $v_- = \infty i$.

Note in particular that $\operatorname{Im} v_+$ and $\operatorname{Im} v_-$ determine the moduli of the annuli $A_+$ and $A_-$ in the corresponding boundary blocks, via (10.3). The triple $\Sigma = (W, v_+, v_-)$ is called the *combinatorial data* (where $v_\pm$ are ignored if undefined).

Note that the pivot sequence itself is forgotten by the combinatorial data, but can be recovered from the width sequence up to the action of $\operatorname{SL}_2(\mathbf{Z})$. This corresponds to keeping track of the image Kleinian group but forgetting the marking, i.e. the representation $\rho$.

To understand limiting configurations, let us generalize slightly the definition of the combinatorial data $\Sigma$ to allow the index set $J$ to be any interval in $\mathbf{Z}$, and to allow the value $w(n) = \infty$. With this definition, we consider a topology on the set of $\Sigma$'s as follows: We say that a sequence $\Sigma_i = (W_i, v_+^i, v_-^i)$ converges to $\Sigma$ if the index sets $J_i$ converge to $J$ as subsets of $\mathbf{Z}$ (that is, for any $q > 0$, eventually $J_i \cap [-q, q] = J \cap [-q, q]$), and for each $n \in J$ either $w_i(n) \to w(n)$ in $\mathbf{Z}$, or $w(n) = \infty$ and $|w_i(n)| \to \infty$. Finally we require if $v_+$ is defined in $\Sigma$ that $v_+^i$ are eventually defined and $v_+^i \to v_+$ in the natural compactification of $[-1, 1] + i\mathbf{R}$ by $\infty i$, and similarly for $v_-$.

*Convergence of model manifolds.* A generalized combinatorial data set $\Sigma$ determines a generalized model manifold $M_\Sigma$. For each $n \in J$, if $w(n)$ is finite we construct a block $\mathcal{B}_n$ just as before, except that the pivots do not have fixed names: for example for an internal block $\mathcal{B}_n$ we build $\mathcal{B}_{\alpha, \beta, \gamma}$ where $\beta$ and $\gamma$ are neighbors of $\alpha$ and $\gamma = D_\alpha^{w(n)} \beta$. When gluing successive blocks we must identify $\partial_1 \mathcal{B}_n$ to $\partial_0 \mathcal{B}_{n+1}$ by the unique orientation-preserving isometry taking the $\alpha$ curve on $\partial_1 \mathcal{B}_n$ to the $\beta$ curve on $\partial_0 \mathcal{B}_{n+1}$, and the $\gamma$ curve on $\partial_1 \mathcal{B}_n$ to the $\alpha$ curve on $\partial_0 \mathcal{B}_{n+1}$.

When $w(n) = \infty$ and the block is to be internal, we construct first the truncated block $\mathcal{B}_n^0$ isometric to $\mathcal{B} \setminus U$. Note that without $U$ there is no natural identification between top and bottom boundaries, and in fact all truncated internal blocks are isometric. There is a unique Margulis tube for a rank-2 parabolic group whose boundary torus is isometric to $U$, and we put the metric of this tube on $U$ minus its core circle, and glue that into the block.

For a boundary block the same thing works, with the following provisos. We use $\operatorname{Im} v_+$ and/or $\operatorname{Im} v_-$ to determine the moduli of $A_+$ and/or $A_-$, via (10.3), and the real parts together with the corresponding widths $w(0)$ or $w(p+1)$ serve to determine the exact geometry of the torus $\partial U$ and the choice of meridian. The case $v_+ = \infty i$ yields a *pair* of annuli each isometric to $S^1 \times [0, \infty)$, as before, and a rank-1 parabolic Margulis tube (even if the corresponding width is infinite), and similarly $v_- = \infty i$. The one exception is the possibility that $\alpha_+ = \alpha_-$, so that there is just one double boundary block $\mathcal{B}_0$, and *both* $v_+ = v_- = \infty i$. In this case we replace $U$ with *two* copies of a rank-1 parabolic Margulis tube, glued to the annulus $A_L$ and $A_R$ respectively. Thus the manifold is homeomorphic to a thrice-punctured sphere cross $[0,1]$. This case is not central to our discussion (it corresponds to a Fuchsian geometric limit), but we include it for completeness.

We can now describe geometric convergence for the model manifolds more clearly, and relate it to convergence for the combinatorial data $\Sigma$.

Let us adopt the following notation. Let $M = \cup_{n \in J} \mathcal{B}_n$ be a model manifold. Let $\mathcal{B}_n^0$ denote $\mathcal{B}_n$ minus the open solid torus $U_n$. Let $\check{\mathcal{B}}_n$ denote $\mathcal{B}_n$ minus the cusp region $Q_{\mathcal{B}_n}$ and let $\check{\mathcal{B}}_n^0$ denote $\mathcal{B}_n^0 \cap \check{\mathcal{B}}_n$. Then let $M^0 = \bigcup_{n \in J} \mathcal{B}_n^0$, and similarly for $\check{M}$, $\check{M}^0$. For the integer $q > 0$ let $M[q]$ be the finite union of blocks $\mathcal{B}_n$ with $|n| \le q$. Let $M^0[q] = M^0 \cap M[q]$ and define $\check{M}[q]$, $\check{M}^0[q]$ similarly. Note that $\check{M}^0[q]$ is compact.

Now consider a sequence $(M_i, x_i)$, with combinatorial data $\Sigma_i$. Let us assume, re-indexing by shifting $J_i$ if necessary, that $0 \in J_i$ for each $i$, and that $x_i \in \mathcal{B}_{i,0}$ (here we denote the $n^{\text{th}}$ block in $M_i$ as $\mathcal{B}_{i,n}$). Suppose that the $\Sigma_i$ converge to $\Sigma$, and let $M_\infty = M_\Sigma$. We claim that, if the $x_i$ are chosen properly, a subsequence of $(M_i, x_i)$ converges geometrically to $(M_\infty, x)$ for some $x \in \mathcal{B}_{\infty,0}$.

In fact, consider $M_\infty^0[q]$ – it is composed of finitely many truncated blocks $\mathcal{B}_n^0$. By our construction, any two truncated internal blocks are isometric, and any two truncated boundary blocks of the same type are bilipschitz equivalent with some uniform bound, except possibly for the long annuli $A_\pm$ in their boundaries. Thus, for $i$ large enough there is a natural embedding $k_i : M_\infty^0[q] \to M_i[q]$ which is an isometry on each truncated internal block (and bilipschitz for boundary blocks minus the long annuli), and preserves the indexing $n$. Thus if we choose $x_i \in \check{\mathcal{B}}_{i,0}^0$ but not in the long boundary annuli, then the $k_i^{-1}(x_i)$ have a convergent subsequence with limit $x \in \check{\mathcal{B}}_0^0$. Furthermore convergence of the $v_\pm$ insures that the long boundary annuli converge geometrically, and their boundary identifications can vary in a compact set. Thus, letting $q \to \infty$, we find that for some subsequence which we continue to index as if it were the whole sequence, the submanifolds $(M_i^0, x_i)$ converge geometrically to $(M_\infty^0, x)$.

Note, the subtlety here is that the maps $k_i$ do not necessarily extend across the solid tori $U_n$ – the meridian curve $\mu_n$ of the solid torus may be different from the meridian of the corresponding torus in $M_i$. However, this information is encoded in the widths $w_i(n)$. In particular, for internal blocks if $w(n) \ne \infty$ then for large enough $i$, $w_i(n) = w(n)$, and then the map $k_i$ extends to an isometry between the solid tori $U_n$ and $U_{i,n}$. If $w(n) = \infty$ then $w_i(n) \to \pm\infty$, the radii of the corresponding Margulis tubes go to infinity, and it is easy to check that the geometric limit (with basepoint on the boundary torus) is exactly the rank-2 parabolic Margulis tube which we have placed in $M$. (This is essentially Jørgensen's original example of a geometric limit of cyclic groups which is not itself cyclic; see [47], [51]). The case of boundary blocks is similar, except that we use convergence of the $v_\pm$ to give convergence of the long annuli. Furthermore the case where $\operatorname{Im} v_+$ or $\operatorname{Im} v_-$ goes to $\infty$ yields a rank-1 parabolic rather than a rank-2, regardless of the behavior of $w_i(n)$.

We conclude that the $(M_i, x_i)$ converge geometrically to $(M_\infty, x)$. Furthermore, let us extend the definition of $k_i$ to maps $k_i : M_\infty^0 \to M_i$ as follows: for each $i$ choose $q = q_i$ to be the largest $q$ so that $J_i \cap [-q, q] = J_\infty \cap [-q, q]$ (so $q_i \to \infty$ as $i \to \infty$). We define $k_i$ on $M_\infty^0[q_i]$ as before, and extend to the rest of $M_\infty^0$ by first collapsing the blocks outside $M_\infty^0[q]$ to $\partial M_\infty^0[q]$. We then also have geometric convergence of $k_i \to k_\infty$, where $k_\infty : M_\infty^0 \to M_\infty$ is the inclusion map.

Note that, except in one case, a geometric limit of model manifolds with basepoints in the thick part must always be homeomorphic to the product $S \times I$ minus some sequence of curves $\gamma_k \times \{k\}$, where $I$ is an interval in $\mathbf{R}$, each $\gamma_k$ is a nonperipheral simple curve on $S$, and $\gamma_k$ and $\gamma_{k+1}$ are homotopically distinct. The exception is the case where for each $i$ there is just one block, so that $\alpha_-^i = \alpha_+^i$, and both $v_-^i$ and $v_+^i$ go to infinity. Then the limit is the exceptional thrice-punctured sphere case described above, and we leave the details to the reader.

*Convergence of model maps.* The following lemma allows us to describe a geometric limit of any sequence of punctured torus groups in terms of a limit of models.

LEMMA 11.2. *Let $N_i$ be a sequence of punctured-torus manifolds with associated models $M_i$ and maps $f_i : (M_i, x_i) \to (\widehat{C}(N_i), f_i(x_i))$. Assuming that the $x_i$ are in the interior of the thick part $\check{M}_i^0$, these maps have a geometrically convergent subsequence, and the limit $f_\infty : (M_\infty, x_\infty) \to (\widehat{C}_\infty, y_\infty)$ is a proper homotopy-equivalence of degree 1, where $\widehat{C}_\infty$ is a subset of $N_\infty$, homotopy equivalent to it and containing its convex core.*

*Proof.* Re-index as above, if necessary, so that the index sets $J_i$ of the combinatorial data $\Sigma_i$ always contain 0, and that $x_i \in \check{\mathcal{B}}_{i,0}^0$. Then after possibly restricting to a subsequence we may assume that $\{\Sigma_i\}$ converges to some $\Sigma$. By the previous discussion, after restricting to another subsequence, we may assume that $(M_i, x_i)$ converges geometrically to some $(M_\infty, x_\infty)$ (keeping $x_i$ in the interior means it can't lie on the long annuli $A_\pm$ of boundary blocks). By Lemma 9.3, the injectivity radius of $N_i$ at $f_i(x_i)$ is uniformly bounded below by $\varepsilon_3$, and therefore we may further assume that $(N_i, f_i(x_i))$ converges geometrically to $(N_\infty, y_\infty)$, and that the subsets $\widehat{C}(N_i)$ converge to some subset $\widehat{C}_\infty$ of $N_\infty$. The maps $f_i$ are uniformly Lipschitz, so after extracting a further subsequence we may assume that they converge to $f_\infty$. We will abuse notation by continuing to index each subsequence as if it were the whole sequence.

Since our manifolds are $K(\pi,1)$'s, showing that $f_\infty$ is a homotopy equivalence reduces to showing that $(f_\infty)_* : \pi_1(M_\infty) \to \pi_1(N_\infty)$ is an isomorphism. The scheme for this proof is due to Thurston, in [86], [85].

On finite blocks, the map $k_i$ defined above eventually factors through $k_\infty$ in this sense: Fix $q > 0$, and $I$ large enough that, for $i > I$ and $n \in J_\infty \cap [-q, q]$, $w_i(n) = w_\infty(n)$ if the latter is not $\infty$. When $w_\infty(n) = \infty$, $U_n$ is missing its core circle and hence retracts to its boundary. It follows that for $i > I$ we can extend $k_i|_{M_\infty^0[q]}$ to an embedding $k_i^q : M_\infty[q] \to M_i[q]$ satisfying $k_i = k_i^q \circ k_\infty$ on $M_i[q]$.

The manifold $M_\infty[q]$ has the form of a product $S \times [0,1]$ minus a number of level curves corresponding to those $n$ where $w_\infty(n) = \infty$, and the map $k_i^q$ has the effect of filling in those missing curves by a Dehn surgery whose meridian curve is determined by $w_i(n)$. Since $w_i(n) \to \infty$ as $i \to \infty$, there is a precise sense in which the map $k_i^q$ is "eventually injective" on $\pi_1$:

LEMMA 11.3. *Let $S$ be a surface of finite type and negative Euler characteristic. Let $M$ be $S \times [s,t]$ minus a number of regular neighborhoods of level simple curves $\gamma_i \times \{n_i\}$, where $n < n_{i+1}$, each $\gamma_i$ is nonperipheral, and successive curves $\gamma_i, \gamma_{i+1}$ have essential intersection in $S$. Let $M_i$ be a sequence of Dehn fillings on $M$ such that the lengths of meridian curves on each torus boundary go to infinity.*

1. *If the images of $\alpha$ in all the $\pi_1(M_i)$ are trivial, then $\alpha$ is already trivial in $\pi_1(M)$.*

2. *If $\alpha \in \pi_1(M)$ is not in a conjugate of one of the boundary torus subgroups, and its images in all the $\pi_1(M_i)$ are $m^{\text{th}}$ powers, then $\alpha$ is itself an $m^{\text{th}}$ power.*

The proof of this lemma is postponed to Section 11.4.

To apply the lemma, note that the condition on successive level curves follows automatically from the way we glue successive blocks. The condition on meridian lengths follows from the widths going to infinity.

We can now show that $(f_\infty)_*$ is injective. Let $\alpha$ be a nontrivial loop in $M_\infty$. Then we can assume it lies in $M_\infty^0[q]$ for some $q > 0$, and applying part (1) of Lemma 11.3 we conclude that, for large enough $i$, $k_i^q(\alpha)$ is nontrivial in $\pi_1(M_i[q])$. In fact it is nontrivial in $\pi_1(M_i)$ since $\pi_1(M_i[q])$ injects in $\pi_1(M_i)$ by Seifert-van Kampen.

Now if $f_\infty(\alpha)$ is trivial, it bounds a disk in $N_\infty$. For high enough $i$ the disk pulls back to $N_i$, and we conclude that $f_i(k_i^q(\alpha))$ is also homotopically trivial. This contradicts the fact that, by construction, $f_i$ are all $\pi_1$-injective, so we conclude that $(f_\infty)_*$ is injective.

We next prove that $(f_\infty)_*$ is surjective. Let $\Gamma = (f_\infty)_*(\pi_1(M_\infty))$ and $\hat{\Gamma} = \pi_1(N_\infty)$, considered as subgroups of $\mathrm{PSL}_2(\mathbf{C})$. Up to suitably normalizing via conjugation in $\mathrm{PSL}_2(\mathbf{C})$, geometric convergence of the manifolds is equivalent to geometric convergence of the groups in the Gromov-Hausdorff topology (see [24]). Thus $\hat{\Gamma}$ is the geometric limit of $(f_i)_*(\pi_1(M_i))$. We need to show that $\Gamma = \hat{\Gamma}$.

We also have a sequence of (surjective but noninjective) representations $\varphi_i = (f_i \circ k_i)_*$ from $\pi_1(M_\infty^0)$ to $\mathrm{PSL}_2(\mathbf{C})$, and $\varphi_\infty$ similarly, such that $\varphi_i \to \varphi_\infty$ and $\varphi_\infty(\pi_1(M_\infty^0)) = \Gamma$.

We first prove that if $\Gamma$ is a finite-index subgroup of $\hat{\Gamma}$ then $\Gamma = \hat{\Gamma}$. Finite index implies that for any $g \in \hat{\Gamma}$, some finite power $g^m$ is in $\Gamma$. Thus $g^m = \varphi_\infty(h)$ for some $h \in \pi_1(M_\infty^0)$.

Choose $q$ large enough so that $h$ can be represented by a loop in $M_\infty^0[q]$, and $h' = (k_\infty)_*(h)$ can thus be assumed to lie in $\pi_1(M_\infty[q])$. Suppose that $h'$ is contained in a (conjugate of) a rank-2 cusp subgroup of $M_\infty$. Since $f_\infty$ maps rank-2 cusps to rank-2 cusps with degree 1, the corresponding cusp subgroup of $\hat{\Gamma}$ is contained in the image of $\varphi_\infty$. Note that $g$ is contained in this cusp group since it commutes with $\varphi_\infty(h)$, so it must already be contained in $\Gamma$.

Now assume that $h'$ is not in a cusp group. By definition of geometric limit of groups, $g = \lim \varphi_i(g_i)$, for some sequence $g_i \in \pi_1(M_\infty^0)$. Thus $\lim_{i\to\infty} \varphi_i(g_i^m) = \varphi_\infty(h)$. Since also $\lim_{i\to\infty} \varphi_i(h) = \varphi_\infty(h)$ and the limit group is discrete, we conclude (recapitulating an argument of Jørgensen-Marden in [51]) that for large enough $i$, $\varphi_i(g_i)^m = \varphi_i(h)$. Note that since $\varphi_i$ may not be injective, we cannot conclude that $g_i^m = h$. However we do know that

$$(11.1) \qquad (k_i)_*(g_i)^m = (k_i)_*(h),$$

since $(f_i)_*$ is injective. Hence $(k_i^q)_*(h') = (k_i)_*(h)$ is eventually an $m^{\mathrm{th}}$ power, and since $h'$ is not in a cusp subgroup we can apply part (2) of Lemma 11.3 to conclude that $h'$ is itself an $m$-th power, i.e. $h' = (g')^m$ for $g' \in \pi_1(M_\infty)$. Since $m^{\mathrm{th}}$ roots are unique in a torsion-free Kleinian group (see also [45, Chap. VI]), we conclude that $(f_\infty)_*(g') = g$, and thus $g \in \Gamma$ after all. Thus $\Gamma = \hat{\Gamma}$ in this case.

Now consider the possibility that $\Gamma$ has infinite index in $\hat{\Gamma}$. In this case we use a variation of Thurston's covering theorem (similar to an argument made in [85]).

Let $N_a = \mathbf{H}^3/\Gamma$, so that there is a covering $\pi : N_a \to N_\infty$ and a lift $f_a : M_\infty \to N_a$ of $f_\infty$ with $\pi \circ f_a = f_\infty$, and $f_a$ is a homotopy equivalence. Since $f_\infty(M_\infty) = \widehat{C}_\infty$, which is the geometric limit of sets containing the convex cores of $C(N_i)$, it contains the convex core of $N_\infty$ (see Kerckhoff-Thurston [55]).

Since $C(N_\infty)$ is a deformation retract of $N_\infty$, every loop in $N_\infty$ can be deformed into $f_\infty(M_\infty)$. It follows that the preimage $\pi^{-1}(f_\infty(M_\infty))$ is equal to $f_a(M_\infty)$, and that if $\pi$ is infinite-to-one on $N_a$ it must already be infinite-to-one on $f_a(M_\infty)$.

However, restricted to each rank-2 cusp of $M_\infty$, the map must clearly be finite-to-one. It follows that there exists a sequence of blocks $\mathcal{B}_{n_j}$ with $|n_j| \to \infty$, whose images in $N_\infty$ meet a compact set. The map $f_\infty$ restricted to, say, $\partial_0 \mathcal{B}_{n_j}$, is a limit of a sequence of our "halfway surfaces", and in particular is a uniformly Lipschitz map $h_{n_j}^\infty : S \to N_\infty$, where the metric on $S$ is the standard metric induced from any of the block boundaries, up to homeomorphism. Since these surfaces meet a compact subset of $N_\infty$ minus its main cusp (and are standard in the cusp) there must be two of them $h_1, h_2$ which admit a homotopy in $N_\infty$ with tracks shorter than the injectivity radius in this compact set. We may lift this homotopy to the covering $N_a$ and conclude that a finite-volume 3-chain in $N_a$ maps to a closed 3-chain in $N_\infty$ (closed relative to the cusps). It follows that $N_\infty$ has finite volume, but this is impossible since it is the geometric limit of a sequence of manifolds with infinite volume. (See Canary [23] for a discussion of the covering theorem and other generalizations.)

We conclude that $\Gamma$ cannot have infinite index in $\hat{\Gamma}$ either, and therefore $(f_\infty)_*$ is surjective.

The same argument also shows that $f_\infty$ is proper, for we know it takes boundary to boundary, and we have shown that the images of blocks $\mathcal{B}_n$ cannot accumulate in a compact subset of $N_\infty$.

We now show that the degree is 1, which is a more subtle issue. Let us first make the following observations about ordering of ends and the orientation convention of Section 3. Recall that an ordering for the ends $(e_-, e_+)$ of $S_0 \times \mathbf{R}$ is determined by a choice of orientation for $S_0$. If $\psi : (S_0, \partial S_0) \to (S_0 \times \mathbf{R}, \partial S_0 \times \mathbf{R})$ is an immersion inducing the identity on $\pi_1(S_0)$, then we can also detect the ordering of the ends in this way: there are exactly two components of $S_0 \times \mathbf{R} \setminus \psi(S_0)$ which are *unbounded* (i.e. whose closures are not compact), and we can order them $(U_-, U_+)$ according to the convention that an oriented path from $U_+$ to $U_-$ has intersection number $+1$ with $\psi(S_0)$. One can easily check that $U_-$ is a neighborhood of $e_-$, and $U_+$ of $e_+$.

Now if $h : (S_0 \times \mathbf{R}, \partial S_0 \times \mathbf{R}) \to (S_0 \times \mathbf{R}, \partial S_0 \times \mathbf{R})$ is proper and induces the identity on $\pi_1$ then we can also check that if it induces a bijection of the ends which is order-preserving then $\deg h = 1$, if it induces an order-reversing bijection then $\deg h = -1$, and if it maps both ends to one then the degree is 0.

For any model map $f : M \to N$ recall that $\check{M}$ denotes the complement in $M$ of the main cusp neighborhood $Q_M$, and similarly $\check{N}$ is the complement in $N$ of $Q_N$, and let $\widehat{C}(\check{N}) = \widehat{C}(N) \cap \check{N}$. Let $\partial' \check{M}$ denote the cylindrical frontier of $Q_M$ in $M$, and let $\partial' \widehat{C}(\check{N})$ similarly denote the frontier of $Q_N \cap \widehat{C}(N)$ in $\widehat{C}(N)$. By our construction, $f$ restricts to a map of pairs $f : (\check{M}, \partial'\check{M}) \to (\widehat{C}(\check{N}), \partial'\widehat{C}(\check{N}))$, and in the case where the width sequence is bi-infinite and the entries $w(n)$ are all finite, each pair may be identified with $(S_0 \times \mathbf{R}, \partial S_0 \times \mathbf{R})$. Let $S_0(n)$ denote a level surface of the form $\partial_0 \check{\mathcal{B}}_n$ in $M$.

Returning to our sequence $f_i \to f_\infty$, consider first the case that the limiting width sequence $w_\infty$ is bi-infinite, and that $w_\infty(n)$ does not take on the value $\infty$. If $\deg f_\infty$ is not 1, then at least one of the ends of $\check{M}_\infty$ must map to the opposite end of $\check{N}_\infty$. This situation is covered by a special case of a theorem proved by Thurston for arbitrary genus [86], in which he also shows that the ends cannot flip in the limit. We will give a slightly different argument here, adapted to our techniques but following the same general logic.

Let $(U_-, U_+)$ denote the unbounded components of $\check{N}_\infty \setminus f_\infty(S_0(0))$, ordered as above by the orientation convention. Then order reversal implies, without loss of generality, that as $j \to +\infty$, $f_\infty(S_0(j))$ is eventually in $U_-$. Fix such a $j$. Let $D$ be the maximal diameter of any bounded component of $\check{N}_\infty \setminus f_\infty(S_0(0))$. Let $X_\infty$ be a neighborhood of $f_\infty(S_0(0))$ in $\check{N}_\infty$ of radius at least $D + 1$, and large enough to contain $f_\infty(S_0(j))$. By the definition of geometric convergence, given $K$ as close as we like to 1, for large enough $i$ there is a region $X_i$ in $\check{N}_i$ and a $K$-bilipschitz homeomorphism $h_i : X_\infty \to X_i$. For large $i$, $w_i(n) = w_\infty(n)$ for $|n| \leq j$, so the embedding $k_i : M_\infty^0[j] \to M_i^0[j]$ extends to an isometry $k_i^j : M_\infty[j] \to M_i[j]$, as in our previous discussion. Furthermore, $h_i^{-1} \circ f_i \circ k_i^j$ is $\varepsilon$-close to $f_\infty$ in the $C^0$ topology, for $\varepsilon$ small. Thus $X_i$ must contain a $(D+1-\varepsilon)/K$-neighborhood of $f_i(S_0(0))$ (see e.g. [21]) so that (for $K$ close to 1, $\varepsilon$ close to 0) the unbounded components $U_\pm^i$ of $\check{N}_i \setminus f_i(S_0(0))$ must each intersect $X_i$ in sets of diameter greater than $D$. It follows that $h_i(U_\pm \cap X_\infty) \subset U_\pm^i$. Because $h_i$ is orientation-preserving, the ordering of $+$ and $-$, as detected by intersection numbers, is preserved.

Thus, $f_i(S_0(j))$ must be contained in $U_-^i$. On the other hand since $f_i$ itself does not reverse end order, for some $k > j$ we have $f_i(S_0(k)) \subset U_+^i$. It follows that some block $\mathcal{B}_{i,m_j}$ with $j < m_j < k$ has image under $f_i$ that meets $f_i(S_0(0))$. Repeating this argument with arbitrarily high $j$ we obtain a sequence of pivots $\alpha_{m_j}^{i_j}$ of the pivot sequence for $N_{i_j}$, with $m_j > j$, such that the geodesics $(\alpha_{m_j}^{i_j})^*$, or their $\varepsilon_0$-Margulis tubes, remain a bounded distance away from $f_i(S_0(0))$. Each $(\alpha_{m_j}^{i_j})^*$ has uniformly bounded length by the pivot theorem, so after restriction to a subsequence their images under $h_i^{-1}$ must converge to a geodesic $\beta$ in $N_\infty$, to which they are eventually homotopic. Since $f_\infty$ is a homotopy equivalence, $\beta$ determines a conjugacy class in $\pi_1(M_\infty) = \pi_1(S)$, which is eventually equal to the conjugacy class of $\alpha_{m_j}^{i_j}$; but this contradicts the fact that, on the other hand, $\alpha_{m_j}^{i_j}$ must converge to $\nu_+^\infty$ (compare Lemma 12.1), which is irrational in this case. We conclude that $f_\infty$ cannot reverse the order of the ends. (In Thurston's argument the last contradiction is obtained in terms of *realizability* of the ending lamination.)

Proving deg $f_\infty = 1$ in the remaining cases is easier. Suppose first that the limiting width sequence is not bi-infinite. Then $\widehat{C}(N_\infty)$, and $\widehat{C}(N_i)$ for high $i$, have nonempty boundary, and $f_i : \partial M_i \to \partial \widehat{C}(N_i)$ is constructed to be a uniformly bilipschitz orientation-preserving homeomorphism, which converges to $f_\infty : \partial M_\infty \to \partial \widehat{C}(N_\infty)$. It follows that the degree is 1.

Finally, if $w_\infty(n) = \infty$ for some $n$, then $N_\infty$ has a rank-2 parabolic cusp which is the image of a corresponding cusp $U_n$ in $M_\infty$. By our construction, $f_\infty$ (and any $f_i$) has degree 1 over this cusp, and hence degree 1 globally since it is proper.

This concludes the proof of Lemma 11.2. □

11.3. *Lifting to a quasi-isometry.* With Lemma 11.2 in hand, we may prove the following corollary, which essentially says that the map $f$ does not contract too much, in a homotopic sense:

LEMMA 11.4. *Let $f : M \to N$ be a model for a punctured torus group. For any $B > 0$ there exists $A > 0$ such that, if $\beta$ is a loop in $N$ through $f(x)$ of length no more than $B$, then there exists a loop $\alpha$ through $x$ in $M$ of length no more than $A$, such that $f(\alpha)$ is homotopic to $\beta$, fixing the basepoint $f(x)$.*

*Proof.* For any fixed $x \in M$ there is such an $A(x, B)$, just because $f$ is a homotopy equivalence and there is a finite number of homotopy classes of loops with a given length bound passing through a given point in a hyperbolic manifold. If there were no universal $A(B)$ we could find a sequence $x_i$ in $M$ so that the best $A(x_i, B)$ go to infinity. Clearly $x_i$ would not go out to a cusp because there $f$ is uniformly bilipschitz. Therefore we may assume that $x_i \in \check{\mathcal{B}}^0_{n_i}$ where $|n_i| \to \infty$. Thus if we consider the sequence of manifolds with shifted basepoint $(M, x_i)$ we may apply Lemma 11.2 to obtain a geometric limit $f_\infty : (M_\infty, x_\infty) \to (N_\infty, y_\infty)$ which is still a homotopy equivalence. Therefore there is some $A$ which works in the limit. Now pulling back large neighborhoods of $x_\infty$ and $y_\infty$ to the approximants, we see that $A + \delta$ for some $\delta > 0$ works there too. (Compare a similar argument in Lemma 4.5 of [69]). □

Finally we claim that a model map $f : M \to N$ is a quasi-isometry. Restricted to each Margulis tube $U_n$, $f$ is by construction a bilipschitz homeomorphism to its image. Since each internal truncated block minus its Margulis tube, $\check{\mathcal{B}}^0_n$, has bounded diameter, we must prove that the separation between the images of any two blocks, $f(\check{\mathcal{B}}_n)$ and $f(\check{\mathcal{B}}_{n'})$, as measured in the truncated manifold $\check{N}$, is roughly proportional to $|n - n'|$. This will prove that the restriction of $f$ to $\check{M} \to \check{N}$ is a quasi-isometry, and by construction the extension to the main cusp neighborhoods will be a quasi-isometry too.

Again let $S_0(n)$ denote $\partial_0 \check{\mathcal{B}}_n$. Let us prove the following claim: there is some $m_0 > 0$ such that, whenever $n \in J$ and $n + m_0 \in J$, there is some $0 < m \le m_0$ such that $f(S_0(n+m))$ separates $f(S_0(n))$ from $e_+$ in $\check{N}$. (Equivalently we can make the claim for $m_0 < 0$ and $e_-$). Suppose this is not the case. Then there are sequences $m_i \to +\infty$ and $n_i \to \pm\infty$ such that for every $0 < m \le m_i$, $f_i(S_0(n_i+m))$ fails to separate $f_i(S_0(n_i))$ from the $e_+$ end of $\check{N}$. In other words, if $(U_-^i, U_+^i)$ are the unbounded components of $\check{N} \setminus f(S_0(n_i))$, ordered by our orientation convention, then we have for all $0 < m \le m_i$ that $f(S_0(n_i + m))$ is not contained in $U_+^i$.

Choosing basepoints $x_i$ in $\mathcal{B}_{n_i}$ and possibly restricting to a subsequence, we obtain by Lemma 11.2 a geometric limit $f_\infty : M_\infty \to N_\infty$, which is a proper degree 1 homotopy equivalence. Note that in the limit the combinatorial data yield a bi-infinite sequence since $n_i \to \pm\infty$ – thus we need not consider boundary blocks.

We claim that $\check{N}_\infty^0 \setminus f_\infty(S_0(0))$ has two unbounded components $U_\pm^\infty$. (This is in contrast to $\check{N}_\infty$ which may have infinitely many ends corresponding to rank-2 cusps.) Recalling that $f_\infty$ restricts to a map of pairs $f_\infty : (\check{M}_\infty^0, \partial \check{M}_\infty^0) \to (\check{N}_\infty^0, \partial \check{N}_\infty^0)$, we can perform a series of integer Dehn fillings on the torus boundaries of $\check{M}_\infty^0$ to obtain a manifold homeomorphic to $S_0 \times \mathbf{R}$. After performing fillings with the same integers on the torus boundaries of $\check{N}_\infty^0$, we can extend $f_\infty$ to a map between the filled manifolds that is still a proper, degree 1 homotopy equivalence. It follows that the filled $\check{N}_\infty^0$ is also $S_0 \times \mathbf{R}$ (see [14], [86]) and hence has two ends, and that the extended $f_\infty$ is, again, order-preserving on the ends. Thus $\check{N}_\infty^0$ also has two ends and we define $U_\pm^\infty$ to be the components of $\check{N}_\infty^0 \setminus f_\infty(S_0(0))$ which are neighborhoods of these ends, ordered again by our intersection-number convention.

We now see that there exists $m_\infty > 0$ such that $f_\infty(S_0(m_\infty)) \subset U_+^\infty$. We may now apply to $\check{N}_\infty^0$ the same argument as in the proof of Lemma 11.2 to see that there is a large neighborhood $X_\infty$ of $f_\infty(S_0(0))$ in $\check{N}_\infty^0$ and a sequence of near isometries $h_i : X_\infty \to X_i$, where $X_i$ are neighborhoods of $f(S_0(n_i))$ in $N$, and $h_i$ takes $U_+^\infty \cap X_\infty$ into $U_+^i$. (We have to make $X_\infty$ include all the bounded components of $\check{N}_\infty^0 \setminus f_\infty(S_0(0))$, as before.) Moreover if we choose $X_\infty$ large enough we may push $f_\infty(S_0(m_\infty))$ via $h_i$ into $X_i$ and conclude that $f(S_0(n_i + m_\infty)) \subset U_+^i$. But for $i$ large enough $m_\infty < m_i$, a contradiction.

This proves the claim that $m_0$ exists as above. It follows that for our original model $f : M \to N$ we can find a sequence $s = (\cdots < n_1 < n_2 < \cdots)$ with $n_{i+1} \le n_i + m_0$ and $\inf s = \inf J$, $\sup s = \sup J$, such that $f(S_0(n_j))$ separates $f(S_0(n_i))$ from $f(S_0(n_k))$ in $\check{N}_\infty$ for all $i < j < k$. All these immersed surfaces are disjoint and any two with adjacent indices are a definite distance apart (again by a geometric limit argument). This forces the separation between any $f(S_0(n))$ and $f(S_0(n'))$ to be at least $a_1|n - n'| - a_2$ for some constants $a_i > 0$.

On the other hand we know $f$ is Lipschitz, so the separation is bounded above by another linear function of $|n - n'|$. This concludes the proof that $f$ is a quasi-isometry.

We now have a map $f$ which satisfies the following conditions: It is a surjective homotopy equivalence, Lipschitz, a quasi-isometry, and Lemma 11.4 holds. It follows by lemma 3.1 of [70] that the lift $\widetilde{f}$ to the universal covers is a quasi-isometry. This completes the proof of Theorem 11.1 (model manifold), modulo a leftover lemma:

*Proof of Lemma* 11.3. A version of this lemma is asserted by Thurston in [85] without proof. Although the statement is purely topological we present here a proof which uses hyperbolic geometry in what seems like an essential way. (We have considered applying the boundary-slope techniques of Gordon-Litherland and Gordon-Luecke et al. [35], [36], [28], but without success.)

In fact, consider the following more general situation. Let $M$ be the interior of a compact 3-manifold $\overline{M}$, such that $M$ admits a complete hyperbolic metric $\sigma$. Suppose that $F$ is a collection of torus components of $\partial \overline{M}$, and let $M_i$ be a sequence of Dehn fillings of $\overline{M}$, obtained by attaching solid tori to the components of $F$ so that the lengths of meridians on each component of $F$ go to infinity as $i \to \infty$.

In a complete nonpositively curved manifold, a null-homotopic curve of length $K$ bounds a disk of diameter at most $K/2$ (use a ruled disk). Let $\alpha$ be a curve in $M$ which is null-homotopic in $M_i$ for infinitely many $i$. Let $K$ denote the length of $\alpha$, and let $H$ denote a union of horoball neighborhoods of the cusps, one for each component of $F$, that exclude a $K/2$-neighborhood of $\alpha$. Then $M_i$ is obtained from $M \setminus H$ by gluing a union of solid tori $H_i$ to the torus boundaries of $H$. Since the meridians of $H_i$ have lengths going to $\infty$ with $i$, it follows from the "$2\pi$-theorem" of Gromov-Thurston (see [37], and Bleiler-Hodgson [13] or Moriah-Rubinstein [72] for a detailed proof) that for large enough $i$ one can put a metric $\sigma_i$ on $M_i$, which has sectional curvatures pinched between some fixed negative constants, and equals the hyperbolic metric $\sigma$ in $M \setminus H$. If $\alpha$ is null-homotopic in $M_i$ it bounds a disk of $\sigma_i$-diameter at most $K/2$, as above, but it follows that this disk is contained in $M \setminus H$, and hence $\alpha$ was already null-homotopic in $M$. This proves the first assertion of the lemma.

For the second assertion, we consider a manifold with sectional curvatures bounded above by some $-\kappa < 0$. Let $A : S^1 \times [0,1] \to N$ be a *ruled annulus*, that is, a map such that $A|_{\{\theta\} \times [0,1]}$ is a geodesic for each $\theta \in S^1$. Then the induced metric on the annulus also has sectional curvatures bounded by $-\kappa$. Standard comparison theorems imply that $\mathrm{Area}(A) \leq C(\kappa)\ell(\partial A)$. Furthermore, given $K, \varepsilon > 0$ there exists $D(\kappa, \varepsilon, K)$ such that, if $\partial A$ has length at most $K$ and $A$ has diameter at least $D$, then the shortest curve in $A$ homotopic to a boundry component has length at most $\varepsilon$.

Let $-\kappa$ be the upper bound on the sectional curvatures of the metrics $\sigma_i$ (in fact $\kappa$ can be taken arbitrarily close to 1 as $i \to \infty$). Let $\varepsilon$ be less than the $\sigma$-length of the shortest noncuspidal curve in $M$. Let $\alpha$ be a curve whose image in $\pi_1(M_i)$ is an $m^{\text{th}}$ power for infinitely many $i$, and which is not homotopic into a cusp. Let $K$ be the $\sigma$-length of $\alpha$, and choose the horoball neighborhoods $H$ so that they avoid a $D(\kappa, \varepsilon, 2K)$-neighborhood of $\alpha$.

If $\alpha_i^*$ is the $\sigma_i$-geodesic representative of $\alpha$ in $M_i$, let $A$ be a ruled annulus with one boundary mapping to $\alpha$ and the other to $\alpha_i^*$. After a small perturbation we may assume $A$ is transverse to $\partial H$, and consider the intersection locus $A^{-1}(\partial H)$ in the annulus. If a component is homotopically trivial in the annulus and its image is a power of a meridian of $H_i$, then the area enclosed by this curve is at least the area of a meridian disk of $H_i$, and this goes to $\infty$ with $i$ (see Moriah-Rubinstein [72]). It follows from the upper bound on $\text{Area}(A)$ that eventually every such component bounds a disk that can be pushed out of $H_i$.

Furthermore if $A$ meets $H$ at all then its diameter is at least $D$, and hence it contains a nontrivial curve of length at most $\varepsilon$. This curve must then be cuspidal, so that $A$ intersects $\partial H$ in a nontrivial curve. Combining this with the previous observation we conclude that $\alpha$ can be deformed to $\partial H$ within $M \setminus H$, a contradiction to the choice of $\alpha$.

It follows that $A$ misses $H$, so that $\alpha$ is homotopic to $\alpha_i^*$ within $M \setminus H$. This means that $\alpha_i^*$ is already the geodesic representative of $\alpha$ in the metric $\sigma$ on $M$. Since $\alpha$ is an $m^{\text{th}}$ power in $M_i$, $\alpha_i^*$ runs $m$ times around the geodesic representing the root, and therefore is an $m^{\text{th}}$ power in $M$ as well. This proves the second assertion.

To apply this to our situation, note that the condition that $\gamma_i$ and $\gamma_{i+1}$ intersect implies that $S \times [s,t]$ minus the curves $\gamma_i \times \{n_i\}$ is acylindrical and hence hyperbolic, by Thurston's geometrization theorem (actually more simply in the case of the punctured torus, one can do it using Klein-Maskit combination theorems, combining a sequence of punctured torus groups with accidental parabolics across totally geodesic thrice-punctured sphere boundaries).

11.5. *The exterior of the convex core.* There is one last detail to tidy up: We have a proper degree 1 map $f : M \to \widehat{C}(N)$, which in the case of two degenerate ends can be written as $f : M \to N$ and lifts to $\widetilde{f} : \widetilde{M} \to \mathbf{H}^3$. Then there is a geometrically finite end, $\widehat{C}(N) \neq N$, and we would like to get a map which covers all of $N$ in this case too.

Recall that $\widehat{C}(N)$ is the union of an $r$-neighborhood of the convex core $C(N)$, and the Margulis tubes of short noninternal curves, if any. There is a submanifold $C(M)$ in $M$ which corresponds to $C(N)$, defined as follows: suppose $\mathcal{B}$ is a boundary block with solid torus $U$. Then $\partial U$ is composed of four annuli $A_0, A_1, A_L$ and $A_R$ as in Section 10 where $A_1$ or $A_0$ (or both) lie in $\partial M$ (in the case of a block with parabolic, one of $A_0$ or $A_1$ is actually two annuli). Thus $Y = \partial U \setminus \partial M$ is a union of one or two annuli. Recalling that $U$ has been given the metric of a Margulis tube, let $C(Y)$ denote the convex hull of $Y$ in $U$ with respect to this metric. Let $C(M)$ be the submanifold of $M$ obtained by replacing $U$ with $C(Y)$ for each of the (0, 1 or 2) boundary blocks.

We can now append to $C(M)$ a collar of the form $\partial C(M) \times [0, \infty)$ with metric

$$(11.2) \qquad ds^2 \cosh^2 t + dt^2,$$

where $ds^2$ is the metric on $\partial C(M)$ and $t \in [0, \infty)$. Let $\overline{M}$ denote $C(M)$ union this collar. We extend $f$ to this collar to obtain $\overline{f} : \overline{M} \to N$ via the following observations:

For each boundary block as above, we note that in $N$, the boundary of the Margulis tube $\mathbf{T}$ corresponding to $U$ meets $C(N)$ in one or two annuli, whose $r$-neighborhood is the image $f(Y)$. Since, in the universal cover, the lifts of $\mathbf{T}$ and $C(N)$ are both convex, it follows that $C(N) \cap \mathbf{T}$ is the convex hull in $\mathbf{T}$ of $C(N) \cap \partial \mathbf{T}$.

One can now apply the methods of Epstein-Marden [30] to see that $f$ restricted to $\partial C(Y)$ is homotopic rel boundary to a bilipschitz map from $\partial C(Y)$ to $\mathbf{T} \cap \partial C_r(N)$. (The constants can in fact be made independent of the group, but this does not matter to us).

Epstein-Marden also show in [30] that the exterior $N \setminus C(N)$ is bilipschitz equivalent to a metric of the form (11.2), where $ds^2$ is the metric on $\partial C(N)$. It therefore follows, after an appropriate interpolation in a collar $\partial C(M) \times [0, 1]$, that $\overline{f}$ can be made a proper degree 1 map from $\overline{M}$ to $N$, which is bilipschitz in $\partial C(M) \times [1, \infty)$, and lifts to a quasi-isometry $\widetilde{\overline{f}} : \widetilde{\overline{M}} \to \mathbf{H}^3$.

## 12. Proofs of the main theorems

12.1. *The ending lamination theorem.* We are now ready to prove Theorem A. We may proceed as in [70]. If $N_1$ and $N_2$ have the same ending invariants we obtain liftable quasi-isometries $f_1, f_2$ from the same model manifold $\overline{M}$ to $N_1$ and $N_2$. Thus the map $\widetilde{f}_1 \circ \widetilde{f}_2^{-1}$, where $\widetilde{f}_2^{-1}$ is any quasi-inverse of $\widetilde{f}_2$, gives a quasi-isometric conjugacy of the two group actions in $\mathbf{H}^3$. This extends at infinity to a quasi-conformal conjugacy of the two group actions on the sphere (see Mostow [73], and in a more general context [34], [27]). Furthermore, if the domain of discontinuity is nonempty the equality of ending invariants implies that the map can be made conformal from $\Omega_1$ to $\Omega_2$ (on the quotient surfaces there is a bounded homotopy to a conformal map, and this can be lifted to a homotopy which is constant on the limit set). Sullivan's theorem [83] states that a quasiconformal conjugacy of finitely generated Kleinian groups which is conformal on the domains of discontinuity is in fact an element of $\mathrm{PSL}_2(\mathbf{C})$; so we are done.

12.2 *The rigidity theorem.* In this section we prove Theorem C. Consider two punctured torus groups $\Gamma_1, \Gamma_2$ whose actions on the sphere are topologically conjugate: there exists a homeomorphism $\psi : \hat{\mathbf{C}} \to \hat{\mathbf{C}}$ such that $\Gamma_2 = \psi \Gamma_1 \psi^{-1}$ in $\mathrm{homeo}(\hat{\mathbf{C}})$. We wish to prove that $\psi$ can be replaced by a quasiconformal or anti-quasiconformal map.

Note first that it suffices to consider orientation-preserving $\psi$ (if not, let $\psi' = R \circ \psi$ where $R$ is a Möbius inversion and apply the argument to $\Gamma_2' = R\Gamma_2 R^{-1}$). Thus we would like to show that $\Gamma_2 = \Phi \Gamma_1 \Phi^{-1}$ where $\Phi$ is quasiconformal.

The conjugacy gives an identification of $\Gamma_1$ and $\Gamma_2$ as groups, so we may view them both as representations $\rho_1, \rho_2$ of a fixed copy of $\pi_1(S)$. Let $(\nu_-^i, \nu_+^i)$ denote the ordered pair of ending invariants of $\rho_i$ where $i = 1, 2$.

Let $\alpha_j$ be a sequence of elements of $\mathcal{C}$, with geodesic representatives $\alpha_j^i$ in $N_i = \mathbf{H}^3/\Gamma_i$ (i=1,2). We can characterize whether the $\alpha_j^i$ leave every compact set in $N_i$ as follows: Let $F_i(\alpha_j)$ be the set of pairs of fixed points $(x, y) \in \hat{\mathbf{C}}$ of the elements in the conjugacy class in $\Gamma_i$ determined by $\alpha_j$. Then $\alpha_j^i$ leaves every compact set if and only if $d_j^i = \sup\{d(x, y) : (x, y) \in F_i(\alpha_j)\}$ go to zero as $j \to \infty$ (where $d(\cdot, \cdot)$ denotes spherical distance). This property is preserved by a homeomorphism of the sphere (compare Ohshika [75]). Thus a sequence converges to an ending lamination of $\Gamma_1$ if and only if it converges to an ending lamination of $\Gamma_2$.

We conclude that any irrational ending invariant of $\Gamma_1$ appears as an invariant of $\Gamma_2$. The same holds for a rational invariant by an easier argument, since $\psi$ must conjugate parabolics to parabolics. Suppose, say, $\nu_+^1 \in \hat{\mathbf{R}}$. If $\nu_+^1 = \nu_+^2$ (in which case, if $\nu_-^1 \in \hat{\mathbf{R}}$ then $\nu_-^1 = \nu_-^2$), we can apply the ending lamination theorem directly to conclude that $\Gamma_1$ and $\Gamma_2$ are quasiconformally conjugate, and we are done.

Suppose then that $\nu_+^1 = \nu_-^2$ (necessarily then $\nu_-^1 = \nu_+^2$ if $\nu_-^1 \in \hat{\mathbf{R}}$), and let us derive a contradiction. If $R$ is a Möbius inversion then the conjugate representation $\rho_2'$ defined by $\rho_2'(g) = R\rho_2(g)R^{-1}$ has reversed ending invariants $(\nu_+^2, \nu_-^2)$. (See Section 3 to see how orientation induces an ordering on $\nu_\pm$.) Applying the ending lamination theorem to $\rho_1$ and $\rho_2'$, we get a qc-homeomorphism $F$ conjugating $\Gamma_1$ to $\Gamma_2' = \rho_2'(\pi_1(S))$ and inducing the same identification as $R \circ \psi$. Thus the composition $\Phi = F^{-1} \circ R \circ \psi$ is an *orientation-reversing* homeomorphism conjugating $\Gamma_1$ to itself, and acting as the identity on the group.

This can of course happen if $\Gamma_1$ is quasi-Fuchsian and $\Phi$ exchanges the components of $\Omega_{\Gamma_1}$, but we claim now that this is the only case. Since the conjugation is the identity on the group, $\Phi$ must restrict to the identity on $\Lambda_{\Gamma_1}$, which is the closure of the set of fixed points of $\Gamma_1$. Thus if $\Lambda_{\Gamma_1} = \hat{\mathbf{C}}$ we have $\Phi = \mathrm{id}$ and it cannot be orientation-reversing.

Suppose $\Omega_{\Gamma_1}$ is nonempty. For any component $O$ of $\Omega_{\Gamma_1}$, $\Phi$ is the identity on $\partial O$ so either $\Phi(O) = O$ or $\Phi(O) = \hat{\mathbf{C}} - \overline{O}$. The latter is impossible when $\Gamma_1$ is not quasi-Fuchsian since $\hat{\mathbf{C}} - \overline{O}$ is empty or contains some of the limit set. The former implies that $\Phi$ descends to an orientation-reversing homeomorphism of the quotient surface $O/\Gamma_1$ which induces the identity on its fundamental group, and this is also impossible. This contradiction implies that, in fact, the ending laminations must have matched up correctly to begin with.

Note that, in case both $\nu_{\pm}^1$ are in $\mathbf{D}$ or in $\hat{\mathbf{Q}}$, the same holds for $\nu_{\pm}^2$, and in this case (the geometrically finite case) the theorem follows from Marden's isomorphism theorem [58].

12.3. *Topology of the deformation space.* In this section we give the proof of Theorem B. Let $\mathcal{R} = \mathcal{R}(\pi_1(S))$ denote the space of representations of $\pi_1(S)$ into $\mathrm{PSL}_2(\mathbf{C})$ with parabolic commutator, modulo conjugation of the image in $\mathrm{PSL}_2(\mathbf{C})$. In the natural topology on $\mathcal{R}$ (known as the "algebraic topology"), the subset $\mathcal{D}$ of discrete, faithful representations is closed (by a theorem of Chuckrow, see [49]) and has nonempty interior (precisely the quasi-Fuchsian representations, by theorems of Marden [58] and Sullivan [84]).

Let $\Delta$ denote the diagonal of $S^1 \times S^1$. The ending invariant construction gives a map

$$\nu : \mathcal{D} \to (\overline{\mathbf{D}} \times \overline{\mathbf{D}}) \setminus \Delta$$

associating to a representation $\rho$ its ordered pair of ending invariants. Points in $\Delta$ are never obtained: the same element cannot be parabolic on both ends of the manifold, and more generally the two ending laminations must be distinct (see Thurston [86] or Bonahon [14]).

It is perhaps surprising that the map $\nu$ is *not continuous*. As pointed out to the author by Dick Canary, one can construct sequences of representations $\rho_i$ which converge to some $\rho$, but for which the invariants $(\nu_-^i, \nu_+^i)$ converge to the diagonal $\Delta$ (for example take some $\alpha \in \mathcal{C}$ and consider the sequence $(D_\alpha^i \nu_0, D_\alpha^{2i} \nu_0)$). The construction, which is related to the Kerckhoff-Thurston example in [55], appears in Anderson-Canary [6]. See also McMullen [66] for a further discussion of this phenomenon.

Nevertheless the following fact does hold:

LEMMA 12.1. *Let $\{\rho_i\}$ be a sequence of marked punctured-torus groups, whose ending invariants $(\nu_-^i, \nu_+^i) = \nu(\rho_i)$ converge to a pair $(\nu_-^\infty, \nu_+^\infty)$ which is not in $\Delta$. Then the $\rho_i$ converge in $\mathcal{D}$ to a representation $\rho$, and its ending invariants are $\nu(\rho) = (\nu_-^\infty, \nu_+^\infty)$.*

*Remark.* The conclusion that a convergent subsequence exists (up to conjugation in $\mathrm{PSL}_2(\mathbf{C})$) is a special case of Thurston's Double Limit Theorem (see Thurston [85] or Otal [77]), though it will follow from our analysis as well. Our real interest here is in the fact that the limiting end invariants have the expected values. □

*Proof.* Suppose first that $\nu_-^\infty, \nu_+^\infty \in \mathbf{D}$. Then eventually $\nu_\pm^i \in \mathbf{D}$ as well, and we have convergence in the space of quasi-Fuchsian groups. In this case the result is already known, by the Ahlfors-Bers theory [10], [12].

More generally, suppose that one invariant, say $\nu_+^\infty$, is in $\mathbf{D}$. This is essentially the case of convergence in a Bers slice. Let $\sigma$ be a metric representing $\nu_+^\infty$ on $S$, and $f_i : (S, \sigma) \to \Omega_+(\rho_i)/\rho_i(\pi_1(S))$ be the extremal quasiconformal map inducing $\rho_i$ on the fundamental group. Since $\nu_+^i \to \nu_+^\infty$, the dilatation $K_i$ of $f_i$ converges to 1. Thus the lifts $\tilde{f}_i : \tilde{S} \to \Omega_+(\rho_i)$, after appropriate normalization, converge on compact sets, by general compactness theorems for quasiconformal maps (see Lehto-Virtanen [57]), to a conformal embedding $f_\infty : \tilde{S} \to \hat{\mathbf{C}}$, and this map induces a limit representation $\rho$ on the fundamental groups. It follows that $f_\infty(\tilde{S}) = \Omega_+(\rho)$, and hence the top invariant of $\rho$ is $\nu_+$.

We now consider the limiting behavior of invariants that escape the interior. Following the discussion in Section 4, we have for each $i = 1, \ldots, \infty$ a pair of points $\alpha_\pm^i \in \partial \mathbf{D}$, an edge set $E^i = E(\alpha_-^i, \alpha_+^i)$, and a pivot sequence with internal pivot subsequence $P^i$ and $P_0^i$, respectively.

Let us assume for the moment that, if one of $\nu_\pm^\infty$ lies in the interior then it is generic in the sense that $\alpha_\pm^\infty$ is uniquely determined. Hence the same will eventually be true for $\nu_\pm^i$. We return to this point at the end.

We claim that every pivot in $P^\infty$ is eventually in $P^i$. Consider an edge $e \in E^\infty$. It separates $\alpha_-^\infty$ from $\alpha_+^\infty$. If $\nu_+^\infty \in \partial \mathbf{D}$, we have $\nu_+^\infty = \alpha_+^\infty$ and hence $\nu_+^i \to \alpha_+^\infty$. It then follows that $\alpha_+^i \to \alpha_+^\infty$. If $\nu_+^\infty$ lies in the interior then (by the assumption of genericity) eventually $\nu_+^i$ is sufficiently close that $\alpha_+^i = \alpha_+^\infty$. The same argument applies to $\alpha_-^\infty$.

Thus in any case we find that eventually $e$ separates $\alpha_-^i$ from $\alpha_+^i$, and hence $e \in E^i$. It follows that every internal pivot of $P^\infty$ eventually lies in $P^i$.

Suppose that $\alpha_+^\infty$ is in $\mathcal{C}$, and hence is the last, noninternal pivot of $P^\infty$. If $\nu_+^\infty \in \mathbf{D}$ then as above, eventually $\alpha_+^i = \alpha_+^\infty$, and in particular $\alpha_+^\infty$ is the last vertex of $P^i$. If $\nu_+^\infty = \alpha_+^\infty$ something interesting can happen: Normalize $\mathbf{D}$ as $\mathbf{H}^2$ where $\alpha_+^\infty = \infty$. In this normalization, $|\nu_+^i| \to \infty$. If $\mathrm{Im}\,\nu_+^i \to \infty$ then eventually $\alpha_+^i = \infty = \alpha_+^\infty$. If not, then $|\mathrm{Re}\,\nu_+^i| \to \infty$, and since $\nu_-^i$ are converging to $\nu_-^\infty \neq \infty$, eventually an arbitrary number of vertical Farey edges separate $\alpha_-^i$ from $\alpha_+^i$, so that $\alpha_+^\infty$ is an *internal* vertex of $P^i$.

The same argument applies to $\alpha_-^\infty$, so our claim is proved.

Suppose now that $P^\infty$ contains at least two pivots (equivalently that $\alpha_-^\infty \neq \alpha_+^\infty$). Let us number them $\alpha_0$ and $\alpha_1$, and let $\mathcal{B}_0^i, \mathcal{B}_1^i$ be the corresponding blocks in the model $M_i$ for large enough $i$. The surface $S^i(1) = \partial_0 \mathcal{B}_1^i = \partial_1 \mathcal{B}_0^i$ has a metric where both $\alpha_0$ and $\alpha_1$ have shortest length, so we may identify these surfaces, by the identity on the underlying $S$, for all $i$. The map $f_i|_{S(1)}$ induces the representation $\rho_i$. Thus, after conjugating the image groups appropriately, or equivalently after choosing a lift of $f_i|_{S(1)}$ to $\mathbf{H}^3$ that maps a fixed basepoint to the origin, we find that the generators of the group associated to $\alpha_0$ and $\alpha_1$ have bounded translation lengths, and hence (by Chuckrow's theorem) there is a convergent subsequence with some limit representation $\rho_\infty$. Let us reindex this as if it were the whole sequence.

Taking a further subsequence, we may also assume by Theorem 11.2 that there is a geometric limit $f_\infty : M_\infty \to N_\infty$, so that $S(1)$ still embeds in $M_\infty$ between the first two blocks, and $f_\infty|_{S(1)}$ induces the representation $\rho_\infty$. We now wish to determine the ending invariants of $\rho_\infty$.

The width sequence for $P^\infty$ is, by the above discussion, exactly reproduced as a subset of the width sequence for $M_\infty$, and the last pivot $\alpha_+$ of $P^\infty$ (if there is one) appears as the first (lowest-numbered) cusp in $M_\infty$ above $S(1)$. Similarly the first pivot $\alpha_-$ of $P^\infty$, if there is one, appears as the last cusp below $S(1)$.

Lift $M_\infty$ to the cover $\widetilde{M}_\infty$ corresponding to $S(1)$, and lift $f_\infty$ to a map from $\widetilde{M}_\infty$ to the corresponding cover $\widetilde{N}_\infty$ of $N_\infty$. The map induced from $\pi_1(S) = \pi_1(\widetilde{M}_\infty)$ to $\pi_1(\widetilde{N}_\infty)$ is just the limit representation $\rho_\infty$. If there is a cusp of $M_\infty$ above $S(1)$ corresponding to $\alpha_+$, it lifts to a rank 1 cusp $\widetilde{T}$ in the top half of $\widetilde{M}_\infty$, where we use our orientation convention to order the unbounded components of $\widetilde{M}_\infty \setminus S(1)$. Since $\widetilde{f}_\infty$ is now a degree 1 proper map between manifolds homeomorphic to $S \times \mathbf{R}$, it follows that $\widetilde{f}_\infty(\widetilde{T})$ is in the unbounded component of $\widetilde{N}_\infty \setminus \widetilde{f}_\infty(S(1))$ which is a neighborhood of its $e_+$ end. Thus the $e_+$ ending invariant of $\rho_\infty$ is exactly $\alpha_+^\infty = \nu_+^\infty$ in this case.

If in fact $\nu_+^\infty$ is irrational then $P^\infty$ is infinite in the forward direction, and so is $M_\infty$. The pivots in $M^\infty$ therefore converge (in the normalization determined by $\alpha_0$ and $\alpha_1$) to $\nu_+^\infty$. In this case the portion of $M_\infty$ above $S(1)$ lifts homeomorphically to $\widetilde{M}_\infty$, and maps to the $e_+$ end of $\widetilde{N}_\infty$. Thus the $e_+$ invariant of $\rho_\infty$ is again $\nu_+^\infty$.

The same argument applies to $e_-$ when $\nu_-^\infty$ is a boundary point.

We have shown that some *subsequence* of the original $\{\rho_i\}$ converges, and that for every limit $\rho$ that can arise, the invariants $\nu(\rho)$ are exactly $(\nu_\pm^\infty)$. The Ending Lamination Theorem (Theorem A) now implies that all of these limits are equal, and hence in fact the entire sequence converges. (We remark, however, that the set of possible *geometric* limits can be quite large, as we have seen.)

We have left out the possibility that $\alpha_-^\infty = \alpha_+^\infty$. Since $\nu_+^\infty$ and $\nu_-^\infty$ are different by assumption, this can only occur if at least one of them, say $\nu_+^\infty$, is in the interior. By the Bers-type argument in the beginning of the proof we therefore have a convergent subsequence, and the rest of the argument goes through in the same way.

Finally, let us return to the minor assumption that $\nu_+^\infty$ (similarly $\nu_-^\infty$) is generic, in the sense that if it lies in $\mathbf{D}$ it is closest to a unique $\alpha_+^\infty \in \mathcal{C}$. If this should be false then there are two or three possible choices, and any one can be made. However in the sequence $\nu_+^i \to \nu_+$, perhaps for each $i$ a unique $\alpha_+^i$ can be chosen, and not the same for each $i$. This is easily resolved by separating the sequence into a finite number of subsequences in each of which the $\alpha_\pm^i$ are eventually constant. The limits for each subsequence will agree, as above. □

The fact that the map $\nu$ is *surjective* is already known: Bers' Simultaneous Uniformization theorem implies that every pair in $\mathbf{D} \times \mathbf{D}$ is obtained and in fact that $\nu$ is a homeomorphism from the quasi-Fuchsian groups to $\mathbf{D} \times \mathbf{D}$. The boundary points of $\overline{\mathbf{D}} \times \overline{\mathbf{D}} \setminus \Delta$ can be reached using work of Maskit [61] and Keen-Maskit-Series [54] for rational points and Thurston [85], [77] for irrational points. See also Ohshika [74] for generalizations. Our Lemma 12.1, which reproduces Thurston's theorem in our setting, also accomplishes this.

The ending lamination theorem serves to show that $\nu$ is *injective*, and hence
$$\nu^{-1} : (\overline{\mathbf{D}} \times \overline{\mathbf{D}}) \setminus \Delta \to \mathcal{D}(\pi_1(S))$$
is well-defined. Although $\nu$ is not a homeomorphism (it is discontinuous, as above), Lemma 12.1 implies exactly that $\nu^{-1}$ is continuous.

This gives the first statement of Theorem B, and the fact that every representation in $\mathcal{D}$ is the limit of quasi-Fuchsian representations. Now let us consider Bers and Maskit slices.

A Bers slice, in our terminology, is the image by $\nu^{-1}$ of a slice $\{x_0\} \times \overline{\mathbf{D}}$, where $x_0 \in \mathbf{D}$. Bers defined these slices for a general surface $S$ (see Gardiner [33]), as a subset of $\mathcal{D}(\pi_1(S))$ where the "bottom" invariant component of $\Omega$ has quotient equal to a fixed point in the Teichmüller space of $S$. He showed that such a slice is compact and conjectured that it was the closure of its interior, identified with the Teichmüller space of $S$.

In our case, since $\nu^{-1}$ is continuous and $\{x_0\} \times \overline{\mathbf{D}}$ is compact, the image must be compact and the restriction of $\nu^{-1}$ is a homeomorphism. Hence Bers' conjecture holds for the punctured torus.

A Maskit slice is a similar construction, where $x_0$ is taken to be a rational point on the boundary. In this case the point $(x_0, x_0)$ must be omitted so the slice is $\nu^{-1}(\{x_0\} \times (\overline{\mathbf{D}} \setminus \{x_0\}))$. (Again, the original definition of Maskit applies for the higher genus case as well, and does not depend on the existence of the irrational end invariants.)

Since the Maskit slice is not compact, we must be more careful in showing that $\nu^{-1}$ restricts to a homeomorphism. It suffices to show that $\nu^{-1}|_{\{x_0\} \times (\overline{\mathbf{D}} \setminus \{x_0\})}$ is *proper*, or in other words given a sequence $(x_0, y_i)$ with $y_i \to x_0$, show that the corresponding representations $\rho_i$ must diverge. Consider again our pivot sequence for $(x_0, y_i)$.

Recall that $y_i(x_0)$ is the image of $y_i$ in a fixed normalization of $\mathbf{H}^2$ sending $x_0$ to $\infty$. Suppose first that $\operatorname{Im} y_i(x_0) \to \infty$. Then $x_0$ is represented by a very short curve on both sides of the convex core for each $\rho_i$. An annulus between these two curves therefore cuts the convex core into a thrice-punctured sphere times an interval, and every curve in $S$ intersecting $x_0$ must have geodesic representative crossing this annulus, so its length will go to infinity as $i \to \infty$. It follows that the translation lengths diverge for at least one of any generator pair, so the sequence of representations diverges. This is in fact the exceptional geometric limit we have mentioned before, in which there is one double-boundary block in the model and the geometric limit is a thrice-punctured sphere group.

The remaining case is that $\operatorname{Im} y_i(x_0)$ remain bounded and hence $\operatorname{Re} y_i(x_0) \to \pm\infty$.

In this case we can apply the pivot sequence analysis to show that there is a sequence of integers $\beta_i \to \pm\infty$ representing neighbors of $x_0$ in $\mathcal{C}$, with $\ell_{\rho_i}(\beta_i)$ bounded. Namely, if $y_i$ is in a triangle $T_i$ with vertex $x_0$, take $\beta_i$ to be one of the other vertices. Otherwise take $\beta_i$ to be the first internal pivot in the pivot sequence for $\rho_i$. (More concretely, $\beta_i$ is the integer part of $\operatorname{Re} y_i(x_0)$, plus or minus 1).

Now letting $\operatorname{tr}_n^i$ denote the trace in the representation $\rho_i$ of the neighbor of $x_0$ indexed by $n$, we may apply the trace identity (7.2) together with the parabolic commutator condition, as in Section 7, to conclude that $\operatorname{tr}_n^i = \operatorname{tr}_0^i + cn$ where $c = \pm 2\sqrt{-1}$. Since $\operatorname{tr}_{\beta_i}^i$ stays bounded we conclude that $\operatorname{tr}_0^i \to \infty$. Thus again there is no fixed choice of generators that converge.

This proves that the Maskit slice is a homeomorphic image of a disk minus a boundary point.

*Boundaries of embeddings.* Both Maskit and Bers slices have well-known embeddings into finite-dimensional vector spaces. In the case of the punctured torus these have one complex dimension (see Wright [89] and Keen-Series [53] for a thorough development and good pictures of this embedding for a Maskit slice). Our results imply in particular that the boundaries of these embeddings are Jordan curves.

STATE UNIVERSITY OF NEW YORK, STONY BROOK, NY
*E-mail address*: yair@math.sunysb.edu

References

[1] W. Abikoff, Kleinian groups – geometrically finite and geometrically perverse, Geometry of Group Representations, AMS Contemp. Math. no. 74, 1988, 1–50.
[2] L. Ahlfors, An extension of Schwarz's lemma, Trans. Amer. Math. Soc. **43** (1938), 359–364.
[3] ———, *Conformal Invariants*: *Topics in Geometric Function Theory*, McGraw-Hill Book Co., New York, 1973.
[4] L. Ahlfors and L. Bers, Riemann's mapping theorem for variable metrics, Ann. of Math. **72** (1960), 385–404.
[5] R. C. Alperin, W. Dicks, and J. Porti, The boundary of the Gieseking tree in hyperbolic three-space, Centre de Recerca Matemàtica Preprint núm 330, April 1996.
[6] J. Anderson and R. Canary, Algebraic limits of Kleinian groups which rearrange the pages of a book, Invent. Math. **126** (1996), 205–214.
[7] W. Ballmann, M. Gromov, and V. Schroeder, *Manifolds of Nonpositive Curvature*, Progress in Math. **61**, Birkhäuser Boston, Inc., 1985.
[8] A. F. Beardon, *The Geometry of Discrete Groups*, Grad. Texts in Math. **91**, Springer-Verlag, New York, 1983.
[9] R. Benedetti and C. Petronio, *Lectures on Hyperbolic Geometry*, Universitext, Springer-Verlag, Berlin, 1992.
[10] L. Bers, Simultaneous uniformization, Bull. Amer. Math. Soc. **66** (1960), 94–97.
[11] ———, On boundaries of Teichmüller spaces and on Kleinian groups I, Ann. of Math. **91** (1970), 570–600.
[12] ———, Spaces of Kleinian groups, Maryland conference in Several Complex Variables I, LNM **155**, Springer-Verlag, New York, 1970, 9–34.
[13] S. Bleiler and C. Hodgson, Spherical space forms and Dehn filling, Topology **35** (1996), 809–833.
[14] F. Bonahon, Bouts des variétés hyperboliques de dimension 3, Ann. of Math. **124** (1986), 71–158.
[15] F. Bonahon and J. P. Otal, Variétés hyperboliques à géodésiques arbitrairement courtes, Bull. London Math. Soc. **20** (1988), 255–261.
[16] B. Bowditch, Markoff triples and quasi-Fuchsian groups, Proc. London Math. Soc. **77** (1998), 697–736.
[17] B. Bowditch and D. B. A. Epstein, Natural triangulations associated to a surface, Topology **27** (1988), 91–117.
[18] J. Brock, Iteration of mapping classes on a Bers slice, Ph.D. thesis, UC Berkeley, 1997.
[19] R. Brooks and J. P. Matelski, Collars in Kleinian groups, Duke Math. J. **49** (1982), 163–182.
[20] P. Buser, *Geometry and Spectra of Compact Riemann Surfaces*, Progress in Math. **106**, Birkhäuser Boston Inc., 1992.
[21] R. Canary and Y. Minsky, On limits of tame hyperbolic 3-manifolds, J. Differential Geom. **43** (1996), 1–41.
[22] R. D. Canary, Ends of hyperbolic 3-manifolds, J. Amer. Math. Soc. **6** (1993), 1–35.
[23] ———, A covering theorem for hyperbolic 3-manifolds and its applications, Topology **35** (1996), 751–778.
[24] R. D. Canary, D. B. A. Epstein, and P. Green, *Notes on Notes of Thurston*, Analytical and Geometric Aspects of Hyperbolic Space, Cambridge University Press, 1987, London Math. Soc. Lecture Notes Series no. 111, 3–92.
[25] J. Cannon and W. Thurston, Group invariant Peano curves, preprint, 1989.
[26] L. Carleson and T. W. Gamelin, *Complex Dynamics*, Universitext Tracts in Math., Springer-Verlag, New York, 1993.
[27] M. Coornaert, T. Delzant, and A. Papadopoulos, Géométrie et theorie de groupes: Les groups hyperboliques de Gromov, LNM **1441**, Springer-Verlag, New York, 1990.

[28] M. Culler, C. McA. Gordon, J. Luecke, and P. B. Shalen, Dehn surgery on knots, Ann. of Math. **125** (1987), 237–300.
[29] A. Douady and C. J. Earle, Conformally natural extension of homeomorphisms of the circle, Acta Math. **157** (1986), 23–48.
[30] D. B. A. Epstein and A. Marden, Convex hulls in hyperbolic space, a theorem of Sullivan, and measured pleated surfaces, *Analytical and Geometric Aspects of Hyperbolic Space*, London Math. Soc. Lecture Notes Series no. 111, Cambridge University Press, 113–253, 1987.
[31] A. Fathi, F. Laudenbach, and V. Poenaru, *Travaux de Thurston sur les surfaces*, vol. 66–67, Astérisque, 1979.
[32] W. Fenchel, *Elementary Geometry in Hyperbolic Space*, Studies in Math. **11**, Walter de Gruyter & Co., New York, 1989.
[33] F. Gardiner, *Teichmüller Theory and Quadratic Differentials*, John Wiley & Sons, Inc., New York, 1987.
[34] É. Ghys and P. de la Harpe, *Sur les Groupes Hyperboliques d'après Mikhael Gromov*, Progress in Math. **83**, Birkhäuser Boston, Inc., 1990.
[35] C. McA. Gordon and R. A. Litherland, Incompressible planar surfaces in 3-manifolds, Topology Appl. **18** (1984), 121–144.
[36] C. McA. Gordon and J. Luecke, Knots are determined by their complements, Bull. Amer. Math. Soc. **20** (1989), 83–87.
[37] M. Gromov and W. Thurston, Pinching constants for hyperbolic manifolds, Invent. Math. **89** (1987), 1–12.
[38] J. Harer, Stability of the homology of the mapping class groups of orientable surfaces, Ann. of Math. **121** (1985), 215–249.
[39] ———, The virtual cohomological dimension of the mapping class group of an orientable surface, Invent. Math. **84** (1986), 157–176.
[40] W. J. Harvey, Boundary structure of the modular group, *Riemann Surfaces and Related Topics*: *Proc. of the* 1978 *Stony Brook Conference* (I. Kra and B. Maskit, eds.), Annals of Math. Studies, no. 97, Princeton University Press, Princeton, NJ, 1981.
[41] A. E. Hatcher and W. P. Thurston, A presentation for the mapping class group of a closed orientable surface, Topology **19** (1980), 221–237.
[42] N. V. Ivanov, Automorphisms of complexes of curves and of Teichmüller spaces, Internat. Math. Res. Notices No. 14, (1997), 651–666, .
[43] ———, Complexes of curves and the Teichmüller modular group, Uspekhi Mat. Nauk **42** (1987), 55–107.
[44] ———, Complexes of curves and Teichmüller spaces, Math. Notes **49** (1991), 479–484.
[45] W. H. Jaco and P. B. Shalen, *Seifert Fibered Spaces in* 3-*manifolds*, Mem. of the Amer. Math. Soc. **21**, no. 220, A.M.S., 1979.
[46] T. Jørgensen, On pairs of once-punctured tori, Unpublished manuscript.
[47] ———, On cyclic groups of Möbius tranformations, Math. Scand. **33** (1973), 250–260.
[48] ———, On discrete groups of Möbius transformations, Amer. J. of Math. **98** (1976), 739–749.
[49] T. Jørgensen and P. Klein, Algebraic convergence of finitely generated Kleinian groups, Quarterly J. of Math. Oxford **33** (1982), 325–332.
[50] T. Jørgensen and A. Marden, Two doubly degenerate groups, Quarterly J. Math. Oxford **30** (1979), 143–156.
[51] ———, Algebraic and geometric convergence of Kleinian groups, Math. Scand. **66** (1990), 47–72.
[52] D. Kazhdan and G. Margulis, A proof of Selberg's conjecture, Math. USSR Sb. **4** (1968), 147–152.
[53] L. Keen and C. Series, Pleating coordinates for the Maskit slice embedding of the Teichmüller space of punctured tori, Topology **32** (1993), 719–749.

[54] L. Keen, B. Maskit, and C. Series, Geometric finiteness and uniqueness for Kleinian groups with circle packing limit sets, J. Reine Angew. Math. **436** (1993), 209–219.
[55] S. Kerckhoff and W. P. Thurston, Noncontinuity of the action of the modular group at Bers' boundary of Teichmüller space, Invent. Math. **100** (1990), 25–47.
[56] C. Kourouniotis, The geometry of bending quasi-Fuchsian groups, *Discrete Groups and Geometry* (W. J. Harvey and C. Maclachlan, eds.), London Math. Soc. Lecture Notes no. 173, 148–164, Cambridge University Press, 1992.
[57] O. Lehto and K. I. Virtanen, *Quasiconformal Mapping*, Springer-Verlag, New York, 1965.
[58] A. Marden, The geometry of finitely generated Kleinian groups, Ann. of Math. **99** (1974), 383–462.
[59] B. Maskit, On boundaries of Teichmüller spaces and on Kleinian groups II, Ann. of Math. **91** (1970), 607–639.
[60] ______, Moduli of marked Riemann surfaces, Bull. Amer. Math. Soc. **80** (1974), 773–777.
[61] ______, Parabolic elements in Kleinian groups, Ann. of Math. **117** (1983), 659–668.
[62] ______, Comparison of hyperbolic and extremal lengths, Ann. Acad. Sci. Fenn. **10** (1985), 381–386.
[63] ______, The canonical splitting of a Kleinian group, IHES preprint, 1992.
[64] H. A. Masur and Y. Minsky, Geometry of the complex of curves I: Hyperbolicity, Stony Brook IMS Preprint #1996/11, Invent. Math., to appear.
[65] D. McCullough, Compact submanifolds of 3-manifolds with boundary, Quart. J. Math. Oxford **37** (1986), 299–307.
[66] C. McMullen, Complex earthquakes and Teichmüller theory, J. Amer. Math. Soc. **11** (1998), 283–320.
[67] ______, Cusps are dense, Ann. of Math. **133** (1991), 217–247.
[68] R. Meyerhoff, A lower bound for the volume of hyperbolic 3-manifolds, Canad. J. Math. **39** (1987), 1038–1056.
[69] Y. Minsky, Teichmüller geodesics and ends of hyperbolic 3-manifolds, Topology **32** (1993), 625–647.
[70] ______, On rigidity, limit sets, and end invariants of hyperbolic 3-manifolds, J. Amer. Math. Soc. **7** (1994), 539–588.
[71] ______, On Thurston's ending lamination conjecture, Conf. Proc. Lecture Notes Geom. Topol. IV, Internat. Press, Cambridge, MA, 1994.
[72] Y. Moriah and H. Rubinstein, Heegaard structures of negatively curved 3-manifolds, Comm. Anal. Geom. **5** (1997), 375–412.
[73] G. D. Mostow, *Strong Rigidity of Locally Symmetric Spaces*, Annals of Math. Studies no. 78, Princeton University Press, Princeton, NJ, 1973.
[74] K. Ohshika, Ending laminations and boundaries for deformation spaces of Kleinian groups, J. London Math. Soc. **42** (1990), 111–121.
[75] ______, Topologically conjugate Kleinian groups, Proc. Amer. Math. Soc. **124** (1996), 739–743.
[76] J.-P. Otal, Sur le nouage des géodésiques dans les variétés hyperboliques, C. R. Acad. Sci. Paris Sér. I Math. **320** (1995), 847–852.
[77] ______, *Le théoréme d'hyperbolisation pour les variétés fibrées de dimension* 3, Astérisque, No. 235, 1996.
[78] J. Parker and C. Series, Bending formulae for convex hull boundaries, J. Anal. Math. **67** (1995), 165–198.
[79] R. Penner and J. Harer, *Combinatorics of Train Tracks*, Annals of Math. Studies no. 125, Princeton University Press, Princeton, NJ, 1992.
[80] G. P. Scott, Compact submanifolds of 3-manifolds, J. London Math. Soc. **7** (1973), 246–250.
[81] C. Series, The geometry of Markoff numbers, Math. Intelligencer **7** (1985), 20–29.

[82] C. SERIES, The modular surface and continued fractions, J. London Math. Soc. **31** (1985), 69–80.
[83] D. SULLIVAN, On the ergodic theory at infinity of an arbitrary discrete group of hyperbolic motions, *Riemann Surfaces and Related Topics*: *Proc. of the* 1978 *Stony Brook Conference*, Annals of Math. Studies, no. 97, Princeton University Press, Princeton, NJ, 1981.
[84] ______, Quasiconformal homeomorphisms and dynamics II: Structural stability implies hyperbolicity for Kleinian groups, Acta Math. **155** (1985), 243–260.
[85] W. THURSTON, Hyperbolic structures on 3-manifolds, II: Surface groups and manifolds which fiber over the circle, E-print: math.GT/9801045 at http://front.math.ucdavis.edu.
[86] ______, The geometry and topology of 3-manifolds, Princeton University Lecture Notes, 1982.
[87] ______, Hyperbolic structures on 3-manifolds, I: Deformation of acylindrical manifolds, Ann. of Math. **124** (1986), 203–246.
[88] ______, *Three-Dimensional Geometry and Topology*, (S. Levy, ed.), Princeton Math. Series **35**, Princeton University Press, Princeton, NJ, 1997.
[89] D. WRIGHT, The shape of the boundary of Maskit's embedding of the Teichmüller space of once-punctured tori, preprint, 1990.